\documentclass[reqno,a4paper,12pt]{amsart}
\usepackage{graphicx}

\usepackage[all,poly]{xy}
\usepackage{amsfonts}
\usepackage[mathcal]{eucal}

 \usepackage{mathtools}

\DeclarePairedDelimiter\floor{\lfloor}{\rfloor}
\usepackage{amssymb}
\usepackage{amsmath}
\usepackage{mathrsfs}
\usepackage{color}
\usepackage[colorlinks]{hyperref}
\usepackage{enumitem}

\definecolor{citecol}{RGB}{145, 1, 1}

\hypersetup{colorlinks=true,citecolor=citecol,linkcolor=blue,linktocpage=true}

\theoremstyle{plain}
\newtheorem {lemma}{Lemma}[section] 
\newtheorem {theorem}[lemma]{Theorem}

\newtheorem {thm}[lemma]{Theorem}
\newtheorem {corollary}[lemma]{Corollary}

\newtheorem {cor}[lemma]{Corollary}

\newtheorem {prop}[lemma]{Proposition}

\theoremstyle{definition}
\newtheorem {remark}[lemma]{Remark}
\newtheorem {rem}[lemma]{Remark}

\newtheorem {example}[lemma]{Example}

\theoremstyle{definition}

\newtheorem{deff}[lemma]{Definition}{}
\newtheorem{conj}[lemma]{Conjecture}

\newcommand{\gr}{\operatorname{gr}}

\author{Pere Ara}\address{
Departament de Matemàtiques, Universitat Autònoma de Barcelona, 08193 Bellaterra (Barcelona), Spain.}
\email{pere.ara@uab.cat}
\author{Tran Quang Do}
\address{Institute of Mathematics, VAST, 18 Hoang Quoc Viet, Cau Giay, Hanoi, Vietnam}
\email{tqdo@math.ac.vn}

\author{Tran Giang  Nam}
\address{Institute of Mathematics, VAST, 18 Hoang Quoc Viet, Cau Giay, Hanoi, Vietnam}
\email{tgnam@math.ac.vn}
\title[Dynamics on graphs with disjoint cycles and applications]{Dynamics on graphs with disjoint cycles and applications} 
\begin{document}

\begin{abstract}
In this article,  we introduce the notion of connected finite graphs with disjoint cycles in normal form  and show that any such graph can be transformed into a normal form graph via a finite sequence of in-splittings and out-splittings. Consequently, we provide number-theoretic criteria for meteor graphs of length three to be strongly shift equivalent, where a meteor graph of length three is a connected finite essential graph consisting of three disjoint cycles which makes a unique chain of cycles of length three. We then prove that meteor graphs of length three whose cycle lengths are pairwise coprime are shift equivalent if and only if they are strongly shift equivalent, if and only if their corresponding Leavitt path algebras are graded Morita equivalent, if and only if their graded $K$-theories, $K^{\gr}_0$, are order-preserving $\mathbb{Z}[x, x^{-1}]$-module isomorphic. As a consequence, Williams' Conjecture and Hazrat's Graded Morita Equivalence Conjecture hold for graphs with disjoint cycles that contain exactly three cycles whose lengths are pairwise coprime.
\end{abstract}
\maketitle

\section{Introduction}

Shifts of finite type are fundamental objects of study in symbolic dynamics. They have applications in several areas, including topological quantum field theory, ergodic theory, statistical mechanics, coding theory, and information theory \cite{lindmarcus}.
A basic problem in the theory is to find computable criteria for determining whether two shifts of finite type are conjugate.

Up to conjugacy, every shift of finite type can be realized as the {\it edge shift} of an {\it essential graph}, that is, a finite connected directed graph $E$ with neither sources nor sinks \cite{lindmarcus}. The edge shift $X_E$ associated with $E$ consists of all bi-infinite paths in $E$, equipped with the natural left shift.
Determining whether two shifts of finite type $X_E$ and $X_F$ are conjugate is, in general, a very difficult problem, see for instance \cite[Example 7.3.13]{lindmarcus}. In his seminal paper \cite{williams}, Williams introduced the notions of {\it shift equivalence} (SE) and {\it strong shift equivalence} (SSE), which provide more tractable algebraic criteria for studying conjugacy. He proved that two shifts of finite type $X_E$ and $X_F$ are conjugate if and only if the adjacency matrices of $E$ and $F$ are strongly shift equivalent. Equivalently, this holds if and only if $E$ can be transformed into $F$ by a finite sequence of in-splittings and out-splittings, and their inverses, called {\it in-amalgamations} and {\it out-amalgamations} (see Theorem \ref{willimove} below).

Shift equivalence is a weaker equivalence relation than strong shift equivalence, and is generally more amenable to computation. Williams' Conjecture \cite{williams, willwrong} asserts that shift equivalence and strong shift equivalence are equivalent for shifts of finite type. The conjecture was disproved by Kim and Roush 25 years later by constructing counterexamples \cite{kimroush99}. Despite this negative result, identifying classes of edge shifts for which shift equivalence and strong shift equivalence coincide remains an important open problem. Even for graphs with only two vertices, the situation is subtle, with several examples exhibiting different behaviours; see, for example,  \cite[Example 7.3.13]{lindmarcus}. Given the difficulty of the problem, any result establishing strong shift equivalence under the sole assumption of shift equivalence is therefore of considerable interest.

Beyond their relevance to dynamical systems, the notions of SSE and SE admit interesting interpretations in several areas of operator algebra theory, such as graph $C^*$-algebras, Leavitt path algebras, and groupoids.
Specifically, it is well known (see, e.g., \cite[Corollary 4.7]{CarRout}) that, for essential graphs $E$ and $F$,  $E$ is strongly shift equivalent to $F$ if and only if there is a gauge-invariant diagonal-preserving stable isomorphism between $C^*(E)\otimes \mathcal{K}$ and $C^*(F)\otimes \mathcal{K}$, if and only if there is a diagonal-preserving graded stable isomorphism between $L_{\mathbb{C}}(E)\otimes M_{\infty}(\mathbb{C})$ and $L_{\mathbb{C}}(F)\otimes M_{\infty}(\mathbb{C})$. On the other hand, the shift equivalence of $E$ and $F$ is completely characterized by the graded Grothendieck groups, $K^{\gr}_0(L(E))$ and  $K^{\gr}_0(L(F))$, of the corresponding Leavitt path algebras (\cite{arapar}). Moreover, Hazrat \cite{hazd} showed that if $E$ and $F$ are strongly shift equivalent, then $L(E)$ and $L(F)$ are graded Morita equivalent, and that if the latter
holds, then $E$ and $F$ are shift equivalent. Therefore, graded Morita equivalence of Leavitt path algebras sits right in between strong shift equivalence and shift equivalence. One of the central open problems connecting these notions is Hazrat's Graded Morita Classification Conjecture, which asserts that $E$ and $F$ are shift equivalent if and only if their corresponding Leavitt path algebras are graded Morita equivalent (see \cite{mathann}).
Given these strong connections, considerable effort has recently been devoted (see, e.g., \cite{ART24, BrixCar, toke2, CarEiOrRes, ChenYang, CortHaz, ef, Smith}) to finding operator-algebraic characterizations of the distinction between SE and SSE, as well as of other dynamical properties and aspects of noncommutative geometry.

The goal of this article is to investigate dynamics on graphs with disjoint cycles, motivated in part by the desire to showcase a new setting in which algebraic and number-theoretic tools can be used to address purely dynamical questions. Another motivation for studying this class of graphs is that it completely characterizes Leavitt path algebras of finite Gelfand-Kirillov dimension, as proved by Alahmadi, Alsulami, Jain, and Zelmanov in \cite{aajz:lpaofgkd}. In this case, the Gelfand-Kirillov dimension of the Leavitt path algebra is a nonnegative integer. In particular, for essential graphs, the corresponding Gelfand-Kirillov dimension is always an odd integer. Motivated by the result of Zelmanov and his collaborators, an essential graph is called a {\it graph of finite Gelfand-Kirillov dimension} if its Leavitt path algebra has
finite Gelfand-Kirillov dimension (see \cite{dohaznam}). We also note that the graphs constructed by Kim and Roush, mentioned above, as well as the graphs introduced in \cite[Example 7.3.13]{lindmarcus}, are examples of graphs with infinite Gelfand-Kirillov dimension. A recent preprint by Vas \cite{vas2025} claims that the pointed ordered graded $K_0$-group is a complete invariant, up to graded isomorphism, for Leavitt path algebras of finite graphs with disjoint cycles. However, this pointed classification problem differs from the unpointed graded Morita classification studied here and does not directly yield a Williams-type result. Although Vas's result provides a finer classification for the same class of graphs, the relationship between the pointed graded classification and the unpointed invariants arising from shift equivalence remains a separate question.

Cordeiro, Hazrat, Gillaspy, and Gon\c{c}alves \cite{CGGH} proved that Williams' conjecture holds for the class of \emph{meteor graphs}, namely, essential graphs consisting of two disjoint cycles together with the paths connecting them, which constitute the simplest family of graphs of Gelfand-Kirillov dimension three. This result was subsequently generalized in \cite{dohaznam}, where Hazrat and the second and third authors proved that Williams' conjecture holds for all graphs of Gelfand-Kirillov dimension three.
Hazrat and Pacheco \cite{HP2024} showed that Williams' conjecture holds for the class of essential graphs with three vertices, no parallel edges, and with no non-trivial hereditary and saturated subsets. Motivated by these results, in this article we investigate Williams' conjecture for graphs of higher Gelfand-Kirillov dimension.
To do so, we introduce the notion of {\it connected finite graphs with disjoint cycles in normal form}, which enables a clearer description of the path structure and simplifies the computation of the graded Grothendieck groups $K^{\gr}_0$ of the associated Leavitt path algebras. Furthermore, we prove that every connected finite graph with disjoint cycles can be transformed into a graph in normal form by a finite sequence of in-splittings and out-splittings (Theorem \ref{thm:normal-form}). Consequently, Williams' conjecture for connected finite graphs with disjoint cycles reduces to the corresponding problem for graphs in normal form.

For graphs with Gelfand-Kirillov dimension at least five, the main difficulty is that new paths may arise when applying in-splittings and out-splittings (see Lemma \ref{lem:Trail-inequality}), and no general method is currently known for controlling this phenomenon. As a first step toward the general case, we therefore focus on the simplest class of graphs with Gelfand-Kirillov dimension five, namely, meteor graphs of length three. A \textit{meteor graph of length three}  is a connected essential graph consisting of three disjoint cycles and having a unique chain of cycles of length three. For this class of graphs, we develop an effective method for controlling the growth functions
(Lemmas \ref{lem:Trail-inequality}, \ref{prop-functor} and \ref{functor}). We then establish number-theoretic criteria for two meteor graphs of length three in normal form $E$ and $F$ to be strongly shift equivalent in the case where $E$ can be transformed into $F$ by a finite sequence of in-splittings and out-splittings (Propositions \ref{maintheo-firstcase} and \ref{maintheo:easydiection}). To extend this result to all meteor graphs of length three, it remains to reduce the general case to the setting considered above. To this end, we first characterize a class of meteor graphs of length three that can be obtained from a graph in normal form by a finite sequence of in-splittings and out-splittings. We identify a class of graphs for which this occurs (Lemma \ref{lm:normal-cover}), and, somewhat surprisingly, the lengths of paths in these graphs are closely related to the Frobenius numbers of the lengths of their cycles. We then develop a procedure that extends an arbitrary meteor graph of length three to a graph of the type described in Lemma \ref{lm:normal-cover} (see Lemmas \ref{lem:first-type-edge-ext}, \ref{lem:second-type-edge-ext}, and \ref{lm:third-edge-ext}). 
These observations enable us to derive number-theoretic criteria characterizing when two meteor graphs of length three are strongly shift equivalent (Theorem \ref{numtheo}). Combining this with Krieger's theorem \cite[Theorem~4.2]{krieger} and Ara-Pardo's result \cite[Theorem 3.10]{arapar}, we show that Williams' Conjecture and Hazrat's Graded Morita Equivalence Conjecture hold for the class of meteor graphs of length three whose cycle lengths are pairwise coprime (Theorem \ref{thm21}).

The paper is organized as follows. In Section \ref{sec2}, we recall the fundamental concepts of $\mathbb Z$-monoids, the monoids $M_E$, $T_E$ associated with a directed graph $E$, in-splittings and out-splittings, and Leavitt path algebras together with their graded Grothendieck groups. These tools will be used for our detailed analysis of the talented monoid $T_E$, a $\mathbb{Z}$-monoid. In Section \ref{sec3}, we introduce the notion of connected finite graphs with disjoint cycles in normal form (Definitions \ref{def:quasi-normal-form} and \ref{def:normal-form}). Through a detailed analysis of in-splittings and out-splittings, we prove that every such graph can be transformed into a normal form graph by a finite sequence of in-splittings and out-splittings (Theorems \ref{thm:quasi-normal-form} and \ref{thm:normal-form}). In Section \ref{sec4}, based on Theorem \ref{thm:normal-form}, we establish number-theoretic criteria characterizing when meteor graphs of length three in normal form are strongly shift equivalent (Theorem \ref{numtheo}). In Section \ref{sec5}, Theorems \ref{thm:normal-form} and \ref{numtheo} allow us to reduce the problem of proving Williams' Conjecture and Hazrat's Graded Morita Equivalence Conjecture for meteor graphs of length three with pairwise coprime cycle lengths to their associated reduction graphs (Lemma \ref{lm:redu-SSE}), for which the graded Grothendieck groups of the corresponding Leavitt path algebras are easier to analyze. Using Theorem \ref{numtheo}, we then prove that both conjectures hold for meteor graphs of length three whose cycle lengths are pairwise coprime (Theorem \ref{thm21}). Combining this result with \cite[Theorem 4.3]{dohaznam}, we immediately conclude that both conjectures also hold for graphs with disjoint cycles containing exactly three cycles of pairwise coprime lengths (Corollary \ref{cor22}).
 

\section{Preliminaries: Monoids, Graphs and Algebras}\label{sec2}

The basic background needed to follow this paper is the one given in \cite[Section 2]{dohaznam}. Here we will recall only the most basic definitions, and refer the reader to \cite{dohaznam} for details.

\subsection{Monoids, \texorpdfstring{\(\mathbb{Z}\)}{ℤ}-monoids and order ideals}

Throughout this paper, we will work with commutative monoids, with additive notation and neutral element denoted by $0$.

The \emph{algebraic preorder} on a commutative monoid $M$ is given by setting $x\leq y$ if $y=x+z$ for some $z\in M$. 
We will deal with \emph{$\mathbb Z$-monoids}, which are commutative monoids endowed with an action $x\mapsto {^n}x$ of the group $\mathbb Z$ on $M$ by monoid automorphisms. A monoid homomorphism $\phi: M\longrightarrow M'$ is called a {\it $\mathbb Z$-monoid homomorphism} if $\phi(^n x) = $$^n\phi(x)$ for all $x\in M$ and $n\in \mathbb Z$. Similar definitions can be given for an action of an abelian group $\Gamma $ on $M$, see \cite{hazbk}.


The set of natural numbers is denoted by $\mathbb{N}=\left\{0,1,2,\ldots\right\}$. Under the usual sum, it is the free monoid generated by a single element. One of the $\mathbb Z$-monoids we encounter in this paper is the following. 

\begin{deff}\label{cycmond}	Let $k$ be a positive integer. The monoid $T=\bigoplus_{i=1}^k \mathbb N$, with the action of $\mathbb Z$ defined by ${}^1(a_1,\dots,a_{k-1},a_k)=(a_k,a_1\dots, a_{k-1})$, is called the \emph{$\mathbb Z$-cyclic monoid of rank $k$}.
\end{deff}

Let $M$ be a  $\mathbb Z$-monoid.  A {\it $\mathbb Z$-order ideal} of $M$ is a submonoid $I$ of $M$ which is closed under the action of $\mathbb Z$ and it is hereditary in the sense that $x\le y$ and
$y\in I$ imply $x\in I$. 

Let $I$ be a $\mathbb Z$-order ideal of $M$. Define an equivalence relation $\sim_I$ on $M$ as follows: For $a, b\in M$, there exist $x, y\in I$ such that $a + x = b +y$. This is a congruence relation and thus one can form the quotient $\mathbb Z$-monoid $M/{\sim}$ which is denoted by $M/I$.

For $\{a_1, \dots a_k\} \subseteq  M$, we denote the  $\mathbb Z$-order ideal of $M$ generated by the elements $a_i$ by $\langle a_1, \dots, a_k \rangle $. It is easy to see that 
\[\langle a_1, \dots, a_k \rangle=\Big \{ x \in M \mid x \leq \sum_{(i_1, \dots, i_k)\in \mathbb Z^k} {}^{i_1} a_{1} + \dots + {}^{i_k} a_{k} \Big \}.\]

\subsection{Graphs and associated monoids}

A (directed) graph $E = (E^0, E^1, s, r)$ consists of two disjoint sets $E^0$ and $E^1$, called \emph{vertices} and \emph{edges} respectively, together with two maps $s, r: E^1 \longrightarrow E^0$.  The vertices $s(e)$ and $r(e)$ are referred to as the \emph{source} and the \emph{range} of the edge~$e$, respectively.  A graph $E$ is called {\it row-finite} if $s^{-1}(v)$ is finite for all $v\in E^0$. It is called {\it finite} if both $E^0$ and $E^1$ are finite.

A {\em sink} in a graph $E$ is a vertex $v \in E^0$ with $s^{-1}(v) = \emptyset$; a {\em source} is a vertex $v \in E^0$ with $r^{-1} (v) =\emptyset$.  A  finite graph $E$ is {\em essential} if it has no sinks and sources.  

A (finite) \emph{path} in a graph $E$ is a string
$p=e_1\cdots e_n$ of edges $e_i\in E^1$ such that $r(e_i) = s(e_{i+1})$ for all $i$. The \emph{length} of the path $p = e_1 \cdots e_n$ is $n$, and is denoted by $|p|$.  
The source and range maps on edges are extended to paths as
\[s(p)=s(e_1)\qquad\text{and}\qquad r(p)=r(e_n).\]
Vertices are  regarded as paths of length $0$, with each vertex coinciding with its source and its range. We denote by $\text{Path}(E)$ the set of all paths in $E$.

A vertex $v\in E^0$ is said to \emph{lie} on a path $p$ if $v$ is the source or the range of one of the edges which comprise $p$. The set of vertices that lie on $p$ is denoted by $p^0$.




Next, we define the \emph{talented monoid} $T_E$ of $E$, which is believed to encode the graded structure of the Leavitt path algebra $L_{K}(E)$ (see Conjecture~\ref{conjehfyhtr}). As will become apparent later in the paper, it also serves as a bridge between symbolic dynamics and the theory of Leavitt path algebras.

\begin{deff}[{\cite[Page 436]{hazli}}]\label{talentedmon}
For a row-finite graph $E$, the \emph{talented monoid} of $E$, denoted by $T_E$, is the commutative 
	monoid generated by $\{v(i) \mid v\in E^0, i\in \mathbb Z\}$, subject to
	\[v(i)=\sum_{e\in s^{-1}(v)}r(e)(i+1)\]
	for every $i \in \mathbb Z$ and every $v\in E^{0}$ that is not a sink. The additive group $\mathbb{Z}$ of integers acts on $T_E$ via monoid automorphisms by shifting indices: For each $n,i\in\mathbb{Z}$ and $v\in E^0$, define ${}^n v(i)=v(i+n)$, which extends to an action of $\mathbb{Z}$ on $T_E$. Throughout the paper we denote the elements $v(0)$ in $T_E$ by $v$.
\end{deff}

The following important result from \cite{arahazrat} will be freely used in Section \ref{sec5}.

\begin{thm}[{\cite[Corollary 5.8]{arahazrat}}]
\label{thm:cancellation} The talented monoid $T_E$ is cancellative for every row-finite graph $E$.
\end{thm}

\subsection{Leavitt path algebras}\label{leviig}
The Leavitt path algebra $L_K(E)$ of a graph $E$ with coefficients in a field $K$ was introduced by Abrams and Aranda Pino in \cite{ap:tlpaoag05}, and independently by the first author, Moreno and Pardo in \cite{amp}. Leavitt path algebras generalize the Leavitt algebras $L_K(1, n)$ introduced in \cite{leav:tmtoar}. They are intimately related to graph $C^*$-algebras (see \cite{CunKri80, r:ga}) and have strong connections with other areas of algebra, including representation theory and even chip-firing games (\cite{AH23,HNam,HNam1}). We refer the reader to \cite{a:lpatfd} and \cite{lpabook} for a detailed history and overview of Leavitt path algebras.

\begin{deff}[{\cite[Definition 1.2.3]{lpabook}}]\label{LPAs}
For a row-finite graph $E = (E^0,E^1,s,r)$ and any  field $K$, the \emph{Leavitt path algebra} $L_{K}(E)$ {\it of the graph}~$E$
\emph{with coefficients in}~$K$ is the $K$-algebra generated
by the union of the set $E^0$  and two disjoint copies $E^1$, say $E^1$ and $\{e^*\mid e\in E^1\}$, satisfying the following relations for all $v, w\in E^0$ and $e, f\in E^1$:
	\begin{itemize}
		\item[(1)] $v w = \delta_{v, w} w$;
		\item[(2)] $s(e) e = e = e r(e)$ and $e^*s(e) = e^* = r(e)e^*$;
		\item[(3)] $e^* f = \delta_{e, f} r(e)$;
		\item[(4)] $v= \sum_{e\in s^{-1}(v)}ee^*$ for any  vertex $v$ that is not a sink;
	\end{itemize}
	where $\delta$ is the Kronecker delta.
\end{deff}
Notice that $L_K(E)$ is a $\mathbb{Z}$-graded $K$-algebra:  
$L_K(E)= \bigoplus_{n\in \mathbb{Z}}L_K(E)_n$,  where for each $n\in \mathbb{Z}$, the homogeneous component of degree $n$ of $L_K(E)$ is the set \ $$ \text{span}_K \{pq^*\mid p, q\in \text{Path}(E), r(p) = r(q), |p|- |q| = n\}.$$ Here, if $q = v \in E^0$, we set $q^*=v$, and if  $q=e_1\cdots e_m$, with $e_i\in E^1$, then $q^*=e_m^*\cdots e_1^*$. 

While graph $C^*$-algebras have been classified using K-theoretic invariants (see \cite{eilers, ef}), the search for a complete invariant for the classification of Leavitt path algebras remains an active area of research.  Hazrat's Graded Morita Equivalence Conjecture~(\cite{mathann,hazd}, \cite[\S 7.3.4]{lpabook})  predicts, roughly speaking, that the graded Grothendieck group $K_0^{\gr}$ classifies, up to graded Morita equivalence, Leavitt path algebras of finite graphs. The conjecture is closely related to Williams' conjecture (see~\S\ref{dynref}).

In order to state Hazrat's Graded Morita Equivalence Conjecture, we first recall the definition of the graded Grothendieck group of a $\Gamma$-graded ring. Given a $\Gamma$-graded ring $A$ with identity and a graded finitely generated projective (right) $A$-module $P$, let $[P]$ denote the class of graded $A$-modules graded isomorphic to $P$. Then the monoid  
\[\mathcal V^{\gr}(A)= \{[P] \mid  P  \text{ is a graded finitely generated projective A-module}\}\]
has a $\Gamma$-monoid structure defined as follows: for $\gamma \in \Gamma$ and $[P]\in \mathcal V^{\gr}(A)$, $\gamma .[P]=[P(\gamma)]$, where $[P(\gamma)]$ is the $\gamma$-twist of $P$.

The group completion of $\mathcal V^{\gr}(A)$ is called the \emph{graded Grothendieck group} and is denoted by $K^{\gr}_0(A)$.  The $\Gamma$-monoid structure on $\mathcal V^{\gr}(A)$ induces 
a $\mathbb{Z}[\Gamma]$-module structure on the group $K^{\gr}_0(A)$. In particular, the graded Grothendieck group of a $\mathbb{Z}$-graded ring has a natural $\mathbb{Z}[x,x^{-1}]$-module structure. 

By \cite[Proposition 5.7]{arahazrat}, $T_E$ is $\mathbb{Z}$-monoid isomorphic to the monoid $\mathcal V^{\gr}(L_K(E))$, where $E$ is a row-finite graph and $K$ is an arbitrary field. 

\begin{conj}[{\sc Hazrat's Graded Morita Equivalence Conjecture}]\label{conjehfyhtr}
Let $E$ and $F$ be finite graphs, and $K$ a field. Then the following statements are equivalent:
	
$(1)$ The Leavitt path algebras $L_K(E)$ and $L_K(F)$ are graded Morita equivalent;
		
		
$(2)$ There is an order-preserving 
		$\mathbb Z[x,x^{-1}]$-module isomorphism
		$K_0^{\gr}(L(E))\rightarrow K_0^{\gr}(L(F))$;

$(3)$ The talented monoids $T_E$ and $T_F$ are $\mathbb{Z}$-isomorphic.	
\end{conj}
  
We refer the reader to \cite{ART24, arapar, guido, guidowillie, guido1,toke, CortHaz, eilers2, ef, vas, vas2025} for some works on the graded Morita equivalence conjecture and related topics.

\subsection{Symbolic Dynamics}\label{dynref}
Williams introduced the notion of shift equivalence for matrices in \cite{williams} (see also \cite[\S7]{lindmarcus}) as part of a program to develop computable methods for deciding whether two shifts of finite type are conjugate. Recall that two square nonnegative integer matrices $A$ and $B$ are called {\it elementary shift equivalent}, and denoted by $A\sim_{ES} B$, if there are nonnegative matrices $R$ and $S$ such that $A=RS$ and $B=SR$. Two square nonnegative integer matrices $A$ and $B$ are called {\it strongly shift equivalent}, denoted by $A\sim_{SSE}B$, if there is a sequence of finite elementary shift equivalences from $A$ to $B$. The weaker notion of shift equivalence is defined as follows. The nonnegative integer matrices $A$ and $B$ are called {\it shift equivalent}, denoted by $A\sim_{SE} B$, if there are nonnegative matrices $R$ and $S$ such that $A^l=RS$ and $B^l=SR$, for some $l\in \mathbb N$, and $AR=RB$ and $SA=BS$. 

The {\em adjacency matrix} $A_E\in M_{E^0}(\mathbb{N})$ of a  graph $E$ provides the link between symbolic dynamics and  graphs.  By definition, $A_E$ is a square matrix with \[ A_E(v,w) = \left| s^{-1}(v) \cap r^{-1}(w)\right| .\]  Conversely,  any $A \in M_n(\mathbb{N})$ can be interpreted as the adjacency matrix of a graph $E$ with $|E^0| = n$, where $E^1$ consists of precisely $A(i,j)$ edges from vertex $i$ to vertex $j$, for all $1 \leq i, j \leq n$.

Identifying a square nonnegative integer matrix with its associated graph (and the graph with its adjacency matrix), Williams showed that two matrices $A_E$ and $A_F$ are strongly shift equivalent if and only if $E$ can be transformed into $F$ by a finite sequence of certain graph moves. We briefly recall these graph moves, as they will be used extensively throughout the paper.

\subsection*{In-splitting}

\begin{deff}[{\cite[Definition 6.3.20]{lpabook}}]\label{def:insplit}
Let $E$ be a graph. For each $v\in E^0$ with $r^{-1}(v)\neq \emptyset$, take a partition $\left\{\mathcal{E}^v_1,\ldots,\mathcal{E}^v_{m(v)}\right\}$ of $r^{-1}(v)$. We form a new graph $F$ as follows.  Set
	\[F^0=\left\{v_i \mid v\in E^0,1\leq i\leq m(v)\right\}\cup\left\{v \mid r^{-1}(v)=\emptyset \right\}\]
	\[F^1=\left\{e_j \mid e\in E^1,1\leq j\leq m(s(e))\right\}\cup\left\{e \mid r^{-1}(s(e))=\emptyset \right\},\]
	and define the source and range maps  as follows: If $r^{-1}(s(e))\neq\emptyset $, and 
	$e\in\mathcal{E}^{r(e)}_i$, then
	\[s(e_j)=s(e)_j,\qquad r(e_j)=r(e)_i, \text{ for all } 1\leq j \leq m(s(e)).\]
	If $r^{-1}(s(e))=\emptyset$, set $s(e)$ as the original source of $e$, and $r(e)=r(e)_i$, where 
	$e\in\mathcal{E}^{r(e)}_i$.
	
The graph $F$ is called an \emph{in-splitting} of $E$, and conversely $E$ is called an \emph{in-amalgamation} of $F$.
\end{deff}

\subsection*{Out-splitting}

The notions dual to those of in-splitting and in-amalgamation are called \emph{out-splitting} and \emph{out-amalgamation}, respectively. Given a graph $E=(E^0,E^1,s,r)$, the \emph{transpose graph} is defined as $E^*=(E^0,E^1,r,s)$.

\begin{deff}[{\cite[Definition 6.3.23]{lpabook}}]\label{def:outsplit}
	A graph $F$ is an \emph{out-splitting} (\emph{out-amalgamation}) of a graph $E$ if $F^*$ is an in-splitting (in-amalgamation) of $E^*$.
	
\end{deff}

We emphasize that in this paper, we require the sets $\mathcal E^i_v$ used in in- or out-splitting to be non-empty.


Let $E$ be an essential graph with the discrete topology on $E^1$. Define the topological edge shift $(X_E, \sigma_E)$ by setting:
$$X_E:= \{x=(x_n)_{n\in \mathbb{Z}}\mid x_n\in E^1 \text{ such that } s(x_n) = r(x_{n+1})\}\subseteq (E^1)^{\mathbb{Z}}$$
where $\sigma_E: X_E\longrightarrow X_E$ is the shift map with $\sigma_E (x)_n = x_{n+1}$ for all $n\in \mathbb{Z}$.\medskip

Two edge shifts $X_E$ and $X_F$ are {\it conjugate}, denoted $X_E\cong X_F$, if there is a homeomorphism $f: X_E\longrightarrow X_F$ such that the following diagram
$$\xymatrix{X_E\ar[d]_{\sigma_E} \ar[r]^f& X_F \ar[d]^{\sigma_F}\\
	X_E \ar[r]^f & X_F}$$
is commutative.

We are now able to present the precise form of Williams' celebrated criterion for the conjugacy of edge shifts  established in \cite{williams}.

\begin{theorem}[{\sc Williams}~\cite{williams}]\label{willimove}
Let $E$ and $F$ be essential graphs, and let $A_E$ and $A_F$ be the adjacency matrices of $E$ and $F$, respectively. Then, the following statements are equivalent:\medskip

$(1)$ $X_E$ is conjugate to $X_F;$\medskip

$(2)$ $A_E\sim_{SSE}A_F$;\medskip

$(3)$ $E$ can be obtained from $F$ by a sequence of in-splittings, out-splittings,  in-amalgamations, and out-amalgamations.
\end{theorem}



 The following conjecture was proposed by Williams \cite{willwrong} in 1974.

\begin{conj}[{\sc Williams' conjecture \cite{willwrong}}]
$A_E\sim_{SE}A_F\ \Longleftrightarrow\ A_E\sim_{SSE}A_F$ for essential graphs $E$ and $F$.
\end{conj}

It was not until 25 years later that Kim and Roush produced a counterexample \cite{kimroush99}. More precisely, Kim and Roush showed that there are essential graphs $E$ and $F$ of order $7$ with $A_E\sim_{SE} A_F$ while $A_E\nsim_{SSE} A_F$. However, identifying classes of edge shifts for which shift equivalence and strong shift equivalence coincide remains an important open problem.

Although strong shift equivalence characterizes conjugacy of edge shifts, there is no general algorithm for deciding whether two square matrices are strongly shift equivalent. In contrast, the weaker notion of shift equivalence is more tractable, as shown by Krieger.  In~\cite{krieger}, he introduced an invariant, now known as \emph{Krieger's dimension group}, for classifying edge shifts up to shift equivalence (see \cite[Theorem~4.2]{krieger} for details). Surprisingly, the first author and Pardo \cite[Theorem 3.10]{arapar} showed that Krieger's dimension group from dynamics can be realized as the graded Grothendieck group of a Leavitt path algebra. This provides a link between the theory of Leavitt path algebras and symbolic dynamics~\cite{arapar,hazbk, hazd, vas2025}.  

The following theorem will be used to prove the main result of Section \ref{sec5} (Theorem \ref{thm21}).

\begin{thm}[{\cite[Theorem~4.2]{krieger} and \cite[Theorem 3.10]{arapar}}]\label{h99}
Let $E$ and $F$ be essential graphs with adjacency matrices $A_E$ and $A_F$, respectively, and let $K$ be an arbitrary field. Then
$A_E$ is shift equivalent to $A_F$ if and only if $K_0^{\gr}(L_K(E)) \cong K_0^{\gr}(L_K(F))$ via an order-preserving $\mathbb Z[x,x^{-1}]$-module isomorphism.
\end{thm}

\section{Normal form for general finite graphs with disjoint cycles}\label{sec3}

In this section, we introduce the notion of connected finite graphs with disjoint cycles in normal form (Definitions \ref{def:quasi-normal-form} and \ref{def:normal-form}) and show that any connected finite graph with disjoint cycles can be transformed into a normal form graph via a finite sequence of in-splittings and out-splittings (Theorem \ref{thm:normal-form}).\medskip

Let $E$ be an arbitrary graph. A path  $c= e_{1} \cdots e_{n}\in \text{Path}(E)$ of positive length is a \textit{cycle based at} the vertex $v$ if $s(c) = r(c) =v$ and the vertices $s(e_1), s(e_2), \hdots, s(e_n)$ are distinct. 


\begin{deff}
A graph $E$ is called a {\it graph with disjoint cycles} if every vertex in $E$ is the base of at most one cycle.    
\end{deff}

The following definition is very useful in the proof of the main result of this section.

 \begin{deff}\label{def3.1}
Let $E$ be a finite graph with disjoint cycles. We denote by $\mathcal{C}_E$  the set of (disjoint) cycles of $E$. Let $C_1$ and $C_2$ be two cycles in $\mathcal{C}_E$. We write $C_1\ge C_2$ in case there is a path in $E$ from a vertex of $C_1$ to a vertex of $C_2$. 
Since the cycles in $E$ are exclusive cycles,  $\le $ is a partial order on $\mathcal {C}_E$, and $(\mathcal {C}_E, \le)$ becomes a finite poset. Hence we can apply to it the usual definitions in a finite poset. Conversely, given a finite poset $(P,\le)$, there exists a finite graph with disjoint cycles $G$ such that $(\mathcal{C}_G,\le ) \cong (P,\le)$. To see this, take the Hasse graph $H$ of the poset $P$, and replace each vertex in $H$ with a $cycle$ of length $1$. Then the corresponding directed graph $G$ has the desired property. 
          
We define the {\it height} of $C\in \mathcal{C}_E$ as the maximum length of a chain
        $$C=C_0 > C_1 >  \cdots > C_l.$$
        Here the length of the above chain is, by definition, the integer $l$.
        Given $A,B\in \mathcal{C}_E$ such that $B\le A$, define the {\it distance} $d(A,B)$ as the height of $A$ in the interval $[B,A]$. That is, the distance $d(A,B)$ is the largest $r$ such that there is a chain
        $$A=C_0 > C_1> \cdots > C_r=B.$$

        A path $\alpha = e_1e_2 \cdots e_n$ in $E$ is called a {\it trail} if  $s(e_1), r(e_n) \in \mathcal{C}_E^0$ and $s(e_i) \notin \mathcal{C}_E^0$ for all $1< i \le  n$, where $\mathcal{C}_E^0 = \cup_{C\in \mathcal{C}_E} C^0$. 
For a trail $\alpha$, its {\it interior} is the set $\alpha^0 \setminus \mathcal{C}_E^0 .$  A vertex $v$ in $E$ is called an {\it interior vertex} if $v$ is in the interior of some trail.

A cycle $C$ of height $i$ in $\mathcal{C}_E$ will be called an {\it $i$-cycle}. Hence a $0$-cycle is just a minimal element of $\mathcal{C}_E$. 
    \end{deff}

 \begin{deff}\label{def3.2} Let $C$ be a cycle in a graph $E$. We say that an edge $e$ is an {\it entrance} of $C$ if $e$ is not in the cycle and $r(e)\in C^0$. Similarly, we say that $e$ is an {\it exit} of $C$ is $e$ is not in the cycle and $s(e)\in C^0$.     \end{deff}

Let $E$ be a connected finite graph with disjoint cycles. 
In the following, we will label the vertices in the cycles $C$ of $E$ by natural numbers $1,2,\dots , |C|$.
We will refer to the label of a vertex $v$ in $C$ as the {\it index} of $v$. 

We now come to a crucial observation. Let $E$ be a connected finite graph with disjoint cycles. If we perform an in-splitting or an out-splitting on $E$, then the cycle structure of the graph remains unchanged; only some of the trails are modified. This motivates the following definition.

\begin{deff}\label{def:cycle-structure}
Let $E$ and $F$ be two connected finite graphs with disjoint cycles. We say that $E$ and $F$ {\it have the same structure of cycles} if there is an isomorphism of posets $\varphi \colon (\mathcal C_E,\le) \to  (\mathcal C_F,\le)$ such that $|C| = |\varphi (C)|$ for all $C\in \mathcal C_E$.  
\end{deff}

The following fact shows that shift equivalence (and in particular strong shift equivalence) captures the cycle structure of connected finite graphs with disjoint cycles.

\begin{prop}[{\cite[Theorem 3.2]{dohaznam}}]
\label{prop:SSE-implies-same structure}
Let $E$ and $F$ be two connected finite graphs with disjoint cycles. If $E$ and $F$ are shift equivalent (that is, their adjacency matrices are shift equivalent), then $E$ and $F$ have the same structure of cycles. 
\end{prop}

Having the same cycle structure defines an equivalence relation on the set of connected finite graphs with disjoint cycles.
For each equivalence class $\mathfrak C$, we fix a collection of cycles $C_1,\dots ,C_n$, together with a fixed labeling of the vertices in each of the cycles $C_i$ as described above.
Every graph $E\in \mathfrak C$ is isomorphic to a graph $E'$ whose set of cycles is precisely $C_1,\dots ,C_n$. We say that $E'$ is a connected finite graph with 
{\it normalized set of disjoint cycles}, and write $$\mathcal C_{E'} = \{C_1,\dots ,C_n\}.$$ 

We next introduce the definition of connected finite graphs with disjoint cycles in quasi normal form.
  
\begin{deff}\label{def:quasi-normal-form}
Let $E$ be a connected finite graph with a normalized set of disjoint cycles $\mathcal{C}_E= \{C_1,C_2,\ldots,$ $C_n\}.$ We say $E$ is in {\it quasi-normal form} if the following conditions are satisfied:

$(1)$ Any two distinct trails in $E$ have no common interior vertices;

$(2)$ All trails between cycles of $E$ start and end in a vertex of index $1$. Concretely, there exist vertices $v_i \in C^0_i$, $(i = 1, \ldots, n)$, such that every trail from $C_i$ to $C_j$ starts at  $v_i$ and ends at  $v_j$, for all $i\neq j$. The vertices of $C_i$ are cyclically arranged as follows $v_i=v^i_1\to v^i_2\to \cdots \to v^i_{|C_i|}\to v_i$, with $\text{index} (v^i_j)=j$, $i=1,\dots , n$, $j=1,\dots , |C_i|$. 
\end{deff}
 
For clarification, we illustrate Definitions \ref{def3.1}, \ref{def3.2} and \ref{def:quasi-normal-form} by presenting the following example.

\begin{example}\label{exa3.4}
Consider the following  graph $E$:

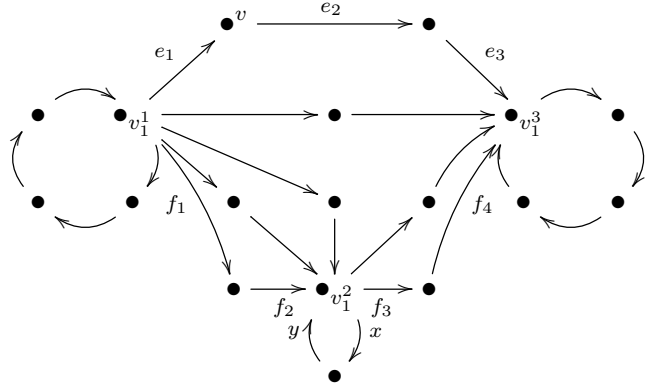
\[\xymatrixrowsep{1.5pc}\xymatrixcolsep{1.5pc}\xymatrix{&&&\bullet^{v}\ar[rr]^{e_2}&&\bullet\ar[dr]^{e_3}\\&\bullet\ar@/^0.8pc/[r]&\bullet_{v^1_1}\ar@/^0.8pc/[d]\ar[ur]^{e_1}\ar[rr]\ar[dr]\ar[drr]\ar@/^0.5pc/[ddr]_{f_1}&&\bullet\ar[rr]&&\bullet_{v^3_1}\ar@/^0.8pc/[r]&\bullet\ar@/^0.8pc/[d]&\\&\bullet\ar@/^0.8pc/[u]&\bullet\ar@/^0.8pc/[l]&\bullet\ar[dr]&\bullet\ar[d]&\bullet\ar@/^0.5pc/[ur]&\bullet\ar@/^0.8pc/[u]&\bullet\ar@/^0.8pc/[l]&\\&&&\bullet\ar[r]_{f_2}&\bullet_{v^2_1}\ar@/^0.8pc/[d]^x\ar[ur]\ar[r]_{f_3}&\bullet\ar@/^0.5pc/[uur]_{f_4}&&&\\&&&&\bullet\ar@/^0.8pc/[u]^y}\]
It is obvious that $E$ is a quasi-normal form graph. We also have that $e_1e_2e_3$ is a trail in $E$,  $v$ is an interior vertex of  $e_1e_2e_3$, and $f_1f_2f_3f_4$ is not a trail in $E$. Let $A$ and $B$ be two cycles in $E$ containing $v^1_1$ and $v^3_1$, respectively.
We then have that $f_4$ is an entrance of  $B$  and $f_1$ is an exit of $A$. Moreover, we have that the distance $d(A, B)$ is equal to $2$.
\end{example}

 
   
   The following theorem gives that every connected finite graph with disjoint cycles can be transformed to a quasi-normal form graph using a finite sequence of in-splittings and out-splittings.
   
  \begin{theorem}\label{thm:quasi-normal-form}
  Any connected finite graph $E$ with disjoint cycles can be transformed, using only graph isomorphisms, and in- and out-splittings, into a graph with disjoint cycles in quasi-normal form.
  \end{theorem} 
   
   \begin{proof} We may assume that $E$ is a connected finite graph with normalized set of disjoint cycles $\mathcal C_E= \{C_1,\dots , C_n\}$. Let $v_i\in C_i$ be the corresponding vertices of index $1$.

 Let $C_i$ be a $0$-cycle in $E$. For a trail $\alpha = e_1 \cdots e_t$  ending at vertex $v$ in $C_i$, with $v\ne v_i$, we perform the in-splitting at  $v$ with a partition  $\{\mathcal{E}_1,\mathcal{E}_2\}$ where  $\mathcal{E}_1 = \{e_t\}$  and $\mathcal{E}_2= r^{-1}(v)\setminus \{e_t\}$. Doing this, we will  increase the length of $\alpha$ by 1 and make the trail end at vertex $v'$ where $v \to v'$ in the cycle $C'_i\cong C_i$, but all other trails ending with an edge different to $e_t$ remain unchanged. Here we identify the transformed cycle $C_i'$ with $C_i$. Hence $v'$ is ``closer" to $v_i$ than $v$, in the sense that the path within $C_i$ from $v'$ to $v_i$ has length $l-1$, where $l$ is the length of the path within $C_i$ from $v$ to $v_i$. If $v'\ne v_i$, we can repeat this process with the new trail, and arrive at the same situation but now with a vertex $v''$ such that the length of the path within $C_i$ from $v''$ to $v_i$ is $l-2$. This process will eventually lead to a graph in which all trails ending in the edge $e_t$ in the original graph are enlarged to trails ending at $v_i$. Moreover all trails not ending with the edge $e_t$ in the original graph remain unchanged.   
 
 Doing this process for all the entrances of $C_i$, we can make all the trails ending at $C_i$ to end at the vertex $v_i$. Applying this process for all the $0$-cycles, we will obtain, using a finite sequence of in-splittings, a new graph such that all the trails ending in a $0$-cycle $C_i$ end at the same vertex $v_i$ of $C_i$.

 After that, let $C_j$ be a $1$-cycle (seen in the transformed graph). For any entrance $e$ of $C_j$, with $r(e)= v\ne v_j$, we can perform the in-splitting at  $v$ with a partition  $\{\mathcal{E}_1,\mathcal{E}_2\}$ where  $\mathcal{E}_1 = \{e\}$  and $\mathcal{E}_2= r^{-1}(v)\setminus \{e\}$. Doing this, we will  increase the length of any trail ending with the edge $e$ by 1 and make the trail end at a  vertex $v'$ where $v \to v'$ in the cycle $C_j$, as before. Note that this process may increase the number of trails ending in a $0$-cycle, but all the trails ending in a $0$-cycle $C_k$ will end in the vertex $v_k$, as desired.

  Repeating this process a finite number of times, we arrive at a graph where all trails ending at $C_j$ end at the same vertex $v_j$. Applying this process for all the $1$-cycles, we arrive at a graph such that all the trails ending at a $t$-cycle $C_i$, for $t=0,1$, end at the selected vertex $v_i\in C_i^0$.  

Now a routine induction argument, based on the above observations, shows that we can reduce our original graph, by a finite sequence of in-splittings, to a graph with the same cycles in which all the trails connecting a cycle $C_i$ to a cycle $C_j$ end at the same selected vertex $v_j\in C_j^0$.

We provide an example of the process with the graph below.
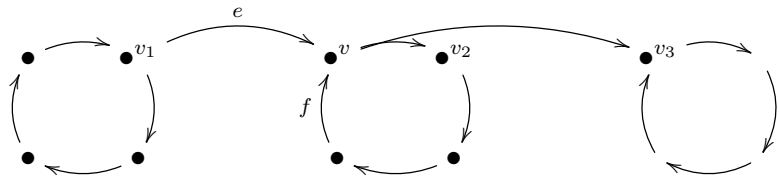
\[\xymatrix{
\bullet \ar@/^/@{->}[r] & \bullet^{v_1} \ar@/^/@{->}[d] \ar@/^1pc/@{->}[rr]^{e} &  & \bullet^{v} \ar@/^/@{->}[r] \ar@/^1pc/@{->}[rrr] & \bullet^{v_2} \ar@/^/@{->}[d] &  & \bullet^{v_3} \ar@/^/@{->}[r] & {} \ar@/^/@{->}[d] \\
\bullet \ar@/^/@{->}[u] & \bullet \ar@/^/@{->}[l] &  & \bullet \ar@/^/@{->}[u]^{f} & \bullet \ar@/^/@{->}[l] &  & {} \ar@/^/@{->}[u] & {} \ar@/^/@{->}[l]
}\]

Now in-splitting at $v$ with $\mathcal{E}_1 = \{e\}, \mathcal{E}_2=\{f\}$, we obtain the graph below:
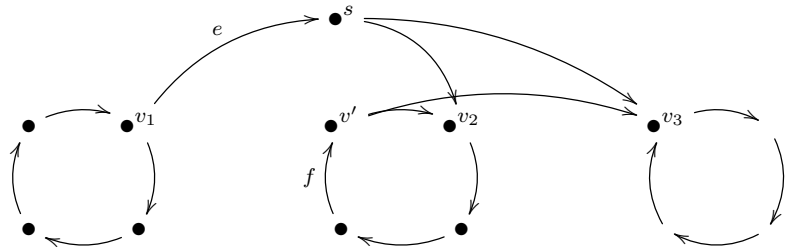
\[\xymatrix{
 &  &  & \bullet^{s} \ar@/^1pc/@{->}[rrrd]  \ar@/^1pc/@{->}[rd]&  &  &  &  \\
\bullet \ar@/^/@{->}[r] & \bullet^{v_1} \ar@/^/@{->}[d] \ar@/^1pc/@{->}[rru]^{e}  &  & \bullet^{v'} \ar@/^/@{->}[r]  \ar@/^1pc/@{->}[rrr] & \bullet^{v_2} \ar@/^/@{->}[d]    &  & \bullet^{v_3} \ar@/^/@{->}[r] & {} \ar@/^/@{->}[d] \\
\bullet   \ar@/^/@{->}[u]  & \bullet \ar@/^/@{->}[l] &  & \bullet  \ar@/^/@{->}[u]^{f} & \bullet \ar@/^/@{->}[l] &  & {} \ar@/^/@{->}[u] & {} \ar@/^/@{->}[l]
}\]
In this graph, we can see that there is a new trail from the first cycle to the last cycle but this trail still ends at vertex $v_3$ in the last cycle. Also, observe that in the transformed graph, the trail from the first to the second cycle ends at the vertex $v_2$.  

Similarly, using out-splittings we can make all trails starting at any cycle $C_i$ to start at vertex $v_i$. 

Now, we will prove that we can make all the trails have mutually disjoint interior vertices.
We denote by $I(E)$ the set of all interior vertices of trails in $E$. 
For $v\in I(E)$, we define two distances $d(v)=d_E(v)$ and $\ell(v)= \ell _E(v)$ as follows. Let $A_{v}$ be the set of all paths starting at $v$, ending at a vertex which lies in $\mathcal{C}_E^0$, and containing no vertices from cycles in $\mathcal{C}_E$ except for the final vertex. Let $B_{v}$ be the set of all paths starting at a vertex which lies in $\mathcal{C}_E^0$, ending at $v$, and contains no vertices from cycles in $\mathcal{C}_E$ except for the initial vertex. We define
\begin{center}
$d(v):= \max\{|\alpha| \mid \alpha \in A_v\}$  and  $\ell(v):= \max\{|\alpha| \mid \alpha \in B_v\}$.
\end{center}

We also consider the following constant
$$K(E) = \max \{ |\gamma | \mid \gamma \text{ is a trail in } E \}.$$

Assume that there is some vertex $u$ in the interior of a trail which is the range of more than one edge. Assuming the prior case, let us temporarily call such a vertex a {\it multi-range vertex}. We denote the set of all multi-range vertices of $E$ by $MR(E)$. Note that $MR(E)\subseteq I(E)$. 
Let $u\in MR(E)$ such that $\ell (u)$ is minimal
among all multi-range vertices. 
Performing an in-splitting at $u$ with a partition  $\{\mathcal{E}^u_1, \mathcal{E}^u_2,\ldots, \mathcal{E}^u_n\}$ of $r^{-1}(u) = \{e_1, \ldots, e_n\}$ where $\mathcal{E}^u_i = \{e_i\}$, we obtain new vertices $u_1,u_2, \ldots u_n$ which are the range of only one edge each. The multi-range vertices in this new graph $E'$ are the same ones as in the original graph, except for $u$, plus all vertices in $I(E)\cap r(s^{-1}(u))$, which are now multi-range vertices. We note that for any $v\in I(E)\cap r(s^{-1}(v))$, we have $\ell_{E'}(v) = \ell _E(v) > \ell_E(u)$.  Moreover the constant $K(E')$ corresponding to the transformed graph $E'$ is the same than the one for $E$, that is, $K(E')=K(E)$. Hence we obtain a new graph $E'$, with the same constant $K$ as the original one, but in which the number of vertices with minimal $\ell$ has decreased in one. Hence after a finite number of in-splittings we will get a new graph $E''$ such that either $MR(E'')= \emptyset$ or  
$$\min \{\ell (v) \mid v \in MR (E'')\} > \min \{\ell (v) \mid v \in MR (E)\}. $$
In addition, we have $K(E)= K(E'')$. Observing that $\ell _{E''} (v) \le K(E'') = K(E)$ for all $v\in I(E'')$, we see that this process must eventually stop, which means that, after applying a finite sequence of in-splittings to $E$ we arrive to a graph $H$ such that $MR(H)= \emptyset$, that is, a graph such that each interior vertex of a trail receives exactly one arrow. 

\begin{figure}[h]
\begin{center}
$\xymatrixrowsep{0.25pc}\xymatrixcolsep{0.001pc}\xymatrix{\bullet\ar[rrd]&&&&\bullet^{v}\ar[rr]\ar[rrd]&&\bullet\ar[ddrrr]&&&&&&&&&&&&&\bullet\ar[rr]&&\bullet_{u_1}\ar[rr]&&\bullet^{v}\ar[rr]\ar[rrd]&&\bullet\ar[ddrrr]&&&&&\\\bullet\ar[rr]&&\bullet^{u}\ar[rru]&&&&\bullet^{q}\ar[ldd]&&&&&&&&&&&&&\bullet\ar[rr]&&\bullet_{u_2}\ar[rru]&&&&\bullet^{q}\ar[ldd]&&&&\\\cdots&&&&&&&&&\bullet_{r}\ar@/^1pc/[rr]&&\bullet_{r_2}\ar@/^1pc/[ll]&\ar@{->}[rrrrrr]^{\text{ In-split at $u$ }}&&&&&&&\cdots&&&&&&&&&\bullet_{r}\ar@/^1pc/[rr]&&\bullet_{r_2}\ar@/^1pc/[ll]\\&&&&&\bullet_{q_2}\ar[rr]&&\bullet_{q_3}\ar[uul]\ar[rru]&&&&&&&&&&&&&&&&&\bullet_{q_2}\ar[rr]&&\bullet_{q_3}\ar[uul]\ar[rru]&&&\\\bullet\ar[rruuu]&&&&&&&&&&&&&&&&&&&\bullet\ar[rr]&&\bullet_{u_n}\ar[rruuuu]&&&&&&&&&&&&}$
\end{center} 
\caption{In-splitting at $u$ makes  a new multi-range vertex $v$ instead of $u$ which $\ell(v) >\ell(u)$.}
\end{figure}
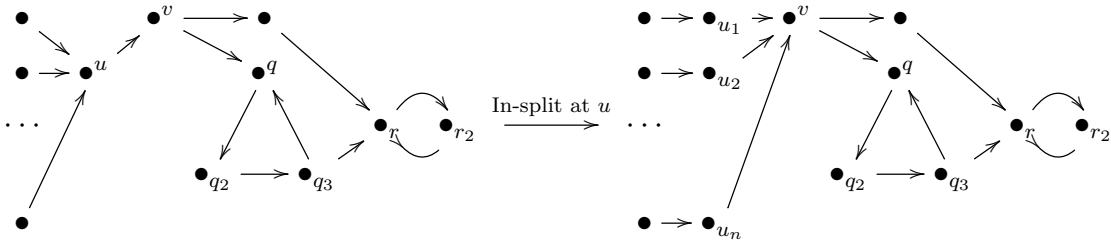

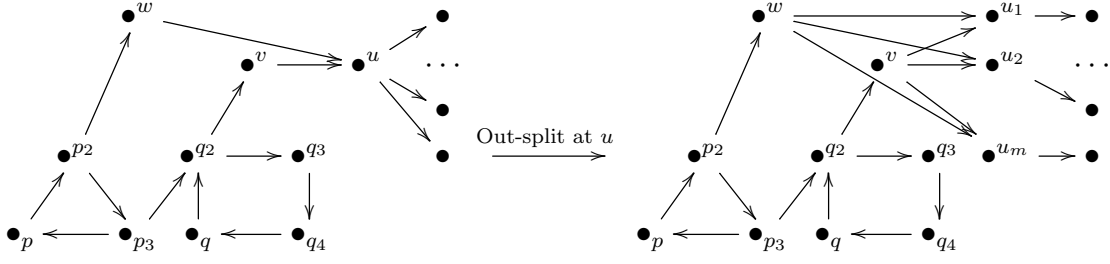
\begin{figure}[h]
\begin{center} 
$\xymatrixrowsep{0.3pc}\xymatrixcolsep{0.1pc}\xymatrix{&&\bullet^w\ar[rrrrd]&&&&&&\bullet&&&&&&&&&&\bullet^w\ar[rrrrd]\ar[rrrr]\ar[rrrrddd]&&&&\bullet^{u_1}\ar[rr]&&\bullet\\&&&&\bullet^v\ar[rr]&&\bullet^u\ar[rru]\ar[rrd]\ar[rrdd]&&\cdots&&&&&&&&&&&&\bullet^v\ar[rr]\ar[rru]\ar[rrdd]&&\bullet^{u_2}\ar[rrd]&&\cdots\\&&&&&&&&\bullet&&&&&&&&&&&&&&&&\bullet\\&\bullet^{p_2}\ar[ruuu]\ar[rdd]&&\bullet^{q_2}\ar[rr]\ar[ruu]&&\bullet^{q_3}\ar[dd]&&&\bullet&\ar@{->}[rrrrrr]^{\text{Out-split at $u$ }}&&&&&&&&\bullet^{p_2}\ar[ruuu]\ar[rdd]&&\bullet^{q_2}\ar[rr]\ar[ruu]&&\bullet^{q_3}\ar[dd]&\bullet^{u_m}\ar[rr]&&\bullet\\&&&&&&&&&&&&&&&&&&&&&\\\bullet_{p}\ar[ruu]&&\bullet_{p_3}\ar[ll]\ar[ruu]&\bullet_{q}\ar[uu]&&\bullet_{q_4}\ar[ll]&&&&&&&&&&&\bullet_{p}\ar[ruu]&&\bullet_{p_3}\ar[ll]\ar[ruu]&\bullet_{q}\ar[uu]&&\bullet_{q_4}\ar[ll]}$\end{center} 
\caption{Out-splitting at $u$ makes  new multi-source vertices $v, w$  instead of $u$ which $d(v), d(w) > d(u)$.}
\end{figure} 

Now a symmetric argument allows to transform our new graph $H$, using a finite sequence of out-splittings, to another graph such that all the vertices in the interior of the trails emit exactly one arrow. Observe that in this process, we will keep the property that there are no multi-range vertices, so our final graph will have the desired property that all trails have disjoint interior. 
\end{proof}

We next introduce the definition of connected finite graphs with disjoint cycles in normal form, which plays an important role in our analysis.   

   \begin{deff}\label{def:normal-form}
       Let $E$ be a connected finite graph with disjoint cycles in quasi-normal form. We call a {\it pair of cycles $(A, B)$ good} if, for any two trails $\alpha$ and $\beta$ from $A$ to $B$, the following conditions are satisfied: \begin{enumerate}
           \item If $|\alpha| - |\beta| \equiv 0 \pmod{ \gcd (|A|,|B|) }$ then  $|\alpha| = |\beta|$; 
           \item If $|\alpha| \equiv t \pmod{ \gcd (|A|, |B|) ) }$ then $|\alpha| \equiv t \pmod{|A||B| }$ for all $0\le t\ < \gcd (|A|, |B|)$. 
       \end{enumerate}   
              We say  $E$ is in {\it normal form} if all pairs of cycles $(A, B)$ in $E$ are good. 
   \end{deff}

 For clarification, we illustrate Definition \ref{def:normal-form} by presenting the following example.  

 \begin{example}\label{exa3.8}
Consider the following graph in quasi-normal form

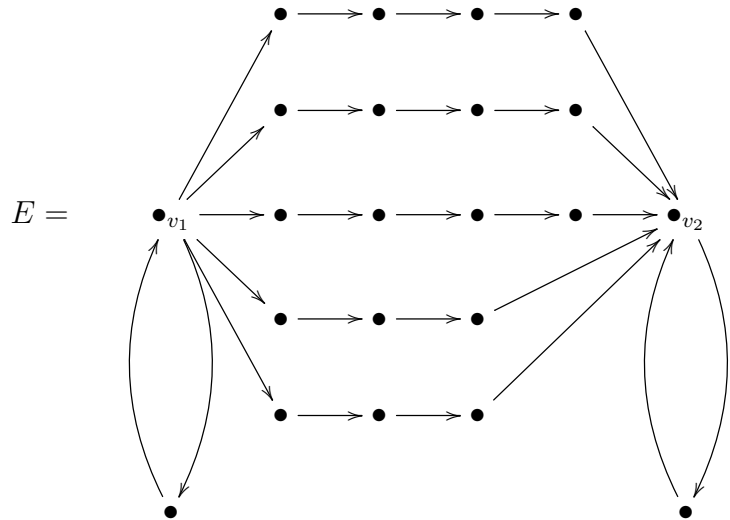
\[\xymatrix{&&\bullet\ar[r]&\bullet\ar[r]&\bullet\ar[r]&\bullet\ar[ddr]&\\
&&\bullet\ar[r]&\bullet\ar[r]&\bullet\ar[r]&\bullet\ar[dr]&\\
E = &\bullet_{v_1}\ar@/^1.3pc/[ddd]\ar[r]\ar[ur]\ar[uur]\ar[dr]\ar[ddr]&\bullet\ar[r]&\bullet\ar[r]&\bullet\ar[r]&\bullet\ar[r]&\bullet_{v_2}\ar@/^1.3pc/[ddd]\\
&&\bullet\ar[r]&\bullet\ar[r]&\bullet\ar[urr]&&\\
&&\bullet\ar[r]&\bullet\ar[r]&\bullet\ar[uurr]&&\\ &\bullet\ar@/^1.3pc/[uuu]&&&&&\bullet\ar@/^1.3pc/[uuu]\\}\]
Let $A$ and $B$ be two cycles in $E$ containing $v_1$ and $v_2$ respectively. Since $|A| = 2 = |B|$, and all trails in $E$ have length four or five, it follows that $E$ is in normal form. 
 \end{example}

For the proof of the main theorem of this section, we will rely on an easy number-theoretic fact. We include a proof for completeness. 

\begin{lemma}
\label{lem:number-theory}
Suppose that $d= \gcd (m,n)$ for positive integers $m,n$. Then given any family $c_1,\dots ,c_k$ of positive integers such that $c_i \equiv j \pmod{d}$ for all $1\le i\le k$, for a fixed $j$ such that $0\le j<d$, there exists a sufficiently large integer $N$ such that, for each $i=1,\dots , k$ there are non-negative integers $a_i,b_i$ such that $N= c_i + a_im+b_in$, and $N-j$ is a multiple of $nm$.
\end{lemma}

\begin{proof}
  Take integers $r,s$ such that $d= rm+sn$. Assume, without loss of generality, that $r>0$ and $s\le 0$.  Let $j$ be as in the statement, so that $c_i -j= z_id$ for non-negative integers $z_i$. Now take 
   	$$z= \mathrm{max} \{ z_1,\dots , z_k\}.$$
   	Then $z\ge 0$ and 
   	$$c_i + z(-s)n + nm = j+ z_irm + (z-z_i)(-s)n +nm = j + a_i'm+b_i'n,$$
   	where $a_i',b_i'$ are non-negative integers, not both zero. Set $M_i = a_i'm+b_i'n>0$, let $M'$ be the least common multiple of $M_1,\dots , M_k$, and set $M= M'nm$. Then we have $M= t_iM_inm$ for some $t_i\ge 1$,
    $i=1,\dots ,k$. 
   	Now define $N:=j+M$. Then $N-j= M$ is a multiple of $nm$, and we have
   	\begin{align*}
   		N-c_i & = (N-j) +(j-c_i) = M +(j-c_i) \\
   		& = t_iM_inm + z(-s)n + nm - M_i \\
   		& = (t_inm-1) M_i + z(-s)n +nm \\
   		& = [(t_inm-1)a_i']m + [(t_inm-1)b_i'+z(-s) + m]n.
   	\end{align*} 
   	Hence $N= c_i + a_im+b_in$ for each $i=1,\dots , k$, with $a_i = (t_inm-1)a_i'\ge 0$ and $b_i= (t_inm-1)b_i'+z(-s)+m > 0$, as desired.   \end{proof}

 The proof of Theorem \ref{thm:normal-form} repeatedly uses the following lemma.

\begin{lemma}\label{lem:basic-for-normal-form}
Let $E$ be in quasi-normal form, and let $\alpha$ be a trail from a cycle $A$ to a cycle $B$. For every $a,b \in \mathbb{N}$, there is a graph $E'$ in quasi-normal form, obtained from $E$ by in- and out-splittings, in which $\alpha$ is replaced by a trail of length $|\alpha| + a|A| + b|B|$, while every other trail between cycles $C,C'$ with $d(C,C') \leq d(A,B)$ is unchanged.
\end{lemma}

\begin{proof}
Let $\alpha = e_1\cdots e_s$ be a trail from $A$ to $B$. Let $v_A$ and $v_B$ be the distinguished vertices in $A$ and $B$ respectively. We first study the effect of doing an out-splitting at $v_A$. Take the partition $\mathcal E ^1=
   	\{e_1\} $, $\mathcal E^2 = s^{-1}(v_A)\setminus \{e_1\}$, of $s^{-1}(v_A)$. Let $E_1'$ be the corresponding out-splitting. Suppose that $A^0=\{v_A,v_2,\dots , v_r \}$ in cyclic order, that is $v_A\rightarrow v_2 \rightarrow \cdots \rightarrow v_r\rightarrow v_A$. Then we have a cycle $A'$ in $E_1'$ isomorphic to $A$,  that we identify with $A$. Note that if $r>1$, then $E_1'$ is no longer in quasi-normal form: there is a new trail $\alpha'$ from $A$ to $B$, with $|\alpha'|= |\alpha | +1$, which replaces $\alpha$, and with source $s(\alpha')= v_r$. All other trails of $E$ from $A$ to $B$ not starting with $e_1$ remain the same in $E_1'$. Since the trails from $A$ to $B$ have mutually disjoint interiors, all them remain invariant except for the trail $\alpha$. Observe that if $\beta$ is a trail from a cycle $A'$ to $A$, then there is an additional trail $\beta'$ from $A'$ to $B$. If $\beta = \beta_1 f$ for an edge $f$, with (necessarily) $r(f)= v_A$, then the new trail from $A'$ to $B$ has the
   	form $\beta' = \beta_1f_1 e_1\cdots e_s$, where the edge $f$ has been ``splitted" in the two edges $f_1$
   	and $f_2$ in the course of the out-splitting. However, it is obvious that $d(A',B) >d(A,B)$. Hence we observe that all trails from a cycle $C$ to a cycle $C'$, with $d(C,C')\le d(A,B)$, except the trail $\alpha$, remain invariant after this out-splitting. Observe that the new trails $\alpha'$ and $\beta'$, where $\beta$ is any trail from $A'$ to $A$, for a cycle $A'$, have no disjoint interiors: the vertex $s_{E_1'}(e_1)$ is an interior vertex of all these trails.
    Similarly, if $\beta_1$ is non-trivial, the vertex $r(\beta_1)$ is an interior vertex of the new trails $\beta'$ and $\beta_1f_2$.
    However, after an application of $s$ in-splittings, at the vertices $s(e_1), s(e_2), \dots , s(e_{s})$, and $|\beta_1|$ out-splittings at the interior vertices of $\beta_1$, as in the second part of the proof of Theorem \ref{thm:quasi-normal-form}, we obtain a new graph $E_1$ where all the trails have mutually disjoint interiors, where the trail $\alpha$ has been replaced by a new trail $\alpha'$ from $A$ to $B$ with $|\alpha'| = |\alpha| +1$, and where all the trails between cycles $C,C'$ with $d(C,C')\le d(A,B)$, except the trail $\alpha$, remain invariant. (Note that here we identify the original trail $\beta$ from $A'$ to $A$ with the new trail $\beta''$ from $A'$ to $A$ obtained by performing the $|\beta_1|$ out-splittings at the interior vertices of $\beta_1$.)    
   	
   	We can now apply an out-splitting to the graph $E_1$ at the vertex $v_r$ of $A$ using the first arrow of the trail $\alpha'$, to get a new trail $\alpha''$ such that $|\alpha''| = |\alpha |+2$ and $s(\alpha '')= v_{r-1}$. Continuing in this way, we arrive to a new graph $E_r$, which will be in quasi-normal form, such that the trail $\alpha$ has been substituted by a new trail $\alpha_r$ of the form
   	$g_1\cdots g_r e_1\cdots e_s$, where the vertices  $r(g_i)$, for $i=1,\dots ,r$ are all interior vertices of $E_r$, and all the other trails from $A$ to $B$ not starting with $e_1$ will remain the same in the new graph $E_r$.  As said before, all the trails from $A$ to $B$ distinct from $\alpha$ will remain the same in $E_r$. Note that $|\alpha_r| = |\alpha | + r$. Moreover, as before, we also get that, after the passage from $E$ to $E_r$, all trails between any two cycles $C$ and $C'$ such that $d(C,C')\le d(A,B)$, except for the initial trail $\alpha$, remain invariant. We highlight the important fact that the graph $E_r$ is again in quasi-normal form, because the new trail $\alpha_r$ starts again in the distinguished vertex $v_A$, and we have performed the corresponding in- and out-splittings at each step of the process.   
   	
   	In conclusion, given a trail $\alpha $ from $A$ to $B$ we can perform a sequence of out-splittings and in-splittings to $E$ in order to enlarge the length of $\alpha$ by any positive multiple of $|A|$, to get another graph in quasi-normal form, in which all the other trails from $A$ to $B$ remain the same as in the original graph. More generally, all the trails from a cycle $C$ to another cycle $C'$, with $d(C,C')\le d(A,B)$, except for the trail $\alpha$, remain the same as in the original graph. 
   	
   	Now the same process can be done, via in-splittings to vertices in $B$ (and some additional in- and out-splittings), to enlarge the length
   	of $\alpha $ by any positive multiple of $|B|$, without modifying the other trails from $A$ to $B$,
   	obtaining a new graph $E'$ in quasi-normal form. Combining both processes, we can replace $\alpha$ by
 another trail $\alpha'$ such that 
 $$|\alpha' | = |\alpha| + a|A| + b|B|$$
 for any non-negative integers $a,b$. Moreover, this process will not change the other trails from $C$ to $C'$, for any two cycles $C,C'$ such that $d(C,C')\le d(A,B)$.   	
   	\end{proof}

 We are now in a position to present the main result of this section which shows that every connected finite graph with disjoint cycles can be transformed into a graph in normal form.
 
\begin{theorem}\label{thm:normal-form}
Any connected finite graph with disjoint cycles $E$ can be transformed to a graph in normal form through a sequence of in-splittings and
out-splittings. 
   \end{theorem}  
   \begin{proof}
By Theorem \ref{thm:quasi-normal-form}, it suffices to assume that $E$ is in quasi-normal form.   
	Consider two cycles $A$ and $B$ in $E$ such that there are trails from $A$ to $B$. Let $d= \gcd (|A|,|B|)$. By Lemma \ref{lem:number-theory}, we can find $d$ large positive integers $n_0,n_1,\ldots ,n_{d-1}$, where $n_i \equiv i \pmod{|A| |B|} $ for all $0 \le i < d$, and such that for each trail $\alpha $ from $A$ to $B$ such that $|\alpha | \equiv i \pmod{d}$, there exist non-negative integers $a_{\alpha},b_{\alpha}$ such that
 $$n_i = |\alpha | + a_{\alpha} |A| + b_{\alpha} |B|.$$
 
  Now, by Lemma \ref{lem:basic-for-normal-form}, for each $0\le i < d$ and each trail $\alpha$ from $A$ to $B$ such that $|\alpha | \equiv i \pmod{d}$, we can enlarge $\alpha$ to a new trail from $A$ to $B$ of length exactly $n_i$, while every other trail between cycles $C,C'$ with $d(C,C')\leq d(A,B)$ is unchanged. Doing this process for every trail from $A$ to $B$, we obtain a new graph $E'$, which is still in quasi-normal form, such that for any two trails $\alpha$, $\beta$ from $A$ to $B$ such that $|\alpha | \equiv |\beta| \pmod{d}$ it happens that $|\alpha |= |\beta|$. In addition, we have that $|\alpha| -t$ is a multiple of $|A| |B|$ whenever $|\alpha| \equiv t \pmod{d} $. Moreover all the trails from a cycle $C$ to a cycle $C'$, with $d(C,C')\le d(A,B)$ and $(C,C')\ne (A,B)$, will remain unchanged after the transition to $E'$. In particular, it follows that in the transformed graph $E'$, the pair of cycles $(A,B)$ is good. 
 
 Now we proceed by induction to show the result. Suppose that $i\ge 0$ and that $E^{(i)}$ is a graph in quasi-normal form obtained from $E$ by a sequence of out- and in-splittings, such that all pairs of cycles $(A,B)$ such that $d(A,B)\le i$ are good. Then performing the above process to the family of all pairs $(C,D)$ such that $d(C,D) = i+1$, we obtain a graph $E^{(i+1)}$, which is still in quasi-normal form, and such that all the pairs of cycles $(C',C'')$ such that $d(C',C'')\le i+1$ are good. 
 
 This completes the induction step, and the proof of the theorem.  \end{proof}


\section{Dynamics of meteor graphs of length three}\label{sec4}
In this section, we introduce the notion of meteor graphs of length three. Based on Theorem \ref{willimove} and Theorem \ref{thm:normal-form}, we then provide number-theoretic criteria for meteor graphs of length three to be strongly shift equivalent (Theorem \ref{numtheo}).\medskip

We begin this section by introducing the definition of meteor graphs of length three.

\begin{deff}\label{def31}
A \textit{meteor graph of length three}  is a connected essential graph $E$ consisting of three disjoint cycles $P_E,Q_E$ and $R_E$ with a unique chain of cycles $P_E > Q_E > R_E$. 
We call $P_E$ the first cycle, $Q_E$ the second cycle and $R_E$ the third cycle. We denote the vertices in $P_E$ by $p_1,p_2, \ldots, p_{|P_E|}$, where $p_i$ is the vertex of index $i$ in $P_E$, the vertices in $Q_E$ by $q_1,q_2,\ldots, q_{|Q_E|}$ and the vertices in $R_E$ by $r_1,r_2,\ldots, r_{|R_E|}$.  
\end{deff}

For example, the graph pictured in Example \ref{exa3.4} is a meteor graph of length three.


The following proposition shows that shift equivalence preserves the cycle structure of meteor graphs of length three.

\begin{prop}\label{propreserve}
Let $E$ and $F$ be essential graphs such that $A_E\sim_{SE} A_F$. Then, if $E$ is a meteor graph of length three with a unique chain of cycles $P_E > Q_E > R_E$, then  $F$ is also a meteor graph of length three with a unique chain of cycles  $P_F > Q_F > R_F$ such that $|P_E| = |P_F|$, $|Q_E| = |Q_F|$ and $|R_E| = |R_F|$.
\end{prop}\label{propo33}
\begin{proof}
It immediately follows from Proposition \ref{prop:SSE-implies-same structure}.   
\end{proof}

Let $E$ be a meteor graph of length three with a unique chain of cycles $P_E > Q_E > R_E$ such that $|P_E| = x^E$,  $|Q_E| = y^E$, and  $|R_E| = z^E$. We denote
$$
d^E = \gcd \left( x^E,z^E \right), d_1^E=\operatorname{gcd}\left(x^E, y^E\right) \text{ and }  d_2^E = \gcd (y^E, z^E) . 
$$
Let $N_1^E(c)$ be the number of trails from $P_E$ to $Q_E$ which have length congruent to $c$ modulo $d_1^E$, for all  $1 \le c \le d_1^E$, $N_2^E(c)$ be the number of trails from $Q_E$ to $R_E$ which have length congruent to $c$ modulo $d{_2^E}$, for all $1 \le c \le d_2^E$, and $N_E(c)$ be the number of trails from $P_E$ to $R_E$ which have length congruent to $c$ modulo $d^E$, for all $1 \le c \le d^E$.  Let $l_{P,Q}^E$ be the number of all trails from $P_E$ to $Q_E$ in $E$, $l_{Q,R}^E$ the number of all trails from $Q_E$ to $R_E$ in $E$, and $l_{P,R}^E$ the number of all trails from $P_E$ to $R_E$ in $E$. We denote

\begin{center}
$ S^E_1(a) =  \sum \limits_{\begin{subarray}{c}
				1 \le k \le \gcd (x^E,y^E) \\
				k \equiv a \pmod{ \gcd (x^E,y^E,z^E)}
		\end{subarray}} N^E_1(k)\,\,$  for all  $\,\, 1 \le a \le  \gcd (x^E,y^E,z^E)$     
\end{center}
and 

\begin{center}
$S^E_2(b) =  \sum \limits_{ \begin{subarray}{c}
				1 \le k \le \gcd (y^E,z^E) \\
				k \equiv b \pmod{ \gcd (x^E,y^E,z^E)}
		\end{subarray} } N^E_2(k) \,\, $ for all  $\,\, 1 \le b \le \gcd(x^E,y^E,z^E).$    
\end{center}

The following key definition plays an important role in our analysis. In the following, we will often use the notation $(x,y)$ for the greatest common divisor $\gcd (x,y)$ of two positive integers $x,y$,
and we will use $[x,y]$ for the least common multiple $\text{lcm} (x,y)$ of $x$ and $y$. Note the equation $(x,y)[x,y]= xy$.

\begin{deff}\label{num-cre-deff} Let $E$ and $F$ be two meteor graphs of length three in normal form. We write $E \approx F$ if the following conditions are satisfied:
	\begin{enumerate}[label=(\roman*)]
		\item $x^E=x^F:=x$, $y^E = y^F:= y$, $z^E = z^F:=z$ 
		\item $N_1^E(c) = N_1^F(c): =N_1(c)  \text{ for all } 1 \le c \le \gcd (x,y) $ 
		\item $N_2^E(c) = N_2^F(c): = N_2(c)   \text{ for all } 1 \le c \le \gcd (y,z) $
		\item There exist integers $g_1,g_2,\ldots ,g_{ \gcd (x,y,z) }  $ and $ l_1,l_2,\ldots ,l_{ \gcd (x,y,z) }   $  such that the following conditions hold:
		
		(1) $ [x,y] \mid l_a $ for all $ 1\le a \le  \gcd (x,y,z)  $, 
		
		(2) $l_a =0 $ if $S_1 (a) = 0  $ for all $1 \le a \le \gcd  (x,y,z) $, 
		
		(3) $ [y,z] \mid g_b $ for all $ 1\le b \le  \gcd (x,y,z)$, 
		
		(4) $g_b = 0 $ if $S_2(b) = 0 $ for all $ 1  \le b \le \gcd  (x,y,z) $, 
		
		(5) $N_F(h)-N_E(h) = \dfrac{ \gcd (x,y,z)}{y \gcd (x,z)} \left(\sum \limits_{\begin{subarray}{c}
				1 \le a,b \le (x,y,z)  \\ 
				a+b \equiv h \pmod{(x,y,z)}
		\end{subarray}} S_1(a) g_b  + \sum \limits_{\begin{subarray}{c}
				1 \le a,b \le (x,y,z)  \\ 
				a+b \equiv h \pmod{(x,y,z)}
		\end{subarray}} l_a S_2(b)\right)$ for all $1\le h\le \gcd (x,z)$,	 
	\end{enumerate}
where $S_1(a):= S_1^E(a) = S_1^F(a)$ and  $S_2(b):= S_2^E(b) = S_2^F(b)$.   
\end{deff}

Note that if $(x,y,z)=1$, then we have $S_1(1)= N_1(1)\ne 0$ and   $S_2(1)= N_2(1) \ne 0$ (since both $E$ and $F$ are meteor graphs of length three). In this case, the conditions $(2)$ and $(4)$ of Definition \ref{num-cre-deff} (iv) are always (vacuously) satisfied. 

It is an easy exercise to check that $\approx$ is an equivalence relation on the set of meteor graphs of length three in normal form. 

In the remainder of this section, we  prove the main result of this section (Theorem \ref{numtheo}) that $E \approx F$ if and only if we can transform $E$ to $F$ by a finite sequence of in-splittings, out-splittings,  in-amalgamations, and out-amalgamations. To establish this result, we first need to control the emergence of new trails from the first cycle to the last cycle during the graph transformations. The following three lemmas provide the necessary tools to address this issue.

\begin{lemma}
\label{lem:Trail-inequality}
Let $E$ be a meteor graph of length three with a unique chain of cycles $P_E > Q_E > R_E$, and let $G$ be the graph obtained from $E$ by either a non-trivial in-splitting or a non-trivial out-splitting at  a vertex $v$ in $E$. Let $\textnormal{Trail}_E(P_E,R_E)$ and $\textnormal{Trail}_G(P_G,R_G)$ be the sets of trails from $P_E$ to $R_E$ and from $P_G$ to $R_G$, respectively. Then, $|\textnormal{Trail}_E(P_E,R_E)| < |\textnormal{Trail}_G(P_G,R_G)|$ if and only if $v\in Q_E^0$, $|s^{-1}(v)| \ge 2$ and $|r^{-1}(v)|\ge 2$. 
\end{lemma}
\begin{proof}
Write $r^{-1}(v) = \{a_1,a_2,\ldots ,a_m\}$ and $s^{-1}(v)= \{b_1,b_2,\ldots, b_n\}$. If $v \in Q^0_E$, we denote by $\text{Trail}^*_E(P_E,R_E)$ the set of all paths $e_1e_2\cdots e_n$ in $E$ satisfying the following conditions:
\begin{enumerate}
\item $s(e_1) \in P_E,$ 
\item $r(e_n) \in R_E,$ 
\item $s(e_2),s(e_3),\ldots , s(e_n) \notin (P_E^0\cup Q_E^0 \cup R_E^0)\setminus \{v\}$, and 
\item there exists a unique number $1 \le j < n$ such that $s(e_{j+1})= r(e_j) = v$.   
\end{enumerate}

We first consider the case of performing an out-splitting at the vertex $v\in E^0$ with a partition $\mathcal{E}_1 \cup  \mathcal{E}_2 \cup  \cdots  \cup\mathcal{E}_h = \{ b_1,b_2, \ldots ,b_n\} $, with $h >1$, so that $n>1$. The resulting graph $G$ has new vertices $v_1,v_2,\ldots ,v_h$ and new edges $\{a_{i,j} \mid 1 \le i \le m, 1 \le j \le h \}$ replacing  $a_1,a_2,\ldots ,a_m$. By the definition of out-splittings, in $G$ we have  $s(b_i)= v_j $ if $b_i \in \mathcal{E}_j$ for all $1 \le i \le n$, and $r(a_{i,j})=v_j$ for all $1 \le i \le m$ and $ 1\le j \le h$. 

If $v \in Q^0_E$, without loss of generality, we may assume that $b_n$ is the unique edge in the cycle $Q_E$ whose source is $v$, and that $a_m$ is the unique edge in $Q_E$ with range $v$. Moreover, we may assume that $b_n\in \mathcal{E}_h$, so that $v_h \in Q^0_G$ and $v_i\notin Q^0_G$ for $1\le i<h$.

Assume first that $m>1$. Note that $\text{Trail}^*_E(P_E,R_E)$ is the set of all paths $\alpha a_i b_j \beta$, where $\alpha a_i \in \text{Trail}_E(P_E,Q_E)$ and $b_j \beta \in \text{Trail}_E(Q_E,R_E)$. In particular 
$i<m$ and $j<n$ whenever $\alpha a_ib_j\beta \in \text{Trail}^*_E(P_E,R_E)$. 

Define
$$  \text{Trail}^{**}_E(P_E,R_E)  = \{\alpha a_i b_j \beta \in \text{Trail}^* (P_E,R_E) \mid b_j\in \mathcal{E}_t \text{ with } 1\le t<h \}.$$
We have a bijection 
$$ \Phi \colon \text{Trail} (P_E,R_E) \sqcup \text{Trail}^{**}_E(P_E,R_E) \longrightarrow \text{Trail}_G(P_G,R_G)$$
given by $\Phi (\gamma) = \gamma$ if $\gamma \in \text{Trail}_E(P_E,R_E)$, and 
$\Phi (\alpha a_i b_j \beta) = \alpha a_{i,t} b_j\beta $, where $1\le t\le h-1$  is the unique integer such that $b_j\in \mathcal{E}_t$ if $\alpha a_i b_j\beta \in  \text{Trail}^{**}(P_E,R_E)$. 
Since $m>1$ and $h>1$, it follows that $\text{Trail}^{**}_E(P_E,R_E)\ne \emptyset$, and hence $|\textnormal{Trail}_E(P_E,R_E)| < |\textnormal{Trail}_G(P_G,R_G)|$.

Whenever $m=1$ we get that $\text{Trail}^*_E(P_E,R_E) = \text{Trail}^{**}_E(P_E,R_E)= \emptyset$, and so  $$|\textnormal{Trail}_E(P_E,R_E)| = |\textnormal{Trail}_G(P_G,R_G)|.$$

A similar analysis can be done in the case where we perform an in-splitting at $v\in Q^0_E$. 

When $v\notin Q^0_E$, we clearly have that  $|\textnormal{Trail}_E(P_E,R_E)| = |\textnormal{Trail}_G(P_G,R_G)|$. Indeed, in this case, adopting the above notation, we have
a bijection
$$\Phi \colon  \text{Trail}_E (P_E,R_E) \longrightarrow \text{Trail}_G(P_G,R_G)$$
defined by $\Phi (\alpha a_ib_j\beta) = \alpha a_{it} b_j\beta$, where $t\in \{1,\dots ,h\}$ is the unique index such that $b_j\in \mathcal E_t$, for all $1\le i\le m$, $1\le j\le n$. . 

\end{proof}

\begin{lemma}\label{prop-functor}
Let $E$ be a meteor graph of length three with a unique chain of cycles $P_E > Q_E > R_E $. Let  $\textnormal{Trail}_E(P_E, Q_E)$,
$\textnormal{Trail}_E(Q_E, R_E)$, and $\textnormal{Trail}_E(P_E, R_E)$ denote the sets of all trails in $E$ connecting $P_E$ to $Q_E$,  $Q_E$ to $R_E$, and $P_E$ to $Q_E$, respectively. For any $\alpha \in \textnormal{Trail}_E(P_E, Q_E)\cup\textnormal{Trail}_E(Q_E, R_E) \cup \textnormal{Trail}_E (P_E,R_E)$ whose source has index $a$ and whose range has index $b$, define \[f_E(\alpha) = \begin{cases}
    b - (a+ |\alpha|) \pmod{(x, y)} \mbox{ if } \alpha \in \textnormal{Trail}_E(P_E, Q_E), \\
     b - (a+ |\alpha|) \pmod{(y, z)} \mbox{ if } \alpha \in \textnormal{Trail}_E(Q_E, R_E), \\
       b - (a+ |\alpha|) \pmod{(x, z)} \mbox{ if } \alpha \in \textnormal{Trail}_E(P_E, R_E),
\end{cases}\] where $x= |P_E|$, $y=|Q_E|$, and $z=|R_E|$. Let $G$ be a graph obtained from $E$ by an in-splitting or an out-splitting at a vertex lying on a cycle, or by a graph move at an interior vertex.  Then there exist  bijections \[
     \Phi^1 _{E,G}: \textnormal{Trail}_E (P_E,Q_E) \to \textnormal{Trail}_G(P_G, Q_G),
    \]    
      \[
     \Phi^2_{E,G}: \textnormal{Trail}_E (Q_E,R_E) \to \textnormal{Trail}_G (Q_G,R_G),
    \] 
    and an injection \[ \Phi^3_{E,G}: \textnormal{Trail}_E (P_E,R_E) \to \textnormal{Trail}_G (P_G,R_G),
    \] such that 
\[f_E(\alpha) = \begin{cases}
	f_G(\Phi^1_{E,G}(\alpha))  \mbox{ if } \alpha \in \textnormal{Trail}_E(P_E, Q_E), \\
   f_G(\Phi^2_{E,G}(\alpha)) \mbox{ if } \alpha \in \textnormal{Trail}_E(Q_E, R_E),\\
f_G(\Phi^3_{E,G}(\alpha)) \mbox{ if } \alpha \in \textnormal{Trail}_E(P_E, R_E).\end{cases}\]     
\end{lemma}	
\begin{proof}
Consider the case where $G$ is obtained from $E$ by performing either an in-splitting or an out-splitting. Let $\alpha$ be an arbitrary element of $\textnormal{Trail}_E(P_E, Q_E)$ whose source has index $a$ on $P_E$ and whose range has index $b$ on $Q_E$. Let $v$ be an arbitrary vertex in $E$, and let $G$ be the graph obtained from $E$ by either an in-splitting or an out-splitting at $v$. If $v\notin \alpha^0$, then $\alpha$ remains unchanged, and we define $$\Phi^1_{E, G}(\alpha) = \alpha\in \textnormal{Trail}_G(P_G, Q_G).$$ In this case, it is obvious that  $$f_E(\alpha) = f_G(\Phi^1_{E, G}(\alpha)).$$ We now consider the remaining cases.

{\it Case} $1$: $v\in \alpha^0\setminus\{s(\alpha), r(\alpha)\}$, and $G$ is obtained from $E$ by an out-splitting  at $v$ with a partition $ \{\mathcal{E}^v_1,\ldots, \mathcal{E}^v_t\} $ of $s^{-1}(v)$. We then have $$G^0 = E^0\setminus\{v\} \cup \{v^1, v^2, \ldots, v^t\}$$ and $$G^1 = E^1\setminus r^{-1}(v) \cup \{e^1, \ldots, e^t\mid e\in r^{-1}(v)\}.$$ Write $\alpha = \alpha_1eg\alpha_2$, where $\alpha_1, \alpha_2\in E^*$, $e, g\in E^1$, $r(\alpha_1)= s(e)$, $r(e) = s(g) =v$, and $r(g) = s(\alpha_2)$. Let $k$ (with $1\le k\le t$) be the unique integer such that $g\in \mathcal{E}^v_k$. We define $$\Phi^1_{E, G}(\alpha) = \alpha_1e^kg\alpha_2\in \textnormal{Trail}_G(P_G, Q_G)).$$ In this case, we note that   $$|\Phi^1_{E, G}(\alpha)| = |\alpha_1e^kg\alpha_2|=|\alpha|,$$ and so
it is clear that  $$f_E(\alpha) = f_G(\Phi^1_{E, G}(\alpha)).$$ 

{\it Case} $2$: $v= s(\alpha)$ and $G$ is obtained from $E$ by an out-splitting  at $v$ with a partition $ \{\mathcal{E}^v_1,\ldots, \mathcal{E}^v_t\} $ of $s^{-1}(v)$. There is a unique edge $e\in E^1$ such that $r(e) = v$ (note that $e$ lies on cycle $P_E$). Then 
$$G^0 = E^0\setminus\{v\} \cup \{v^1, v^2, \ldots, v^t\}$$ and $$G^1 = E^1\setminus \{e\} \cup \{e^1, \ldots, e^t\}.$$ Let $f$ be the edge on $P_E$ with $s(f) = v$. We assume, without loss of generality, that $f\in \mathcal E^v_t$. Write $\alpha = g\beta$, where $g\in E^1$, $\beta\in E^*$, $s(g) = v$, and $r(g) = s(\beta).$ Let $l$ ($1\le l\le t$) be the unique integer such that $g\in \mathcal{E}^v_l$. 

\textit{Case 2.1:} If $1 \le l < t$ We define $$\Phi^1_{E, G}(\alpha) = 
  e^lg\beta  .
$$
We can easily see that $\Phi^1_{E, G}(\alpha)\in \textnormal{Trail}_G(P_G, Q_G)$.  Then 
 $$|\Phi^1_{E, G}(\alpha)| = |e_lg\beta| = |\alpha| +1,\  r_G(\Phi^1_{E, G}(\alpha))=r_G(e^lg\beta) = r_E(\alpha)$$ and $$s_G(\Phi^1_{E, G}(\alpha)) = s_G(e^lg\beta)$$ is  the vertex with index $a-1 \pmod{x}$. Therefore,
\begin{align*}
f_G(\Phi^1_{E, G}(\alpha))&=b - (a -1 + |\alpha| +1) \pmod{(x, y)}\\
&= b - (a + |\alpha|) \pmod{(x, y)}\\
& = f_E(\alpha).
\end{align*}
$$\xymatrixrowsep{0.5pc}\xymatrixcolsep{0.5pc}\xymatrix{&&\bullet\ar@/_0.8pc/@{-->}[ddll]&&&&&&...&&&&&&&&\bullet\ar@/_0.8pc/@{-->}[ddll]&&&&&&...&&&&&&\\&&&&&&\bullet\ar@/_0.7pc/[urr]&&&&&&&&&&&&&&\bullet\ar@/_0.7pc/[urr]&&&&&&&&\\...\ar@/_0.8pc/@{-->}[ddrr]&&&& \bullet^v\ar@/_0.8pc/[uull]^f\ar@/^0.8pc/[urr]^g&&&&&\ar@{->}[rrrr]&&&&&...\ar@/_0.8pc/@{-->}[ddrr]&&&& \bullet^{v^h}\ar@/_0.8pc/[uull]^f&&&&&\\&&&&&&&&&&&&&&&&&&&\bullet^{v^l}\ar@/^1.4pc/[uur]^g&&&&&&&&&\\&&\bullet^u\ar@/_0.8pc/[uurr]^e&&&&&&&&&&&&&&\bullet^u\ar@/_0.8pc/[uurr]^{e^h}\ar@/_1.5pc/[urrr]^{e^l}&&&&&&&&&&&&}$$
\begin{figure}[h]
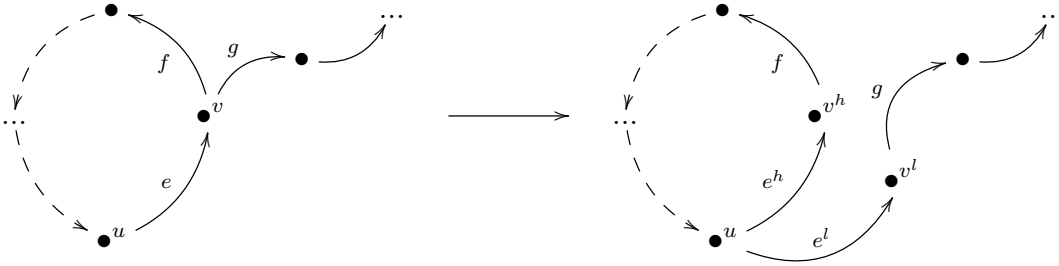
\caption{Out-splitting at $v$.}
\end{figure} 

\textit{Case 2.2:} If $l=t.$ We define 
 $$\Phi^1_{E, G}(\alpha) = 
  g\beta  = \alpha   .
$$
Then we can easily see that $f_G( \Phi^1_{E, G}(\alpha) ) = f_E(\alpha).$


{\it Case} $3$: $v\in \alpha^0\setminus\{s(\alpha), r(\alpha)\}$ and $G$ is obtained from $E$ by an in-splitting  at $v$ with a partition $ \{\mathcal{E}^v_1,\ldots, \mathcal{E}^v_t\} $ of $r^{-1}(v)$.  We then have $$G^0 = E^0\setminus\{v\} \cup \{v^1, v^2, \ldots, v^t\}$$ and $$G^1 = E^1\setminus s^{-1}(v) \cup \{g^1, \ldots, g^t\mid g\in s^{-1}(v)\}.$$ Write $\alpha = \alpha_1eg\alpha_2$, where $\alpha_1, \alpha_2\in E^*$, $e, g\in E^1$, $r(\alpha_1)= s(e)$, $r(e) = s(g) =v$, and $r(g) = s(\alpha_2)$. Let $k$ ($1\le k\le t$) be the unique integer such that $e\in \mathcal{E}^v_k$. We define $$\Phi^1_{E, G}(\alpha) = \alpha_1eg^k\alpha_2\in \textnormal{Trail}_G(P_G, Q_G).$$ In this case, we note that $$|\Phi^1_{E, G}(\alpha)| = |\alpha_1eg^k\alpha_2|=|\alpha|,$$ and hence
it is clear that  $$f_E(\alpha) = f_G(\Phi^1_{E, G}(\alpha)).$$

{\it Case} $4$: $v= r(\alpha)$ and $G$ is obtained from $E$ by an in-splitting  at $v$ with a partition $ \{\mathcal{E}^v_1,\ldots, \mathcal{E}^v_t\} $ of $r^{-1}(v)$. There is a unique edge $f$ lying on cycle $Q_E$ such that $r(f) = v$. We assume, without loss of generality, that $f\in \mathcal E^v_t$. Then 
$$G^0 = E^0\setminus\{v\} \cup \{v^1, v^2, \ldots, v^t\}$$ and $$G^1 = E^1\setminus s^{-1}(v) \cup \{e^1, \ldots, e^t\mid e\in s^{-1}(v)\}.$$ 
Let $g$ be the edge on $Q_1^E$ with $s(g) = v$. Write $\alpha = \beta e$, where $e\in E^1$, $\beta\in E^*$, $r(e) = v$, and $s(e) = r(\beta).$ Let $l$ ($1\le l\le t$) be the unique integer such that $e\in \mathcal{E}^v_l$.

\textit{Case 4.1:} If $1 \le l <t$.  We define $$\Phi^1_{E, G}(\alpha) = 
\beta e g^l, \quad   1\le l<t$$
We can easily see that $\Phi^1_{E, G}(\alpha)\in \textnormal{Trail}_G(P_G, Q_G)$. 
We then have $|\Phi^1_{E, G}(\alpha)| = |\beta eg^l| = |\alpha| +1$, $s_G(\Phi^1_{E, G}(\alpha))=s_G(\beta eg^l) = s_E(\alpha)$ and $r_G(\Phi^1_{E, G}(\alpha)) = r_G(\beta e g^l)$ is  the vertex with index $b +1 \pmod{y }$. Therefore,
\begin{align*}
f_G(\Phi^1_{E, G}(\alpha))&=b + 1- (a + |\alpha| +1) \pmod{(x, y)}\\
&= b - (a + |\alpha|) \pmod{(x, y)}\\
& = f_E(\alpha).
\end{align*}
$$\xymatrixrowsep{0.5pc}\xymatrixcolsep{0.5pc}\xymatrix{...\ar@/_0.8pc/[drr]&&&&&&\bullet\ar@/_0.8pc/[ddll]_f&&&&&&&&...\ar@/_0.8pc/[drr]&&&&&&\bullet\ar@/_0.8pc/[ddll]_f&&\\&&\bullet\ar@/^0.8pc/[drr]^e&&&&&&&&&&&&&&\bullet\ar@/^1.3pc/[ddr]^e&&&&&&\\&&&&\bullet_{v}\ar@/_0.8pc/[ddrr]_g&&&&...\ar@/_0.8pc/@{-->}[uull]&\ar@{->}[rrrr]&&&&&&&&&\bullet_{v^h}\ar@/_0.8pc/[ddrr]^{g^h}&&&&...\ar@/_0.8pc/@{-->}[uull]\\&&&&&&&&&&&&&&&&&\bullet^{v^l}\ar@/_1.5pc/[drrr]_{g^l}&&&&&\\&&&&&&\bullet^u\ar@/_0.8pc/@{-->}[uurr]&&&&&&&&&&&&&&\bullet^u\ar@/_0.8pc/@{-->}[uurr]&&}$$
\begin{figure}[h]
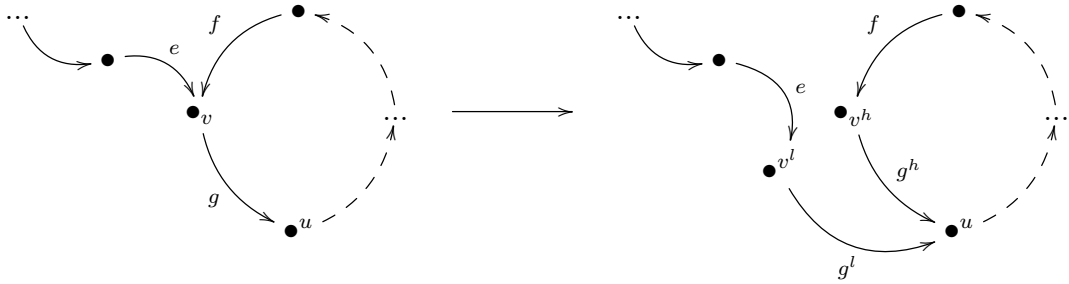
\caption{In-splitting at $v$.}
\end{figure} 

\textit{Case 4.2:} If $l=t$. We define 
$$\Phi^1_{E, G}(\alpha) = 
  \beta e = \alpha  $$
 We can easily see that $\Phi_{E,G}^1 (\alpha) \in \text{Trail}_G(P_G,Q_G)$ so $f_G(\Phi_{E,G}^1 (\alpha) ) = f_E(\alpha).$


Therefore, we have defined a bijection $\Phi^1_{E,G}$ such that 
\begin{center}
$f_E(\alpha) = f_G(\Phi^1_{E,G} (\alpha))$ for every $\alpha \in \text{Trail}_E(P_E, Q_E)$.
\end{center}
By a similar argument, we can define a bijection $\Phi^2_{E,G}$ and an injection $\Phi^3  _{E,G}$ such that 
\begin{center}
$f_E(\alpha) = f_G(\Phi^2_{E,G} (\alpha))$ for every $\alpha\in \text{Trail}_E(Q_E, R_E)$    
\end{center}
 and 
\begin{center}
$ f_G(\Phi^3_{E,G} (\alpha) ) = f_E(\alpha)$ for every $\alpha \in \text{Trail}_E (P_E,R_E)$.  
\end{center} 
Note that the injection $\Phi^3_{E,G}$ was constructed in the course of the proof of Lemma \ref{lem:Trail-inequality}. Moreover, if an out-splitting or an in-splitting is performed at an interior vertex, then the maps $\Phi_{E,G}^1,\Phi_{E,G}^2,\Phi_{E,G}^3$ are bijections. Consequently, the lemma holds when $G$ is obtained from $E$ by a graph move at an interior vertex, thus finishing the proof.
\end{proof}	

We note that the proof of Lemma \ref{prop-functor} yields the following. For each $\alpha \in \text{Trail}_E(P_E, Q_E)$, if an in-splitting is performed at $r(\alpha)$ to obtain $F$ and $\Phi_{E,F}^1(\alpha) \neq \alpha$, then $|\Phi_{E,F}^1(\alpha)| = |\alpha| + 1$, while $s(\Phi_{E,F}^1(\alpha)) \equiv s(\alpha) \pmod{x}$. 

Similarly, for each $\beta \in \text{Trail}_E(Q_E, R_E)$, if an in-splitting is performed at $s(\beta)$ to obtain $F$ and $\Phi_{E,F}^2(\beta) \ne \beta$, then $|\Phi_{E,F}^2(\beta)| = |\beta| + 1$, while $r(\Phi_{E,F}^2(\beta)) \equiv r(\beta) \pmod{z}$. 

Combining this observation with Lemma \ref{prop-functor}, we obtain the following important result.
    
\begin{lemma}\label{functor}
Let $E$ be a meteor graph of length three with a unique chain of cycles $P_E > Q_E > R_E $. Let  $\textnormal{Trail}_E(P_E, Q_E)$  denote the set of all trails of $E$ connecting $P_E$ to $Q_E$, and let $\textnormal{Trail}_E(Q_E, R_E)$ denote the set of all trails of $E$ connecting $Q_E$ to $R_E$. For any $\alpha \in \textnormal{Trail}_E(P_E, Q_E)$ whose source has index $a$, and any $\beta \in \textnormal{Trail}_E(Q_E, R_E)$  whose range has index $b$, define \[f_E(\alpha,\beta) = 
    b - (a+ |\alpha| + |\beta| ) \pmod{(x, z)}   ,
\] where $x= |P_E|,y=|Q_E|$, and $z=|R_E|$. Let $G$ be a graph obtained from $E$ by an in-splitting or an out-splitting at a vertex lying on a cycle, or by a graph move at an interior vertex.   Then there exist bijections \[
     \Phi^1 _{E,G}: \textnormal{Trail}_E (P_E,Q_E) \to \textnormal{Trail}_G(P_G, Q_G) 
    \]    
    and  \[
     \Phi^2_{E,G}: \textnormal{Trail}_E (Q_E,R_E) \to \textnormal{Trail}_G (Q_G,R_G)
    \]  such that, for all $\alpha\in \textnormal{Trail}_E (P_E,Q_E)$ and $\beta\in \textnormal{Trail}_E (Q_E,R_E)$,
\[f_E(\alpha,\beta) = \begin{cases}
	f_G(\Phi^1_{E,G}(\alpha), \Phi^2_{E,G}(\beta))+1\pmod{(x,z)}\text{ if we in-split at $r(\alpha)$ and $\Phi_{E,G}^1(\alpha) \ne \alpha,$}  \\
  f_G(\Phi^1_{E,G}(\alpha), \Phi^2_{E,G}(\beta))+1\pmod{(x,z)}\text{ if we out-split at $s(\beta)$ and  $\Phi_{E,G}^2(\beta) \ne \beta$,}  \\  
f_G(\Phi^1_{E,G}(\alpha), \Phi^2_{E,G}(\beta))  \pmod{(x,z)}\text{\, \, \,  otherwise.} \end{cases}\]     
\end{lemma}	
\begin{proof}
To prove the lemma, it suffices to consider the case where $G$ is obtained from $E$ by an in-splitting or an out-splitting. Let $\Phi^1_{E,G}$ and $ \Phi_{E,G}^2$ denote the bijections defined in Lemma \ref{prop-functor}. We now compare  $f_G(\Phi^1_{E,G}(\alpha), \Phi ^2 _{E,G} (\beta))$ and   $f_E(\alpha,\beta)$ by considering the following cases.

 \textit{Case 1:} Out-splitting at $ s(\alpha)$.

 \textit{Case 1.1:} $\Phi_{E,G} ^1(\alpha) \ne \alpha$.  Then $\Phi_{E,G} ^1(\alpha)$  has its length increase by $1$, and its source has index $a-1 \pmod{x}$. Therefore, \[f_G(\Phi^1_{E,G}(\alpha), \Phi^2_{E,G}(\beta)) = f_E(\alpha,\beta) \pmod{ (x,z) } . \]

 \textit{Case 1.2:} $\Phi_{E,G} ^1(\alpha) = \alpha$. Then we can easily see that \[f_G(\Phi^1_{E,G}(\alpha), \Phi^2_{E,G}(\beta)) = f_E(\alpha,\beta) \pmod{ (x,z) }.\] 
 
  \textit{Case 2:} In-splitting at $ r(\beta)$. 

\textit{Case 2.1:} $\Phi_{E,G}^2(\beta) \ne \beta$.
  Then $\Phi_{E,G} ^2(\beta)$ have its length increase by $1$, and its range has index $b+1 \pmod{z}$, and so  \[f_G(\Phi^1_{E,G}(\alpha), \Phi^2_{E,G}(\beta)) = f_E(\alpha,\beta) \pmod{ (x,z) }.\]  

\textit{Case 2.2:} $\Phi_{E,G}^2(\beta) = \beta$. Then we can easily see that \[f_G(\Phi^1_{E,G}(\alpha), \Phi^2_{E,G}(\beta)) = f_E(\alpha,\beta) \pmod{ (x,z) }.\]

  \textit{Case 3:} Using a graph move at a vertex $v \notin \{s(\alpha),r(\alpha),s(\beta),r(\beta) \}$. Then $\Phi^1_{E,G}(\alpha)$ has the same length, source index, and range index as $\alpha$. The same holds for $\Phi^2_{E,G} (\beta)$ and $\beta$. This implies that \[f_G(\Phi^1_{E,G}(\alpha), \Phi^2_{E,G}(\beta)) = f_E(\alpha,\beta) \pmod{ (x,z) }.\] 
  
  \textit{Case 4:} In-splitting at  $r(\alpha)$. 
  
  \textit{Case 4.1:} $\Phi_{E,G}^1 (\alpha) \ne \alpha.$ Then $\Phi_{E,G} ^1(\alpha)$ have its length increase by $1$, while its source retains index $a$. Therefore, \[f_G(\Phi^1_{E,G}(\alpha), \Phi^2_{E,G}(\beta)) = f_E(\alpha,\beta) - 1 \pmod{ (x,z) }.\] 

 \textit{Case 4.2:} $\Phi_{E,G}^1 (\alpha)  = \alpha.$ Then we can easily see that \[f_G(\Phi^1_{E,G}(\alpha), \Phi^2_{E,G}(\beta)) = f_E(\alpha,\beta) \pmod{ (x,z) }.\]   

  \textit{Case 5:} Out-splitting at  $s(\beta)$. 
  
  \textit{Case 5.1:}  $\Phi_{E,G}^2(\beta) \ne  \beta$.  Then $\Phi_{E,G}^2 (\beta)$ have its length increase by $1$, while its range retains index $b$, and so  \[f_G(\Phi^1_{E,G}(\alpha), \Phi^2_{E,G}(\beta)) = f_E(\alpha,\beta) - 1 \pmod{ (x,z) }.\] 

  \textit{Case 5.2:} $\Phi_{E,G}^2(\beta) =   \beta$. Then we can easily see that \[f_G(\Phi^1_{E,G}(\alpha), \Phi^2_{E,G}(\beta)) = f_E(\alpha,\beta) \pmod{ (x,z) }.\]
  
Therefore, in any case, we obtain the statement, thus finishing the proof.
\end{proof}

Using Lemmas \ref{lem:Trail-inequality} and \ref{prop-functor}, we describe a method for numbering the new trails from the first cycle to the last cycle after applying graph moves. We first give a rule for numbering the trails from the first cycle to the middle cycle, and another rule for numbering the trails from the middle cycle to the last cycle. It is important to take into account that the cycle structure of the graph keeps intact along all the graph moves, and the changes are performed only on the trails of the graph. Hence we always assume that we are working with a meteor graph of length three with normalized set of disjoint cycles (see Section \ref{sec3} for details).

\begin{rem}\label{rem10}
Let $E$ be a meteor graph of length three and $G$ the graph obtained from $E$ by applying either an in-splitting or an out-splitting at  a vertex $v$ in $E$. Let $l=|\textnormal{Trail}_E(P_E,Q_E)|$. Number the trails from $P_E$ to $Q_E$ by $1,2,\ldots ,l$. If $\alpha$ has number $c$ in $E$, then $\Phi^1(\alpha)$ is assigned the same number $c$ in $G$. A similar rule applies to the trails from $Q_E$ to $R_E$. 
\end{rem}

We now provide a rule to number the new trails from the first cycle to the last cycle.

\begin{rem}\label{rem11}
Let $E$ be a meteor graph of length three, and let $G$ be the graph obtained from $E$ by applying either an in-splitting or an out-splitting at a vertex $v$ in $E$. We say that a trail $\alpha$ from $P_G$ to $R_G$ is a {\it new trail} if $\Phi^{-1}(\alpha)$ is not a trail from $P_E$ to $R_E$, where $\Phi$ is defined in the proof of Lemma \ref{lem:Trail-inequality}. Note that every new trail $\alpha$ from $P_G$ to $R_G$ can be written in the form $\Phi (\alpha_1\alpha_2)$,  where $\alpha_1$ is a trail from $P_E$ to $Q_E$ with range $v$, and $\alpha_2$ is a trail from $Q_E$ to $R_E$ with source $v$. By the numbering rule in  Remark \ref{rem10}, suppose that $\alpha_1$ has number $i$ and $\alpha_2$ has number $j$. Then we assign to the trail $\alpha$ the pair $(i,j).$
\end{rem} 

For clarification, we illustrate the above remarks by presenting the following example.

\begin{example}\label{exa4.9}
Let $E$ be the following graph
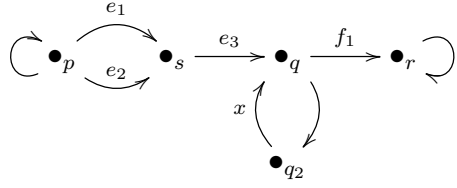
\[
\xymatrix{\bullet_{p}\ar@(ld,lu)\ar@/^1pc/[r]^{e_1}\ar@/_1pc/[r]^{e_2}&\bullet_{s}\ar[r]^{e_3}&\bullet_{q}\ar[r]^{f_1}\ar@/^1pc/[d]&\bullet_{r}\ar@(ru,rd)\\&&\bullet_{q_2}\ar@/^1pc/[u]^{x}}\]There are two trails from $P$ to $Q$, namely $\alpha_1 = e_1e_3$ and $\alpha_2 = e_2e_3$. We assign the number $1$ to $\alpha_1$ and $2$ to $\alpha_2$. There is one trail from  $Q$ to $R$, namely $\beta_1 = f_1$, which we number $1$. There are two paths from $P$ to $R$ passing through the intermediate vertex, namely $\gamma_1=e_1e_3f_1$ and $\gamma_2=e_2e_3f_1$. Note that this graph has no direct trail from  $P$ to $R$.  

Let $G$ be the graph obtained from $E$ by performing an
in-splitting at $q$ with partition $\mathcal{E} _1 = \{e_3\}$ and $ \mathcal{E}_2 = \{x\}.$ 
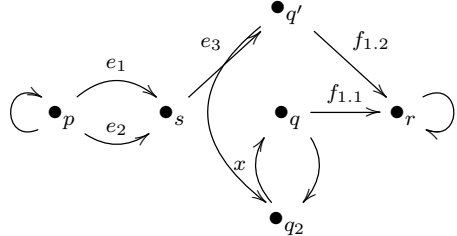
\[\xymatrix{&&\bullet_{q'}\ar@/_2.5pc/[dd]\ar[dr]^{f_{1.2}}\\\bullet_{p}\ar@(ld,lu)\ar@/^1pc/[r]^{e_1}\ar@/_1pc/[r]^{e_2}&\bullet_{s}\ar[ru]^{e_3}&\bullet_{q}\ar[r]^{f_{1.1}}\ar@/^1pc/[d]&\bullet_{r}\ar@(ru,rd)\\&&\bullet_{q_2}\ar@/^1pc/[u]^{x}}\] Then there are two new trails from the cycle $P$ to cycle $R$, namely $\gamma_1' = e_1e_3 f_{1,2}$ and $ \gamma_2'= e_2 e_3 f_{1,2}$. By the numbering rules given in Remarks \ref{rem10} and \ref{rem11}, we assign to $\gamma_1'$ the label $(1,1)$ and  to $\gamma_2'$ the label $(2,1)$. 
\end{example}

We are in a position to prove the first case of the main theorem.

\begin{prop}\label{maintheo-firstcase}
 Let $E$ and $F$ be meteor graphs of length three in normal form such that $F$ is obtained from $E$  by a finite sequence of in-splittings, out-splittings, in-amalgamations at interior vertices, and out-amalgamations at interior vertices. Then $E \approx F$.
\end{prop}
\begin{proof}
Since $F$ is obtained from $E$  by a finite sequence of in-splittings, out-splittings, in-amalgamations at interior vertices, and out-amalgamations at interior vertices, it follows from Theorem \ref{willimove}  that  $A_E\sim_{SSE} A_F$. By our conventions, $E$ and $F$ have chains of cycles $P_E >  Q_E > R_E$ and $P_F >  Q_F > R_F$, respectively, with 
\begin{center}
 $|P_E| = |P_F|:=x$, $|Q_E| = |Q_F|:=y$, and $|R_E| = |R_F|:=z$.    
\end{center}
Indeed we even have that $P_E=P_F$, $Q_E=Q_F$ and $R_E=R_F$ is a normalized set of disjoint cycles. 
Assume that there exists a sequence of meteor graphs of length three $S_1,S_2,\ldots, S_k$ such that  
\[ E=  S_1 \to S_2 \to \cdots \to S_k= F, \]
where $A \to B$ indicates that $B$ is obtained from $A$ via one of the following operations: in-
splitting, out-splitting, in-amalgamation at an interior vertex, or out-amalgamation at an interior vertex.
By Lemma \ref{prop-functor}, we obtain a  bijection \[
     \Phi^1 _{E,F}: \textnormal{Trail}_E (P_E,Q_E) \to \textnormal{Trail}_F(P_F, Q_F) 
    \]    
     a bijection \[
     \Phi^2_{E,F}: \textnormal{Trail}_E (Q_E,R_E) \to \textnormal{Trail}_F(Q_F,R_F)
    \] 
    and an injection \[ \Phi^3_{E,F}: \textnormal{Trail}_E (P_E,R_E) \to \textnormal{Trail}_F (P_F,R_F)
    \] such that 
\[f_E(\alpha) = \begin{cases}
	f_F(\Phi^1_{E,F}(\alpha))  \mbox{ if } \alpha \in \textnormal{Trail}_E(P_E, Q_E) \\
   f_F(\Phi^2_{E,F}(\alpha)) \mbox{ if } \alpha \in \textnormal{Trail}_E(Q_E, R_E)\\
f_F(\Phi^3_{E,F}(\alpha)) \mbox{ if } \alpha \in \textnormal{Trail}_E(P_E, R_E),\end{cases}\]    
where $f_E$ and $f_F$ are defined in Lemma \ref{prop-functor}. 

Consider a trail $\alpha$ from  $Q_E$ to  $R_E$ that has length $c$, and suppose that $\Phi^2_{E,F}(\alpha)$ has length $c'$. Then in $F$, by  Lemma \ref{prop-functor},  we have $$f_E(\alpha) = -c \equiv - c' = f_F(\Phi^2_{E,F}(\alpha)) \pmod{(y,z)},$$ since  $\alpha$ and $\Phi^2_{E,F} (\alpha)$ start at $q_1^E=q_1^F$, respectively, and end at $r_1^E=r_1^F$, respectively, in  $E$ and $F$. This implies that $c \equiv c' \pmod{(y,z)}$. Therefore, we have $$N_2^E(t) = N_2^F(t) \mbox{ for all } 1 \le t \le (y,z).$$  By the same argument, we also obtain that  $$N_1^E(t) = N_1^F(t) ~ \mbox{for all } 1 \le t \le (x,y).$$

We now compute $N_F(s) - N_E(s) $ for all $ 1 \le s \le (x,z)$. To do so, we number the new trails generated by graph moves according to the method described in Remark \ref{rem11}. Let $l_{P,Q}^E$ denote the number of all trails from $P_E$ to $Q_E$ in $E$. Similarly, let $l_{Q,R}^E$ denote the number of all trails from $Q_E$ to $R_E$ in $E$, and let $l_{P,R}^E$ denote the number of all trails from $P_E$ to $R_E$ in $E$.  

Let $\alpha$ be a trail from $P_E$ to $Q_E$ with numbering $c$, and let $\beta$ be a trail from $Q_E$ to $R_E$ with numbering $d$, where $1\le c \le  l^E_{P,Q} $ and $ 1 \le d \le l^E_{Q,R}$. By Remark \ref{rem10}, we may identify $\alpha$ with $\Phi^1_{E, S_i}(\alpha)$ in the graph $S_i$. Let $m_\alpha$ be the number of in-splittings (in the above sequence transforming $E$ to $F$) that are performed at $r(\alpha)$ and change the range vertex of $\alpha$; denote the set of these moves by $\mathcal{M}_\alpha$, so that $|\mathcal{M}_\alpha | = m_\alpha$. Similarly, let $n_\beta$ be the number of out-splittings (in the sequence transforming $E$ to $F$) that are performed at the source vertex of $\beta$ and change this source vertex; denote the set of these moves by $\mathcal{N}_\alpha$. 

Initially, the trail $\alpha$ ends at $q_1$ in $E$, and the trail $\beta$ starts at  $q_1$ in E. After performing $m_\alpha$ moves, the trail $\alpha$ ends at $q_{1+m_\alpha}$ in $F$, while the trail $\beta$ starts at $q_{1-n_\beta}$ in $F$. Since $F$ is in normal form, the trail $\alpha$ starts at $p_1$ and ends at $q_1$, and the trail $\beta$ starts at $q_1$ and ends at $r_1$ in $F$. Therefore, $y \mid m_\alpha$ and $y \mid n_\beta$.

Next, let $a_\alpha $ be the number of out-splittings (in the sequence transforming $E$ to $F$) performed at $s(\alpha)$ that increase the length of the trail $\alpha $ by $1$ and cause the trail to start at the preceding vertex in the first cycle. By the same argument, we have $ x \mid a_{\alpha}.$  
Using now that both $E$ and $F$ are in normal form, it follows from condition (2) in Definition \ref{def:normal-form} that 
the trail $\alpha$ has the same length in $E$ and $F$ modulo $xy$.  
The difference between these two lengths is $a_\alpha + m_\alpha$, so $xy \mid a_\alpha + m_\alpha$. It follows that $x \mid m_\alpha$, and hence $[x,y] \mid m_\alpha.$ Similarly, we obtain that $[y,z] \mid n_\beta.$

We note that a new trail with numbering $(c,d)$ is generated when applying a graph move on $S_i$ to obtain  $S_{i+1}$ only if we perform either an in-splitting in $\mathcal{M}_\alpha$ or an out-splitting in $\mathcal{N}_\beta$ at the vertex $v$, where $v\in Q_{S_i}$. After this operation, the number of trails with numbering $(c,d)$ may increase by $1$, or remain the same. Suppose the former happens, and denote this new trail by $\delta$. Notice that we have  $f_{S_{i+1}} (\delta) = f_{S_i} (\alpha,\beta)$. When we apply a graph move in $\mathcal{M}_\alpha$, the terminal vertex of $\alpha$
shifts to a vertex whose index increases by $1$. Similarly, when we apply a graph move in $\mathcal{N}_\beta$, the initial vertex of $\beta$ shifts to a vertex whose index decreases by $1$. Hence, if $\alpha$ ends and $\beta$ starts at the same vertex, then after performing $y$ graph moves in $\mathcal{M}_\alpha \cup \mathcal{N}_\beta$, $\alpha$ will again end and $\beta$ will again start at the same vertex. Moreover, observe that each time we perform a move in $\mathcal{M}_\alpha\cup \mathcal{N}_\beta$ to obtain $S_{i+1}$ from $S_i$, Lemma \ref{functor} implies that $$f_{S_{i+1}} (\alpha,\beta) =f_{S_i} (\alpha ,\beta) - 1.$$ Therefore, for a graph $S_i$ in which the terminal vertex of $\alpha$ coincides with the initial vertex of $\beta$, and for the graph $S_j $ where $\alpha$ and $\beta$  first meet again after $S_i$, we must have $$f_{S_j} (\alpha,\beta) = f_{S_i} (\alpha ,\beta)  - y,$$ 
because exactly $y$ graph moves from $\mathcal{M}_\alpha \cup \mathcal{N}_\beta$ are performed in the transition from $S_i$ to $S_j$.

From these observations, there exist $\dfrac{m_\alpha+n_\beta}{y}$ new trails $\delta_1,\delta_2, \ldots $ in $F$ such that $$f_F(\delta_1), f_F(\delta_2),\ldots $$ are equal to $$f_E(\alpha,\beta),f_E(\alpha,\beta) -y , \ldots , f_E(\alpha,\beta) - y\left(\dfrac{m_\alpha+n_\beta}{y} - 1 \right),$$ respectively. Since $x \mid m_\alpha$ and $z \mid n_\beta$, we have $ (x,z ) \mid  m _\alpha +n_\beta = y \left(  \dfrac{m_\alpha+n_\beta}{y} \right)$. Therefore, for each integer $c$ satisfying  $1 \le c \le  \dfrac{(x,z)}{(x,y,z)}$, there are exactly $$\dfrac{\frac{m_\alpha+n_\beta}{y}}{\frac{(x,z)}{((x,z),y)}} = \dfrac{(x,y,z)}{y(x,z)} (m_\alpha+n_\beta)$$ new trails $\delta$ such that $$f_F(\delta) \equiv  f_E(\alpha,\beta) + ((x,z),y) c  = f_E(\alpha,\beta) + (x,y,z) c \pmod{(x,z)}.$$ 
Furthermore, we have $$f_F(\delta) = r(\delta) - (s(\delta)+|\delta|)  = -|\delta|  \pmod{(x,z)},$$ since $F$ is in normal form, and $$f_E(\alpha,\beta) = r(\beta) - s(\alpha) - |\alpha| - |\beta| = -|\alpha| - |\beta|  \pmod{ (x,z) },$$ since $E$ is normal form. Therefore, we obtain that $$  |\delta| \equiv  |\alpha| + |\beta|  \pmod{(x,y, z)},$$ i.e.,  $\delta$ has length  $ |\alpha| + |\beta|  \pmod{(x,y,z) } $. 

We also note that for every $1\le h \le (x,z)$ with $h \equiv |\alpha| + |\beta| \pmod{(x,y,z)} $, there is a unique number $c$ such that $1 \le c \le \dfrac{(x,z)}{(x,y,z)}$ and \[ -h \equiv   f_E(\alpha,\beta) + (x,y,z) c \pmod{(x,z)}. \]
From these observations, we conclude that \[
N_F(h) - N_E(h) = N_F(g) - N_E(g) = \sum \limits_{  |\alpha| + |\beta|  \equiv h \pmod{(x,y,z)}  }  \dfrac{(x,y,z)}{y(x,z)}(m_\alpha +n_\beta)  \tag{*}
\]
for all $ 1 \le h,g \le (x,z)$ such that  $ h \equiv g \pmod{(x,y,z)} $. Notice that in $(*)$, the number of occurrences of the term $m_\alpha$ is equal to the number of elements of the set $$\{\beta \in \text{Trail}_E(Q_E,R_E)\mid |\alpha |+ |\beta | \equiv h \pmod{(x,y,z) } \}.$$ Hence, the number of times $m_\alpha$ appears in $(*)$ is 
\[ \sum \limits_{1 \le c \le (y,z),  c + |\alpha| \equiv h  \pmod{(x,y,z)}} N_2(c) = S_2(d),\]
where $d$ is the unique number such that $1\le d\le (x,y,z)$ and $d\equiv h -|\alpha | \pmod{(x,y,z)}$. Here recall that $N_2(c)= N_2^E(c)=N_2^F(c)$ denotes the number of trails whose lengths are congruent to $c\pmod{(y,z)}$. 

Let $l_a := \sum \limits_{ |\alpha | \equiv a \pmod{(x,y,z)} }m_\alpha $ for all $1 \le a \le (x,y,z) $ and let $g_b := \sum \limits_{|\beta | \equiv b \pmod{(x,y,z)} } n_\beta $ for all $1 \le b \le (x,y,z)$. Note that  $ [x,y] \mid l_a $ for all $ 1\le a \le (x,y,z)  $, $ [y,z] \mid g_b $ for all $ 1\le b \le (x,y,z)$, since $ [x,y] \mid m_\alpha$ and $[y,z] \mid n_\beta $. We also have $l_a = 0 $ if $S_1(a) = 0$, and $g_b= 0$ if $S_2(b) = 0 $.  Applying these observations to $(*)$, we obtain that 
\[ 
N_F(h)-N_E(h) = \dfrac{(x,y,z)}{y(x,z)} \left( \sum \limits_{\begin{subarray}{c}
    1 \le a \le (x,y,z)  \\ 
    1 \le b \le (x,y,z) \\
    a+b \equiv h \pmod{(x,y,z)}
\end{subarray} } S_1(a) g_b  + \sum \limits_{\begin{subarray}{c}
    1 \le a \le (x,y,z)  \\ 
    1 \le b \le (x,y,z) \\
    a+b \equiv h \pmod{(x,y,z)}
\end{subarray} } l_a S_2(b)   \right)  \tag{**} \]
for all $1 \le h \le (x,z),$ which shows that $ E \approx F$, thus finishing the proof.
\end{proof}

\begin{remark}
    \label{rem:positive-l-and-m}
Let $E$ and $F$ be two graphs satisfying the hypothesis of Proposition \ref{maintheo-firstcase}. In the proof of the proposition, we have indeed shown that $E\approx F$ with non-negative constants $l_a,g_b$ witnessing (**). In the next proposition, we will use a kind of converse of this result
in order to provide a proof of the implication $E\approx F \implies A_E\sim_{SSE} A_F$.
    \end{remark}

The following definition is very useful for the proofs of the main results of this section, culminating in Theorem \ref{numtheo}.

\begin{deff}
Let $E$ be a finite graph with disjoint cycles, and let $A$ and $B$ be two cycles in $E$. Let $\alpha = e_1e_2\cdots e_n$ be a trail from $A$ to $B$, where $e_1$ is an exit of $A$ and $e_n$ is an entrance of $B$. Performing an out-splitting at $s\left(e_1\right)$ with partition $\mathcal{E}_1=\left\{e_1\right\}$ and $\mathcal{E}_2=\left\{e \mid s(e)=s\left(e_1\right)\right\} \backslash \mathcal{E}_1$, we obtain a new graph $E^{\prime}$.  We refer to this operation as the {\it out-shift at $\alpha$}.  Similarly, performing in-splitting at $r\left(e_n\right)$ with partition $\mathcal{E}_1=\left\{e_n\right\}$ and $\mathcal{E}_2=\left\{e \mid r(e)=r\left(e_n\right)\right\} \backslash \mathcal{E}_1$, we obtain a new graph $E^{\prime}$.  We refer to this operation as the {\it in-shift at $\alpha$}.    
\end{deff}

We now have the necessary tools to prove the easy direction (the sufficiency) of the main theorem.

\begin{prop}\label{maintheo:easydiection}
Let $E$ and $F$ be meteor graphs of length three in normal form with $E\approx F$, and let
$A_E$ and $A_F$ denote the adjacency matrices of $E$ and $F$, respectively. Then $A_E\sim_{SSE} A_F$.
\end{prop}
\begin{proof}
Since $E \approx F$, there exist integers  $l_1,l_2,\ldots , l_{(x,y,z)} $ and $g_1,g_2 , \ldots ,g_{(x,y,z) } $ such that $l_a$ are multiples of $[x,y]$ for all $ 1\le a \le (x,y,z)$, and $g_b$ are multiples of  $[y,z]$ for all $ 1\le b \le (x,y,z)$. Moreover, $l_a = 0$ for all $1 \le a \le (x,y,z)$ with  $S_1(a) = 0 $, and $g_b= 0$ for all $1 \le b \le (x,y,z)$ with  $S_2(b) = 0 $. Furthermore, \[
N_F(h) -  N_E(h)  =  \dfrac{(x,y,z)}{y(x,z)} \left( \sum \limits_{\begin{subarray}{c}
    1 \le a \le (x,y,z)  \\ 
    1 \le b \le (x,y,z) \\
    a+b \equiv h \pmod{(x,y,z)}
\end{subarray} } S_1(a) g_b  + \sum \limits_{\begin{subarray}{c}
    1 \le a \le (x,y,z)  \\ 
    1 \le b \le (x,y,z) \\
    a+b \equiv h \pmod{(x,y,z)}
\end{subarray} } l_a S_2(b)    \right) 
\]  for all $1 \le h \le (x,z).$ 
We write $l_a$ in the form
\[l_a=l_a^E-l_a^F,\]
where
\begin{align*} (l_a^E,l_a^F)= \begin{cases} (l_a,0), & \text{if } l_a\ge 0,\\  (0,-l_a), & \text{if } l_a<0.\end{cases}
\end{align*}
Similarly, we express $g_b$ in the form $$g_b= g_b^E - g_b^F,$$ where 
\begin{align*} (g_b^E,g_b^F)= \begin{cases} (g_b,0), & \text{if } g_b\ge 0,\\  (0,-g_b), & \text{if } g_b<0.\end{cases}
\end{align*}
It follows that  $l_1^E,l_2^E,\ldots , l_{(x,y,z)}^E$, $l_1^F,l_2^F,\ldots , l_{(x,y,z)}^F$, $g_1^E,g_2^E, \ldots ,g_{(x,y,z)}^E$, and $g_1^F,g_2^F, \ldots ,g_{(x,y,z)}^F$ are natural numbers such that $l^E_a$ and $l^F_a$ are multiples of $[x,y]$ for all $ 1\le a \le (x,y,z) $, $g_b^E$ and $g_b^F$ are multiples of  $[y,z]$ for all $ 1\le b \le (x,y,z)$, $l_a^E = l_a^F = 0$ for all $1 \le a \le (x,y,z)$ with  $S_1(a) = 0 $,  $g_b^E=g_b^F=0$ for all $1 \le b \le (x,y,z)$ with  $S_2(b) = 0 $, and \[
N_E(h)  +  \dfrac{(x,y,z)}{y(x,z)} \left( \sum \limits_{\begin{subarray}{c}
    1 \le a \le (x,y,z)  \\ 
    1 \le b \le (x,y,z) \\
    a+b \equiv h \pmod{(x,y,z)}
\end{subarray} } S_1(a) g_b^E  + \sum \limits_{\begin{subarray}{c}
    1 \le a \le (x,y,z)  \\ 
    1 \le b \le (x,y,z) \\
    a+b \equiv h \pmod{(x,y,z)}
\end{subarray} } l_a^E S_2(b)    \right) 
\]   \[
= N_F(h)  +  \dfrac{(x,y,z)}{y(x,z)} \left( \sum \limits_{\begin{subarray}{c}
    1 \le a \le (x,y,z)  \\ 
    1 \le b \le (x,y,z) \\
    a+b \equiv h \pmod{(x,y,z)}
\end{subarray} } S_1(a) g_b^F + \sum \limits_{\begin{subarray}{c}
    1 \le a \le (x,y,z)  \\ 
    1 \le b \le (x,y,z) \\
    a+b \equiv h \pmod{(x,y,z)}
\end{subarray} } l_a^F S_2(b)    \right)
\]
for all $1 \le h \le (x,z).$ 

We claim that there exists a normal form graph $G$, obtained from $E$ by a finite sequence of graph moves, for which  \[ N_G(h) =   N_E(h)  +  \dfrac{(x,y,z)}{y(x,z)} \left( \sum \limits_{\begin{subarray}{c}
    1 \le a \le (x,y,z)  \\ 
    1 \le b \le (x,y,z) \\
    a+b \equiv h \pmod{(x,y,z)}
\end{subarray} } S_1(a) g_b^E  +\sum \limits_{\begin{subarray}{c}
    1 \le a \le (x,y,z)  \\ 
    1 \le b \le (x,y,z) \\
    a+b \equiv h \pmod{(x,y,z)}
\end{subarray} } l_a^E S_2(b)    \right)  \] for all $ 1 \le h \le (x,z)$.

We now obtain $G$ from $E$ by applying a sequence of graph moves as follows. We denote by  $\{a_1, a_2, \ldots, a_r\}$ the set of all integers $a$ such that $1 \le a \le (x,y,z)$ with $S_1(a) \ne 0$. Let $t_1$ be a positive integer such that $t_1$ is a multiple of $xy$ and $t_1 \ge \max\{l^E_{a_i}\mid 1\le i\le r\}$. For each $1\le i\le r$, choose an arbitrary trail $\alpha_i$ from $P_E$ to $Q_E$ whose length is congruent to $a_i$ modulo $(x,y,z) $. 
 
First, perform in-shifts at $\alpha_1$ exactly $ l_{a_1}^E$ times and out-shifts at $\alpha_1$ exactly $t_1-l_{a_1}^E$ to obtain a graph $E_1$ from $E$. The graph $E_1$ differs from $E$ only in that new trails from $P_{E_1}$ to $R_{E_1}$ are created and the length of $\alpha_1$ increases by $t_1$. In particular,  $S^E_1 = S^{E_1}_1$ (since $t_1$ is a multiple of both $x$ and $y$) and $S^E_2 = S^{E_1}_2$.   Next, perform in-shifts at $\alpha_2$ exactly $ l_{a_2}^E$ times and out-shifts at $\alpha_2$ exactly $t_1-l_{a_2}^E$ to obtain graph $E_2$ from $E_1$. The graph $E_2$ differs from $E_1$ only in that new trails from  $P_{E_2}$ to $R_{E_2}$ are created  and the length of $\alpha_2$ increases by $t_2$. In particular,  $S^E_1 = S^{E_1}_1 = S^{E_2}_1$ and $S^E_2 = S^{E_1}_2=S^{E_2}_2$. Repeating this process $r$ times, we arrive at  graph $E_r$, which differs from $E$ only in that new trails from $P_{E_r}$ to $R_{E_r}$ are created and the length of each $\alpha_i$ increases by $t_1$. In particular, $S^E_1 = S^{E_r}_1$ and $S^E_2 =S^{E_r}_2$. Since some trails in $E_r$ may have common interior vertices, Theorem \ref{thm:quasi-normal-form}  allows us to transform $E_r$ into a graph in quasi-normal form without increasing the number of trails.

Second, we denote by  $\{b_1, b_2, \ldots, b_s\}$ the set of all integers $b$ such that  $1 \le b \le (x,y,z)$ with $S_2(b) \ne 0$. Let $t_2$ be a positive integer such that $t_2$ is a multiple of $yz$ and $t_2\ge \max\{g^E_{b_i}\mid 1\le i\le s\}$. For each $1\le i\le s$, choose an arbitrary trail $\beta_i$ in $E_r$ from $Q_{E_r}$ to $R_{E_r}$ whose length is congruent to $b_i$  modulo $(x,y,z) $. Perform out-shifts at $\beta_1$ exactly $ g^E_{b_1}$ times and in-shifts at $\beta_1$ exactly $t_2- g^E_{b_1}$ times to obtain  graph $E_{r+1}$ from $E_r$. The graph $E_{r+1}$ differs from $E_r$ only in that new trails from $P_{E_{r+1}}$ to $R_{E_{r+1}}$ are created and the length of $\beta_1$ increases by $t_2$. In particular, $S^E_1 = S^{E_r}_1= S^{E_{r+1}}_1$ and $S^E_2 =S^{E_r}_2 = S^{E_{r+1}}_2$ (since $g_{b_1}^E$ is a multiple of  $yz$). Similarly, applying out-shifts at $\beta_2$ exactly $g^E_{b_2}$ times and in-shifts at $\beta_2$ exactly $t_2- g^E_{b_2}$ times, we obtain $E_{r+2}$ from $E_{r+1}$. Repeating this process $s$ times, we arrive at  graph $E_{r+s}$ which differs from $E_r$ only in that there are new trails from $P_{E_{r+s}}$ to $R_{E_{r+s}}$ are created and the length of each $\beta_i$ increases by $t_2$. In particular, $S^E_1 = S^{E_r}_1= S^{E_{r+s}}_1$ and $S^E_2 =S^{E_r}_2= S^{E_{r+s}}_2$. Since trails in $E_{r+s}$ may have common interior vertices, by applying Theorem \ref{thm:quasi-normal-form}, we can transform $E_{r+s}$ into a graph in quasi-normal form without increasing the number of trails.

Next, for each $1 \le a_i \le (x,y,z)$, perform out-shifts at every trail $\alpha \ne \alpha_i$ whose length is congruent to that of $\alpha_i$ modulo $(x,y)$, exactly $t_1$ times, in order to make $(P_E,Q_E)$ good. Similarly, for each $1 \le b_i \le (x,y,z)$, perform in-shifts at every trail $\beta \ne \beta_i$ whose length is congruent to that $\beta_i$ modulo $(y,z)$, exactly $t_2$ times, in order to make $(Q_E,R_E)$ good. 

Finally, by Theorem \ref{thm:normal-form}, we may apply out-shifts and in-shifts along trails from $P_E$ to $R_E$ to make $(P_E, R_E)$ good. 
 
After all the moves, we obtain the desired graph $G$. Since $x,y \mid [x,y] \mid l_a^E$ and $y,z \mid [y,z] \mid  g_b^E$,  in $G $ each trail from $P_G$ to $Q_G$ still starts at vertex  $p_1$ and ends at  $q_1$, each trail from $Q_G$ to $R_G$ still starts at  $q_1$ and ends at  $r_1.$ By (the proof of) Proposition \ref{maintheo-firstcase}, we have  \[
    N_G(h) =  N_E(h)  +  \dfrac{(x,y,z)}{y(x,z)} \left( \sum \limits_{\begin{subarray}{c}
    1 \le a \le (x,y,z)  \\ 
    1 \le b \le (x,y,z) \\
    a+b \equiv h \pmod{(x,y,z)}
\end{subarray} } S_1(a) g_b^E  +\sum \limits_{\begin{subarray}{c}
    1 \le a \le (x,y,z)  \\ 
    1 \le b \le (x,y,z) \\
    a+b \equiv h \pmod{(x,y,z)}
\end{subarray} } l_a^E S_2(b)    \right) 
    \] for all $1 \le h \le (x,z).$ 
    
By applying the same process to $F$, we can obtain a normal form graph $G'$ from $F$, where
\[ N_{G'}(h)= N_F(h)  +  \dfrac{(x,y,z)}{y(x,z)} \left( \sum \limits_{\begin{subarray}{c}
    1 \le a \le (x,y,z)  \\ 
    1 \le b \le (x,y,z) \\
    a+b \equiv h \pmod{(x,y,z)}
\end{subarray} } S_1(a) g_b^F + \sum \limits_{\begin{subarray}{c}
    1 \le a \le (x,y,z)  \\ 
    1 \le b \le (x,y,z) \\
    a+b \equiv h \pmod{(x,y,z)}
\end{subarray} } l_a^F S_2(b)    \right)\]
and $N_1^{G'} (c) = N_1^F(c) $ for all $1 \le c \le (x,y)$ and $N_2^{G'} (c) = N_2^F(c) $ for all $1 \le c \le (y,z)$. Therefore, we have  normal form graphs $G$ and $G'$ such that \[
\begin{cases}
    N_1^{G'} (i) = N_1^G(i) \text{ for all } 1 \le i \le (x,y) \\
    N_2^{G'}(i) = N_2^{G}(i) \text{ for all } 1 \le i \le (y,z) \\
    N_{G'}(i) = N_{G} (i) \text{ for all } 1 \le i \le (x,z).
\end{cases}
\]
Notice that, for all the trails from $P_G$ to $Q_G$ whose lengths are congruent to a natural number $i$ modulo $(x,y)$, we can perform $xy$ out-shifts on each trail, increasing each trail's length by $xy$ while keeping $G$ in normal form.  Similarly, for all trails from $P_G$ to $R_G$ whose lengths are congruent to $i$ modulo $(x,z)$, we can perform $xz$ out-shifts on each trail $xz$,  increasing each trail's length by $xz$ without leaving the normal form. For all trails from $Q_G$ to $R_G$ whose lengths are congruent to $i$ modulo $(y,z)$, we can perform $yz$ in-shifts on each trail, increasing each trail's length by $yz$ while preserving the normal form. The same process can be applied to $G'$. 

For each $1\le i \le (x, y)$, we denote by $\ell^G_{1, i}$ and $\ell^{G'}_{1, i}$ the common length of the trails from $P_G$ to $Q_G$ and from $P_{G'}$ to $Q_{G'}$, respectively, whose lengths are congruent to $i$ modulo $(x,y)$. Similarly, for each $1\le i \le (y, z)$, we denote by $\ell^G_{2, i}$ and $\ell^{G'}_{2, i}$ the common length of the trails from $Q_G$ to $R_G$ and from $Q_{G'}$ to $R_{G'}$, respectively, whose lengths are congruent to $i$ modulo $(y,z)$. Finally, for each $1\le i \le (x, z)$, we denote by $\ell^G_{3, i}$ and $\ell^{G'}_{3, i}$ the common length of the trails from $P_G$ to $R_G$ and from $P_{G'}$ to $R_{G'}$, respectively, whose lengths are congruent to $i$ modulo $(x,z)$. Note that by condition (2) in Definition \ref{def:normal-form}, $\ell^G_{1,i} \equiv \ell^{G'}_{1,i} \pmod{xy}$ for all $1\le i\le (x,y)$, and similar congruences hold also for $\ell^G_{2,i}$ and $\ell^{G'}_{2,i}$, and for $\ell^G_{3,i}$ and $\ell^{G'}_{3,i}$. Let 
\begin{center}
$\ell_{1,i} := \max \{\ell^G_{1, i}, \ell^{G'}_{1, i}\}$, $\ell_{2,j} := \max \{\ell^G_{2, j}, \ell^{G'}_{2, j}\}$, and $\ell_{3,k} := \max \{\ell^G_{3, k}, \ell^{G'}_{3, k}\}$
\end{center}
for all $1\le i \le (x, y)$, $1\le j \le (y, z)$, and $1\le k\le (x, z)$.

Let $H$ be the meteor graph of length three in normal form, with a unique chain of cycles 
$P_H > Q_H > R_H$, satisfying the following conditions:

$(1)$ $x = |P_H|$, $y = |Q_H|$, and $z = |R_H|$;

$(2)$ $N^H_1(i) = N^G_1(i) = N^{G'}_1(i)$ for all $1\le i \le (x, y)$, $N^H_2(i) = N^G_2(i) = N^{G'}_2(i)$ for all $1\le i \le (y, z)$, and $N_H(i) = N_G(i) = N_{G'}(i)$ for all $1\le i \le (x, z)$;

$(3)$ For each $1\le i\le (x,y)$, all trails from \(P_H\) to \(Q_H\) whose lengths are congruent to \(i\) modulo \((x,y)\) have length \(\ell_{1,i}\). Similarly, for each $1\le i\le (y,z)$, all trails from $Q_H$ to $R_H$ whose lengths are congruent to \(i\) modulo \((y,z)\) have length \(\ell_{2,i}\), and for each \(1\le i\le (x,z)\), all trails from \(P_H\) to \(R_H\) whose lengths are congruent to \(i\) modulo \((x,z)\) have length \(\ell_{3,i}\).

By the preceding argument, there exist sequences of graph moves on \(G\) and \(G'\) that transform both graphs into \(H\). Therefore, we conclude that $A_E \approx_{SSE} A_F$, thus finishing the proof.
\end{proof}

The remainder of this section is devoted to prove the main theorem by induction on the structure of a graph transformation sequence from $E$ to $F$. Proposition~\ref{maintheo-firstcase} establishes the result when $E$ can be transformed into $F$ through a sequence of in-splittings, out-splittings,
in-amalgamations at interior vertices, and out-amalgamations at interior vertices. It therefore remains to reduce the general case to this setting. To this end, we first characterize the meteor graphs of length three that can be obtained from a normal form graph by a sequence of in-splittings, out-splittings, in-amalgamations at interior vertices, and out-amalgamations at interior vertices. The following lemma identifies a class of such graphs. 

\begin{lemma}\label{lm:normal-cover}
Let $E$ be a meteor graph of length three with a unique chain of cycles 
$P_E > Q_E > R_E$, where $|P_E| = x$, $|Q_E| = y$, and $|R_E| = z$. Let $\mathrm{Trail}(E)$ denote the set of all trails in $E$. Assume that 
\[
|\alpha| \ge \max\{2xy + 2x + 2y,\; 2yz + 2y + 2z,\; 2zx +2 z + 2x\}.
\] for all $\alpha\in \mathrm{Trail}(E)$, and that all trails in $\mathrm{Trail}(E)$ have no common interior vertices.
Then there exist a meteor graph $G$ of length three in normal form and a meteor graph $F$ of length three such that the following conditions are satisfied:

$(1)$ $F$ is obtained from $G$ by a sequence of in-splittings and out-splittings;

$(2)$ $E$ is a subgraph of $F$;
    
$(3)$    $\mathrm{Trail}(F) \setminus \mathrm{Trail}(E) \subset \mathrm{Trail}_{F}(P_{F}, R_{F});$  

$(4)$ Any two distinct trails in $\mathrm{Trail}(F)$ have no common interior vertices;

$(5)$ All trails in $\mathrm{Trail}(F) \setminus \mathrm{Trail}(E)$ start and end in a vertex of index $1$; 

$(6)$ For all trails $\alpha$ and $\beta \in  \mathrm{Trail}(F) \setminus \mathrm{Trail}(E)$, if $$|\alpha| \equiv |\beta| \pmod{(x,z)},$$ then $|\alpha| = |\beta|$;

$(7)$ For a trail $\alpha \in \mathrm{Trail} (F)\setminus \mathrm{Trail}(E)$ and any $1\le i \le (x,z)$, if $$|\alpha| \equiv i \pmod{(x,z)},$$ then $$|\alpha| \equiv i \pmod{xz}.$$
\end{lemma}
\begin{proof}
 Let $a_i$ denote the number of all trails $\alpha$ from $P_E$ to $Q_E$ such that
\[
f_E(\alpha) \equiv -i \pmod{(x,y)}, \quad 1 \le i \le (x,y).
\]
Let $b_i$ denote the number of all trails $\alpha$ from $Q_E$ to $R_E$ such that
\[
f_E(\alpha) \equiv -i \pmod{(y,z)}, \quad 1 \le i \le (y,z).
\]
Similarly, let $c_i$ denote the number of all trails $\alpha$ from $P_E$ to $R_E$ such that
\[
f_E(\alpha) \equiv -i \pmod{(x,z)}, \quad 1 \le i \le (x,z).
\]
We choose a meteor graph $G$ of length three in normal form such that there are exactly $a_i$ trails from $p_1$ to $q_1$ of length $i$ for all $1 \le i < (x,y)$ and $a_{(x,y)} $ trails from $p_1$ to $q_1$ with length $xy$, exactly $b_i$ trails from $q_1$ to $r_1$ of length $i$ for all $1 \le i < (y,z)$ and $b_{(y,z)}$ trails from $q_1$ to $r_1$ with length $yz$, and exactly $c_i$ trails from $p_1$ to $r_1$ of length $i$ for each $1 \le i <(x,z)$ and $c_{(x,z)}$ trails from $p_1$ to $r_1$ with length $xz$.

For each $1 \le i \le (x,y)$, fix a bijection between the set of all trails of length $i$  (or $xy$) in $G$ from $P_G$ to $Q_G$ and the set of all trails $\alpha$ in $E$ from $P_E$ to $Q_E$ such that $f_E(\alpha) \equiv -i \pmod{(x,y)}$. Choose similar bijections for the other two kinds of trails in $G$ and $E$ respectively. For each trail $\alpha$ in $G$, we denote its corresponding trail in $E$ by $\Phi(\alpha)$.

Now let $\alpha$ be a trail in $G$, and suppose that $\Phi(\alpha)$ starts at $p_a$ and ends at $q_b$ in $E$. Perform out-shifts at $\alpha$ consecutively $x-a+1$ times to move its source to $p_a$, and in-shifts consecutively $b-1$ times to move its range to $q_b$. Observe that
\[
f_G(\alpha) \equiv -i \equiv f_E(\Phi(\alpha)) \pmod{(x,y)}.
\]
Thus, once $\alpha$ and $\Phi(\alpha)$ have the same source and range, the difference in their lengths (denoted by $l$) is divisible by $(x,y)$. Moreover,
\[
l \ge 2xy + 2x + 2y - i- (x-a+1) - (b-1) -xy\ge xy.
\]
We claim\footnote{The claim holds for every number $l$ that is greater than the Frobenius number of $x$ and $y$.}  that there exist non-negative integers $e$ and $f$ such that
\[ex + fy = l.\]
Indeed, since $(x,y) \mid l$, there exists a positive integer with $1 \le e \le  y$ such that $$ex \equiv l \pmod{y}.$$ Define $f= \dfrac{l-ex}{y}$. Since $l \ge xy$, we have that $f$ is a non-negative integer, and hence $ex+fy=l$, as desired. 

We then perform $ex$ out-shifts and $fy$ in-shifts on $\alpha$ so that its length matches that of $\Phi(\alpha)$.

Applying this procedure to all trails in $G$, we obtain a graph $F$ such that $E$ is a subgraph of $F$, and
\[
\mathrm{Trail}(F) \setminus \mathrm{Trail}(E) \subset \mathrm{Trail}_{F}(P_{F}, R_{F}).
\] Applying the method from the proof of Theorem \ref{thm:quasi-normal-form} to all trails in $\mathrm{Trail}(F)  $ shows that we can manage to obtain a new graph, also denoted by $F$, such that the new graph has pairwise disjoint interiors and all trails in $\mathrm{Trail}(F) \setminus \mathrm{Trail}(E)$ start and end in a vertex of index $1$. Hence the graph $F$ constructed so far satisfies (1)-(5).  

Now, using that the graph $F$ satisfies (3) and (4), we may apply the method introduced in the last paragraph of the proof of Theorem \ref{thm:normal-form} to all the trails in $\mathrm{Trail}(F) \setminus \mathrm{Trail}(E)$ to ensure that $F$ also satisfies conditions $(6)$ and $(7)$. This concludes the proof. 
\end{proof}

Consequently, we get the following useful corollary.

\begin{corollary}\label{cor:normal-cover}
Let $E$ be a meteor graph of length three in quasi-normal form with a unique chain of cycles 
$P_E > Q_E > R_E$, such that $(P_E, Q_E)$ and $(Q_E, R_E)$ are good, where $|P_E| = x$, $|Q_E| = y$, and $|R_E| = z$. Let $\mathrm{Trail}(E)$ denote the set of all trails in $E$. Assume that 
\[
|\alpha| \ge \max\{2xy + 2x + 2y,\; 2yz + 2y + 2z,\; 2zx +2 z + 2x\}.
\] for all $\alpha\in \mathrm{Trail}(E)$. 
Then there exists a meteor graph $G$ of length three in normal form  such that $E$ is obtained from $G$ by a sequence of in-splittings and out-splittings.
\end{corollary}
\begin{proof}
Let $c_i$ denote the number of all trails $\alpha$ from $P_E$ to $R_E$ such that
\[
f_E(\alpha) \equiv -i \pmod{(x,z)}, \quad 1 \le i \le (x,z).
\]
We choose a meteor graph $G$ of length three in normal form such that, for each $1 \le i \le (x,z)$, there are exactly $c_i$ trails from $p_1$ to $r_1$ of length $i$. Furthermore, the subgraphs of $G$ and $E$, obtained by deleting all trails from the first cycle to the third cycle, are isomorphic. Then, applying the algorithm introduced in Lemma \ref{lm:normal-cover} to each trail in $\text{Trail}_G(P_G, R_G)$, 
we obtain $E$ from $G$ through a finite sequence of in-splittings and out-splittings, thus finishing the proof.  
\end{proof}

We next develop a procedure that extends an arbitrary meteor graph of length three to a graph of the type described in Lemma \ref{lm:normal-cover}, while preserving the in-splitting and out-splitting operators. This construction is carried out in Lemmas \ref{lem:first-type-edge-ext}, \ref{lem:second-type-edge-ext}, and \ref{lm:third-edge-ext}.
We now introduce the notion of first-type edge extensions for meteor graphs of length three.

We say that an edge $e$ of a meteor graph of length three $E$ is a {\it non-cycle edge} if $e$ is not contained in any cycle of $E$.

\begin{deff}\label{deff:first-type-edge-ext}
Let $E$ be a meteor graph of length three with a unique chain of cycles 
$P_E > Q_E > R_E$, and let $c$ be a positive integer. 

Let $S \subseteq E^1$ be a set of non-cycle edges satisfying the following condition:
for every trail $\alpha$ starting at a vertex in $P^0_E$, there exists exactly one edge $e \in S$ such that $e \in \alpha$, while for any trail $\alpha$ whose source is not in $P^0_E$, one has $\alpha \cap S = \varnothing$. 
We refer to this as {\it Condition} $(1)$.

For each $e \in S$, let
\[
V_e = \{v_0^e, v_1^e, \ldots, v_{c-1}^e\}
\text{ and }
E_e = \{e_1^e, e_2^e, \ldots, e_{c+1}^e\},
\]
where the edges satisfy
\[s(e_1^e) = s(e),\
s(e_i^e) = r(e_{i-1}^e) = v_{i-2}^e \text{ for all }  2 \le i \le c+1, \text{ and }
r(e_{c+1}^e) = r(e).\] 
We form a new graph by adjoining the vertices in $V_e$ and the edges in $E_e$ to $E$, and replacing each edge $e \in S$ with the path $$p_e=e_1^e e_2^e \cdots e_{c+1}^e.$$ We denote the resulting graph by $E_{(S,c,1)}$.
\end{deff}

For clarification, we illustrate Definition \ref{deff:first-type-edge-ext} by presenting the following example.

\begin{example}\label{exa:first-type-edge-ext}
Consider the graph $E$ from Example \ref{exa4.9}, and let $S = \{e_1, e_2\}$. The graph $E_{(S, 1, 1)}$ is then as follows:
\end{example}
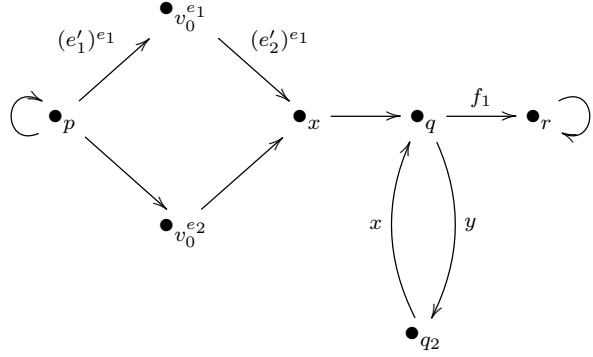
\[\xymatrix{
 &  & \bullet_{v_0^{e_1}} \ar@{->}[rd]^{(e_2')^{e_1}} &  &  &  \\
{} & \bullet_p \ar@{->}[ru]^{(e_1')^{e_1}} \ar@{->}[rd] \ar@(ld,lu) &  & \bullet_x \ar@{->}[r] & \bullet_q \ar@/^1pc/@{->}[dd]^{y} \ar@{->}[r]^{f_1} & \bullet_r  \ar@(ru,rd) \\
 &  & \bullet_{v_0^{e_2}} \ar@{->}[ru] &  &  &  \\
 &  &  &  & \bullet_{q_2} \ar@/^1pc/@{->}[uu]^{x} & 
}\]

It is worth mentioning the following remark.

\begin{remark}\label{rem:first-type-edge-ext}
We observe that $E_{(S,c,1)}$ admits the representation $$(((E_{(S_1,1,1)})_{(S_2,1,1)}) \ldots  ) _{(S_c,1,1)},$$ where $S_1=S$ and, for each $1 \le i \le c-1,$ $$S_{i+1}= \{e_1^\alpha \mid \alpha \in S_i \}.$$

Moreover, for any trail $f_1 \alpha f_2$ in $E$ that starts at a vertex in $P^0_E$, where $\alpha \in S$ and $f_1, f_2 \in E^*$, we define the corresponding trail in $E_{(S,c,1)}$ by
\[f_1\, e_1^\alpha e_2^\alpha \cdots e_{c+1}^\alpha\, f_2.\]
Consequently, every trail in $E$ starting from $P_E$ has its length increased by $c$ when viewed in
 $E_{(S,c,1)}$. 
\end{remark}

\begin{lemma}\label{lem:first-type-edge-ext}
Let $E$ and $F$ be two meteor graphs of length three such that $F$ can be obtained from $E$ by an in-splitting or an out-splitting. Let  $S\subseteq E^1$ be a set of non-cycle edges satisfying Condition $(1)$ in E. Then there exists a set $S' \subseteq F^1$ satisfying Condition $(1)$ in $F$ such that the graph $F_{(S',1,1)}$ can be obtained from $E_{(S,1,1)}$ by a finite sequence of in-splittings and out-splittings. In particular, if $S$ is the set of all exits for the first cycle in $E$, then $S'$ is precisely the set of all exits for the first cycle in $F$.
\end{lemma}
\begin{proof}
Suppose that $F$ is obtained from $E$ by performing an out-splitting at $v$ with respect to the partition 
\[
s^{-1}(v) = \mathcal{E}_1 \cup \mathcal{E}_2 \cup \cdots \cup \mathcal{E}_n.
\]
Let $S\subseteq E^1$ be a set of non-cycle edges satisfying Condition $(1)$. We consider the following cases:

\noindent\textit{Case 1}: $r^{-1}(v) \cap S = \varnothing$. In this case, $S$ also satisfies Condition $(1)$ in $F$. Consequently, performing an out-splitting at $v$ in $E_{(S,1,1)}$ with respect to the same partition as in $E$ yields the graph $F_{(S,1,1)}$.
    
\noindent\textit{Case} $2$: $T:= r^{-1}(v) \cap S \neq \varnothing$. In the graph $F$, each edge $t \in T$ is replaced by a family of edges 
\[   T' = \{ t_i \mid t \in T,\; 1 \le i \le n \}.    \]
Hence, 
    \[    S' = (S \setminus T) \cup T'    \]
satisfies Condition $(1)$ in $F$. We now show that $F_{(S',1,1)}$ can be obtained from $E_{(S,1,1)}$ via a sequence of out-splittings.

First, we perform an out-splitting at $v$ in $E_{(S,1,1)}$ using the same partition as in $E$, and denote the resulting graph by $A$. Observe that for each edge $e_2^\alpha$ with $\alpha \in T$ in $E_{(S,1,1)}$, it is replaced in $A$ by a family of edges $\{(e_2^\alpha)_i\}_{i=1}^n$ satisfying
    \[
    s\big((e_2^\alpha)_i\big) = s(e_2^\alpha) = v_0^\alpha, 
    \qquad 
    r\big((e_2^\alpha)_i\big) = v_i, \quad 1 \le i \le n.
    \]

Next, in $A$, we perform out-splittings at each vertex $v_0^\alpha$ for all $\alpha \in T$, using the partition 
    \[
    \mathcal{E}_1 \cup \mathcal{E}_2 \cup \cdots \cup \mathcal{E}_n,
    \quad \text{where } \mathcal{E}_i = \{(e_2^\alpha)_i\}, \; 1 \le i \le n.
    \]
    Denote the resulting graph by $B$. In $B$, each vertex $v_0^\alpha$ is replaced by vertices $(v_0^\alpha)_1, (v_0^\alpha)_2, \ldots,$ $ (v_0^\alpha)_n$, and each edge $e_1^\alpha$ is replaced by edges $(e_1^\alpha)_i$ for $1 \le i \le n$ such that
    \[
    s\big((e_1^\alpha)_i\big) = s(\alpha), 
    \qquad 
    r\big((e_1^\alpha)_i\big) = (v_0^\alpha)_i, 
    \quad 1 \le i \le n.
    \]
    
Finally, identifying each vertex $(v_0^\alpha)_i$ in $B$ with $v_0^{\alpha_i}$ in $F$, and each edge $(e_j^\alpha)_i$ with $e_j^{\alpha_i}$, for all $\alpha \in T$, $1\le j \le 2$ and $1 \le i \le n$, we obtain a bijective correspondence between $B$ and $F_{(S',1,1)}$. Therefore, $F_{(S',1,1)}$ is obtained from $E_{(S,1,1)}$ via a finite sequence of out-splittings. 

We can handle the case where $F$ is obtained from $E$ by performing an in-splitting at a vertex of $E$ using a symmetric argument. We leave the details to the reader.
\end{proof}

The following example illustrates the main ideas of Lemma \ref{lem:first-type-edge-ext}.

\begin{example} Consider a graph $E$ in the following form:
   
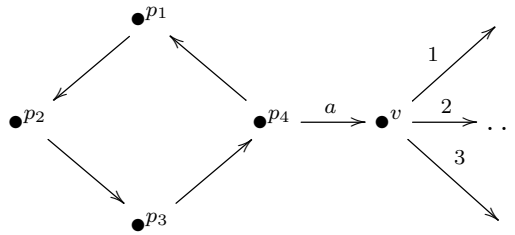
\[\xymatrix{
 & \bullet^{p_1} \ar@{->}[ld] &  &  & {} \\
\bullet^{p_2} \ar@{->}[rd] &  & \bullet^{p_4} \ar@{->}[lu] \ar@{->}[r]^{a} & \bullet^v \ar@{->}[ru]^{1} \ar@{->}[r]^{2} \ar@{->}[rd]^{3} & \ldots \\
 & \bullet^{p_3} \ar@{->}[ru] &  &  & {}
}\]
We take $S= \{a\}$ in $E$ to obtain $E_{(S,1,1)}$
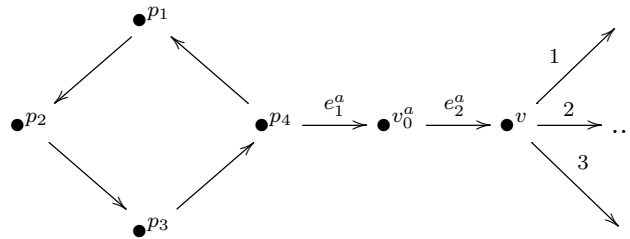
\[\xymatrix{
 & \bullet^{p_1} \ar@{->}[ld] &  &  &  & {} \\
\bullet^{p_2} \ar@{->}[rd] &  & \bullet^{p_4} \ar@{->}[lu] \ar@{->}[r]^{e_1^a} & \bullet^{v_0^a} \ar@{->}[r]^{e_2^a} & \bullet^v \ar@{->}[ru]^{1} \ar@{->}[r]^{2} \ar@{->}[rd]^{3} & ... \\
 & \bullet^{p_3} \ar@{->}[ru] &  &  &  & {}
}\]
By performing an out-splitting at $v$ with respect to the partition $\mathcal{E}_1= \{1,2\}$ and $ \mathcal{E}_2= \{3\}$, we obtain the graph $F$:
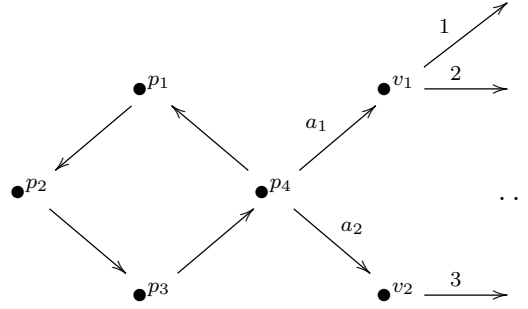
\[\xymatrix{
 &  &  &  & {} \\
 & \bullet^{p_1} \ar@{->}[ld] &  & \bullet^{v_1} \ar@{->}[ru]^{1} \ar@{->}[r]^{2} & {} \\
\bullet^{p_2} \ar@{->}[rd] &  & \bullet^{p_4} \ar@{->}[lu] \ar@{->}[ru]^{a_1} \ar@{->}[rd]^{a_2} &  & \ldots \\
 & \bullet^{p_3} \ar@{->}[ru] &  & \bullet^{v_2} \ar@{->}[r]^{3} & {}
}\]
Following the algorithm, we obtain $S'=\{a_1,a_2\}$. The graph $F_{(S',1,1)}$ is then given as follows:
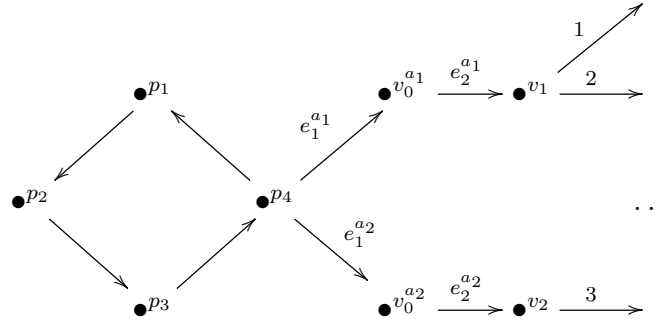
\[
\xymatrix{
 &  &  &  &  & {} \\
 & \bullet^{p_1} \ar@{->}[ld] &  & \bullet^{v_0^{a_1}} \ar@{->}[r]^{e_2^{a_1}} & \bullet^{v_1} \ar@{->}[ru]^{1} \ar@{->}[r]^{2} & {} \\
\bullet^{p_2} \ar@{->}[rd] &  & \bullet^{p_4} \ar@{->}[lu] \ar@{->}[ru]^{e_1^{a_1}} \ar@{->}[rd]^{e_1^{a_2}} &  &  & \ldots  \\
 & \bullet^{p_3} \ar@{->}[ru] &  & \bullet^{v_0^{a_2}} \ar@{->}[r]^{e_2^{a_2}} & \bullet^{v_2} \ar@{->}[r]^{3} & {}
}
\]
By the algorithm in Lemma \ref{lem:first-type-edge-ext}, we perform an out-splitting at $v$  in $E_{(S,1,1)}$ to  obtain $A$. 

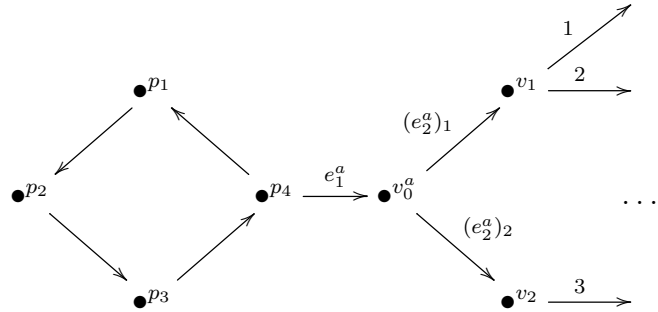
\[\xymatrix{
 &  &  &  &  & {} \\
 & \bullet^{p_1} \ar@{->}[ld] &  &  & \bullet^{v_1} \ar@{->}[ru]^{1} \ar@{->}[r]^{2} & {} \\
\bullet^{p_2} \ar@{->}[rd] &  & \bullet^{p_4} \ar@{->}[lu] \ar@{->}[r]^{e_1^a} & \bullet^{v_0^a} \ar@{->}[ru]^{(e_2^a)_1} \ar@{->}[rd]^{(e_2^a)_2} &  & \ldots \\
 & \bullet^{p_3} \ar@{->}[ru] &  &  & \bullet^{v_2} \ar@{->}[r]^{3} & {}
}\]
Next, we perform an out-splitting at $v_0^a$ in  $A$ to obtain $B$ 
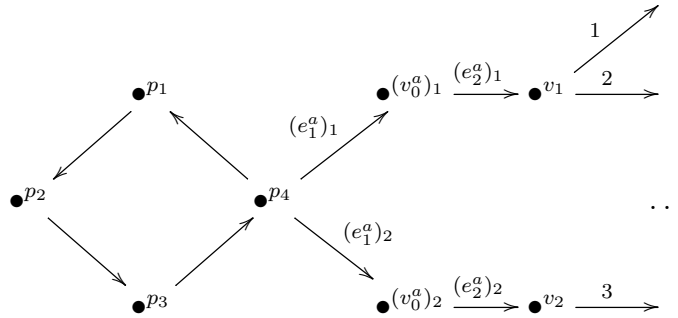
\[
\xymatrix{
 &  &  &  &  & {} \\
 & \bullet^{p_1} \ar@{->}[ld] &  & \bullet^{(v_0^a)_1} \ar@{->}[r]^{(e_2^a)_1} & \bullet^{v_1} \ar@{->}[ru]^{1} \ar@{->}[r]^{2} & {} \\
\bullet^{p_2} \ar@{->}[rd] &  & \bullet^{p_4} \ar@{->}[lu] \ar@{->}[ru]^{(e_1^a)_1} \ar@{->}[rd]^{(e_1^a)_2} &  &  & \ldots \\
 & \bullet^{p_3} \ar@{->}[ru] &  & \bullet^{(v_0^a)_2} \ar@{->}[r]^{(e_2^a)_2} & \bullet^{v_2} \ar@{->}[r]^{3} & {}
}
\]
\end{example}

By Remark \ref{rem:first-type-edge-ext} and Lemma \ref{lem:first-type-edge-ext}, we immediately obtain the following useful corollary. 

\begin{corollary}\label{cor:first-type-edge-ext}
Let $E$ and $F$ be two meteor graphs of length three such that $F$ can be obtained from $E$ by an in-splitting or an out-splitting. Let $c$ be a positive integer, and let $S$ be a set of non-cycle edges satisfying Condition $(1)$ in $E$.  Then there exists a set $S'$ in $F$ satisfying Condition $(1)$  such that $F_{(S',c,1)}$ can be obtained from $E_{(S,c,1)}$ via a finite sequence of in-splittings and out-splittings.
In particular, if $S$ is the set of all exits for the first cycle in $E$, then $S'$ is precisely the set of all exits for the first cycle in $F$.
\end{corollary}

We next introduce the notion of second-type edge extensions for meteor graphs of length three.

\begin{deff}\label{deff:second-type-edge-ext}
Let $E$ be a meteor graph of length three with a unique chain of cycles $P_E > Q_E > R_E$, and let $c$ be a positive integer. 

Let $S \subseteq E^1$ be a set of non-cycle edges satisfying the following condition:
for every trail $\alpha$ ending at a vertex in $R^0_E$, there exists a unique edge $e \in S$ such that $e \in \alpha$, while for any trail $\alpha$ not ending at a vertex in $R^0_E$, we have $\alpha \cap S = \varnothing$. 
We refer to this as {\it Condition} $(2)$.

For each $\beta \in S$, let
\[
V_\beta = \{v_0^\beta, v_1^\beta, \ldots, v_{c-1}^\beta\}
\text{ and }
E_\beta = \{e_1^\beta, e_2^\beta, \ldots, e_{c+1}^\beta\},
\]
where the edges satisfy
\[s(e_1^\beta) = s(\beta),\
s(e_i^\beta) = r(e_{i-1}^\beta) = v_{i-2}^\beta  \text{ for all } 2 \le i \le c+1, \text{ and }
r(e_{c+1}^\beta) = r(\beta).\]
We form a new graph by adjoining the vertices in $V_\beta$ and the edges in $E_\beta$ to $E$, and replacing each edge $\beta \in S$ with the  path $$p_\beta= e_1^\beta e_2^\beta \cdots e_{c+1}^\beta.$$The resulting graph is denoted by $E_{(S,c, 2)}$.
\end{deff}

We note that for any trail $f_1 \beta f_2$ ending at a vertex in $R^0_E$, where $\beta \in S$ and $f_1, f_2 \in E^*$, we define the corresponding trail in $E_{(S,c, 2)}$ by
\[f_1\, e_1^\beta e_2^\beta \cdots e_{c+1}^\beta\, f_2.\]
Consequently, every trail in $E$ ending at a vertex in $R^0_E$ has its length increased by $c$ when viewed in $E_{(S,c,2)}$.

For clarification, we illustrate Definition \ref{deff:second-type-edge-ext} by presenting the following example.

\begin{example}
Consider the graph $E$ from Example \ref{exa4.9}, and let $S = \{f_1\}$. The graph $E_{(S, 1, 2)}$ is then as follows: 
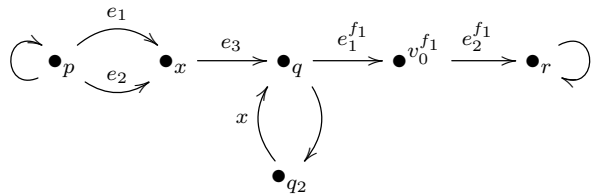
\[\xymatrix{\bullet_{p}\ar@(ld,lu)\ar@/^1pc/[r]^{e_1}\ar@/_1pc/[r]^{e_2}&\bullet_{x}\ar[r]^{e_3}&\bullet_{q}\ar[r]^{e^{f_1}_1}\ar@/^1pc/[d]&\bullet^{v_0^{f_1}}\ar[r]^{e^{f_1}_2}&\bullet_{r}\ar@(ru,rd)\\&&\bullet_{q_2}\ar@/^1pc/[u]^{x}}\]
\end{example}

By an argument analogous to that used in Corollary \ref{cor:first-type-edge-ext}, we obtain the following lemma.

\begin{lemma}\label{lem:second-type-edge-ext}
Let $E$ and $F$ be two meteor graphs of length three such that $F$ can be obtained from $E$ by an in-splitting or an out-splitting. Let $c$ be a positive integer, and let $S$ be a set of non-cycle edges satisfying Condition $(2)$ in $E$.  Then there exists a set $S'$ in $F$ satisfying Condition $(2)$  such that $F_{(S',c,2)}$ can be obtained from $E_{(S,c,2)}$ via a finite sequence of in-splittings and out-splittings.
In particular, if $S$ is the set of all entrances for the third cycle in $E$, then $S'$ is precisely the set of all entrances for the third cycle in $F$.
\end{lemma}

Motivated by Lemma \ref{lm:normal-cover}(3), we introduce the notion of third-type edge extensions for meteor graphs of length three. We first recall that a {\it path graph} is a graph of the form: $$\xymatrix{  \bullet^{v_{m}} \ar[r]^{e_{m}} & \cdots\ar[r]^{e_3}& \bullet^{v_2} \ar[r]^{e_2} & \bullet^{v_1} \ar[r]^{e_1} &v_0.}$$

\begin{deff}\label{deff:third-edge-ext}
Let $E$ be a meteor graph of length three with a unique chain of cycles $P_E > Q_E > R_E$, and let $S$ be a finite collection of path graphs.
Define $E_{(S, 3)}$ to be the meteor graph obtained from $E$ as follows: for each $\alpha\in S$, attach a path of length $|\alpha|$ of the form
$$\xymatrix{  \bullet^{v_{|\alpha|}^{\alpha}} \ar[r]^{e^{\alpha}_{|\alpha|}} & \cdots\ar[r]^{e_3^{\alpha}}& \bullet^{v^{\alpha}_2} \ar[r]^{e_2^{\alpha}} & \bullet^{v^{\alpha}_1} \ar[r]^{e^{\alpha}_1} &v_0}$$ by identifying $v_{|\alpha|}^{\alpha}$ with $p_1$ and $v_0$ with $r_1$, where $p_1\in P_E^0$ and $r_1\in R_E^0$ are the vertices of index $1$.
\end{deff}

\begin{lemma}\label{lm:third-edge-ext}
Let $E$ and $F$ be two meteor graphs of length three such that $F$ is obtained from $E$ by either an in-splitting or an out-splitting, and let  $S$ be a finite collection of path graphs. Then $F_{(S, 3)}$ is  obtained from $E_{(S, 3)}$ by the corresponding in-splitting or out-splitting. 
\end{lemma}
\begin{proof}
We consider the case where $F$ is obtained from $E$ by performing an out-splitting at a vertex $v\in E^0$. 
Consider the following cases:

{\it Case} $1$: $v\neq p_1$. In this case, we have
\[s_E^{-1}(v) = s_{E_{(S, 3)}}^{-1}(v).\]
Hence, performing the same out-splitting at $v$ in $E_{(S, 3)}$ yields $F_{(S, 3)}$.

{\it Case} $2$: $v= p_1$. Suppose that an out-splitting is performed at $p_1$ in $E$ with respect to a partition
\[s^{-1}(p_1) = \mathcal{E}_1 \cup \mathcal{E}_2 \cup \cdots \cup \mathcal{E}_n,\]
where $\mathcal{E}_1 \cap (P_E)^1 \neq \varnothing$. Define a corresponding partition in $E_{(S, 3)}$ by
\begin{center}
$s^{-1}(p_1) = \mathcal{E}_1' \cup \mathcal{E}_2 \cup \cdots \cup \mathcal{E}_n,$ where
$\mathcal{E}_1' = \mathcal{E}_1 \cup \{ e^{\alpha}_{|\alpha|} \mid \alpha\in S\}.$
\end{center}
Performing the out-splitting with respect to this partition yields $F_{(S, 3)}$. The in-splitting case follows by a symmetric argument, completing the proof.
\end{proof}

We are now in a position to present the  main result of this section, which provides a number-theoretic criterion for meteor graphs of length three in normal form to be strongly shift equivalent.

\begin{theorem}\label{numtheo}
Let $E$ and $F$ be meteor graphs of length three in normal form, and let $A_E$ and $A_F$ be the adjacency matrices of $E$ and $F$, respectively.  Then, $A_E\sim_{SSE} A_F$ if and only if  $E \approx F.$	
\end{theorem}

\begin{proof}
$(\Longleftarrow)$ It immediately follows from Proposition \ref{maintheo:easydiection}.

$(\Longrightarrow)$ Let $E$ and $F$ be  meteor graphs of length three in normal form such that  $A_E\sim_{SSE} A_F$. By Proposition \ref{propreserve}, $E$ and $F$ have chains of cycles $P_E >  Q_E > R_E$ and $P_F >  Q_F > R_F$, respectively, with $|P_E| = |P_F|:=x$, $|Q_E| = |Q_F|:=y$, and $|R_E| = |R_F|:=z$. Since $A_E\sim_{SSE} A_F$, it follows from Theorem \ref{willimove} that $E$ can be transformed into $F$ via a finite sequence of in-splittings, out-splittings,  in-amalgamations, and out-amalgamations. Thus, there exists a finite sequence of meteor graphs $S_0,S_1,S_2,\ldots ,S_m$ such that  $
 E= S_0$,  $S_m = F$, and for each $i$, the graph $S_{i+1}$ is obtained from $S_i$ by applying one of the following operations: in-splitting, out-splitting, in-amalgamation, or out-amalgamation. Note that if  $S_{i+1}$ is obtained from  $S_i$ by an in-amalgamation (respectively, out-amalgamation), then $S_i$ can be obtained from $S_{i+1}$ by an in-splitting (respectively, out-splitting). Therefore,  the above sequence may be written in the form 
\[
\scalebox{0.75}{
$\xymatrixrowsep{0.4pc}\xymatrix{E\ar[r]&S_1\ar[r]&\cdots\ar[r]&S_{i_1}&S_{i_1+1}\ar[l]&\ar[l]\cdots&\ar[l]S_j\ar[r]&S_{j+1}\ar[r]&\cdots\ar[r]&S_{i_2}&S_{i_2+1}\ar[l]&\ar[l]\cdots&\ar[l]S_k\ar[ld]\\
&&&&&&&&&&&\ar@{..}[llllddddd]&\\
&&&&&&&&&&&&\\
&&&&&&&&&&&&\\
&&&&&&&&&&&&\\ 
&&&&&&&&&&&&&\\
&&&&&&&&&&&&&\\
&&&&&&\ar[ru]S_t\ar[r]&S_{t+1}\ar[r]&\cdots\ar[r]&S_{i_n}&S_{i_n+1}\ar[l]&\ar[l]\cdots&\ar[l]S_m=F,}$
}
\] 
where $A \to B$ denotes that $B$ is obtained from $A$ via one of the following operations: in-
splitting, out-splitting, in-amalgamation at an interior vertex, or out-amalgamation at an interior vertex.  
Although $n$ is not unique, we may choose a transformation from $E$ to $F$ such that $n$ is minimized. We refer to this minimal value of $n$ as the {\it transform degree of $(E, F)$}, denoted by $\text{tr.deg}(E, F)$.

In this diagram, we admit $i_1= 0$, i.e. $E=S_{i_1}$, and $m=i_{n}$, i.e. $F=S_{i_n}$. However, we assume that 
$i_1<j<i_2<k<i_3 < \cdots < i_n$. Observe that any sequence witnessing $A_E\sim_{SSE} A_F$ can be represented in this form.
Moreover, if necessary, we may insert one or both of the short sequences 
\begin{center}
$E\to E_1\to E_2=E$ and $F= F_2\leftarrow F_1\leftarrow F$,    
\end{center}
where $E_1$ (resp. $F_1$) is obtained from $E$ (resp. $F$) by an in-amalgamation at an interior vertex and $E_2$ (resp. $F_2$) is obtained from $E_1$ (resp. $F_1$) by the inverse move, inserted before $E$ and after $F$, respectively, and thus we may always assume that $i_1 > 0$ and $i_n < m$. 

We use induction on $\text{tr.deg}(E, F)$ to establish the result. If  $\text{tr.deg}(E, F) = 1$, then 
there exists a sequence of meteor graphs of length three $S_1,S_2,\ldots, S_k, S_{k+1}, S_{k+2},\ldots, S_m$ such that  
\begin{equation}
 E \to  S_1 \to S_2 \to \cdots \to S_k\leftarrow S_{k+1} \leftarrow  S_{k+2} \leftarrow\cdots \leftarrow S_m = F, 
 \end{equation}
where $A \to B$ indicates that $B$ is obtained from $A$ via one of the following operations: in-
splitting, out-splitting, in-amalgamation at an interior vertex, or out-amalgamation at an interior vertex.

By Theorem \ref{thm:normal-form}, there exists a sequence of meteor graphs of length three $Z_1,Z_2,\ldots ,Z_t$ such that  \[
 S_k =Z_1 \to Z_2 \to \cdots \to Z_t
\] where $Z_t$ is in normal form and each $Z_i \to Z_{i+1}$ denotes either an in-splitting or an out-splitting.
Therefore, without loss of generality,  we may assume that $S_k$ is a meteor graph of length three in normal form.  
From this observation, and since $\approx$ is an equivalence relation on the set of meteor graphs with length three in normal form, it suffices to show that $E \approx F$ whenever $F = S_k$ in $(1)$. Hence, by Proposition \ref{maintheo-firstcase}, it follows that $ E \approx F.$ 

Assume inductively that $H\approx H'$ whenever $H$ and $H'$ are meteor graphs of length three in normal form, $A_H\approx_{SSE} A_{H'}$, and the transform degree $\mathrm{tr.deg} (H,H')$ is less than $n$, where $n\geq 2$. 

Let $E$ and $F$ be meteor graphs of length three in normal form such that $A_E\approx_{SSE} A_F$ and $\mathrm{tr. deg}(E,F) = n$. The construction below yields graphs $E'$, $G$, and $F'$ in normal form such that 
$$\mathrm{tr.deg}(E', G)=1,\quad  \mathrm{tr.deg}(G,F') \le n-1,$$
and moreover 
$$N_{F'}-N_{E'} = N_F - N_E,\quad \text{ and }\quad N_i^{E}=N_i^{E'}=N_i^{F'}=N_i^F \text{ for } i=1,2.$$
Then it follows from the induction hypothesis that $E'\approx G$ and $G\approx F'$, so that $E'\approx F'$. Then the above identities immediately imply that $E\approx F$. 

Consider a transformation from $E$ to $F$ of the following form:
\[
\scalebox{0.75}{
$\xymatrixrowsep{0.4pc}\xymatrix{E\ar[r]&S_1\ar[r]&\cdots\ar[r]&S_{i_1}&S_{i_1+1}\ar[l]&\ar[l]\cdots&\ar[l]S_j\ar[r]&S_{j+1}\ar[r]&\cdots\ar[r]&S_{i_2}&S_{i_2+1}\ar[l]&\ar[l]\cdots&\ar[l]S_k\ar[ld]\\
&&&&&&&&&&&\ar@{..}[llllddddd]&\\
&&&&&&&&&&&&\\
&&&&&&&&&&&&\\
&&&&&&&&&&&&\\ 
&&&&&&&&&&&&&\\
&&&&&&&&&&&&&\\
&&&&&&\ar[ru]S_t\ar[r]&S_{t+1}\ar[r]&\cdots\ar[r]&S_{i_n}&S_{i_n+1}\ar[l]&\ar[l]\cdots&\ar[l]S_m=F,}$
}
\] 
where $A \to B$ denotes that $B$ is obtained from $A$ via one of the following operations: in-
splitting, out-splitting, in-amalgamation at an interior vertex, or out-amalgamation at an interior vertex. 
By Theorem \ref{thm:quasi-normal-form}, there exists a sequence of meteor graphs of length three, $Y_0, Y_1, Y_2, \ldots, Y_h$, such that  \[
 S_j= Y_0 \to Y_1 \to Y_2 \to \cdots \to Y_h,\] where all trails in $\mathrm{Trail}(Y_h)$ have no common interior vertices, and each transformation $Y_i \to Y_{i+1}$ is either an in-splitting or an out-splitting at an interior vertex. Consequently, there is a reverse transformation from $Y_h$ to $S_j$ given by
\[Y_h \to \cdots \to Y_2 \to Y_1 \to Y_0=S_j,\] where each step $Y_{i+1} \to Y_i$ is either an in-amalgamation or an out-amalgamation at an interior vertex.

We denote by $X_E$, $X_F$, $X_{S_i}$ ($1\le i\le m$), and $X_{Y_l}$ ($1\le l\le h$) the sets of all exits for the first cycle in $E$, $F$, $S_i$, and $Y_l$, respectively. Let $c$ be a positive integer that is a multiple of $xyz$ and satisfies $$c\ge \max\{xy + 2x + 2y,\; yz + 2y + 2z,\; zx +2 z + 2x\}.$$ Let 
\begin{center}
 $E^{(1)} = E_{(X_E, c, 1)}$, 
$F^{(1)} = E_{(X_F, c, 1)}$, $S^{(1)}_i = E_{(X_{S_i}, c, 1)}$ ($1\le i\le m$),   
\end{center}
 and $$Y^{(1)}_l = E_{(X_{Y_l}, c, 1)}\ (1\le l\le h)$$ denote the first-type edge extensions of $E$, $F$, $S_i$, and $Y_j$, respectively. Since $E$ and $F$ are in normal form, it follows that $E^{(1)}$ and $F^{(1)}$ are also in normal form. By Corollary \ref{cor:first-type-edge-ext}, we obtain the following transformation of meteor graphs of length three: 
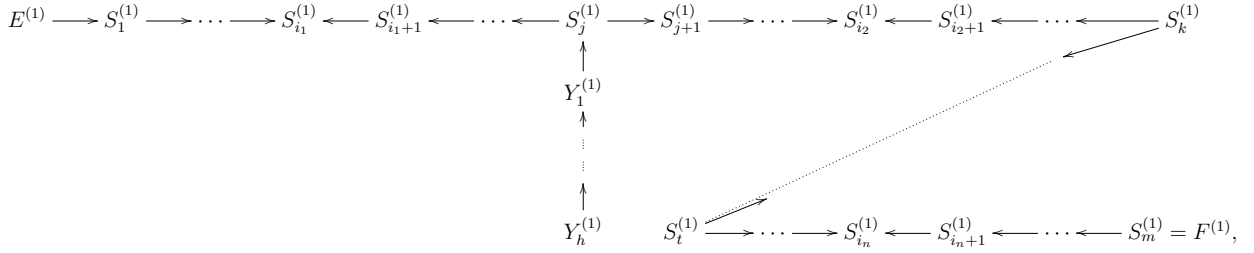
\[
\scalebox{0.7}{
$\xymatrixrowsep{0.4pc}\xymatrix{E^{(1)}\ar[r]&S^{(1)}_1\ar[r]&\cdots\ar[r]&S^{(1)}_{i_1}&S^{(1)}_{i_1+1}\ar[l]&\ar[l]\cdots&\ar[l]S^{(1)}_j\ar[r]&S^{(1)}_{j+1}\ar[r]&\cdots\ar[r]&S^{(1)}_{i_2}&S^{(1)}_{i_2+1}\ar[l]&\ar[l]\cdots&\ar[l]S^{(1)}_k\ar[ld]\\
&&&&&&&&&&&\ar@{.}[lllldddddd]&\\
&&&&&&\ar[uu]Y^{(1)}_1&&&&&&\\
&&&&&&\ar[u]&&&&&&\\
&&&&&&\ar@{..}[u]&&&&&&\\ 
&&&&&&\ar@{..}[u]&&&&&&&\\
&&&&&&&&&&&&&\\
&&&&&&\ar[uu]Y^{(1)}_h&\ar[ru]S^{(1)}_t\ar[r]&\cdots\ar[r]&S^{(1)}_{i_n}&S^{(1)}_{i_n+1}\ar[l]&\ar[l]\cdots&\ar[l]S^{(1)}_m=F^{(1)},}$
}
\] where each $A \to B$ denotes that $B$ is obtained from $A$ via one of the following operations: in-
splitting, out-splitting, in-amalgamation at an interior vertex, or out-amalgamation at an interior vertex.

We denote by $Z_{E^{(1)}}$, $Z_{F^{(1)}}$, $Z_{S^{(1)}_i}$ ($1\le i\le m$), and $Z_{Y^{(1)}_l}$ ($1\le l\le h$) the sets of all entrances for the third cycle in $E^{(1)}$, $F^{(1)}$, $S^{(1)}_i$, and $Y^{(1)}_l$, respectively. Let $d$ be a positive integer that is a multiple of $xyz$ and satisfies $$d\ge \max\{2xy + 2x + 2y,\; 2yz + 2y + 2z,\;2 zx +2 z + 2x\}.$$ 
Let 
\begin{center}
$E^{(2)} = E^{(1)}_{(Z_{E^{(1)}}, d, 2)}$, 
$F^{(2)} = F^{(1)}_{(Z_{F^{(1)}}, d, 2)}$, $S^{(2)}_i = (S^{(1)}_i)_{(Z_{S^{(1)}_i}, d, 2)}$ ($1\le i\le m$), 
\end{center} and $$Y^{(2)}_l = (Y^{(1)}_l)_{(Z_{Y^{(1)}_l}, d, 2)}\ (1\le l\le h)$$ denote the second-type edge extensions of $E^{(1)}$, $F^{(1)}$, $S^{(1)}_i$, and $Y^{(1)}_j$, respectively. Since $E^{(1)}$ and $F^{(1)}$ are in normal form, it follows that $E^{(2)}$ and $F^{(2)}$ are also in normal form. By Lemma \ref{lem:second-type-edge-ext}, we obtain the following transformation of meteor graphs of length three: 
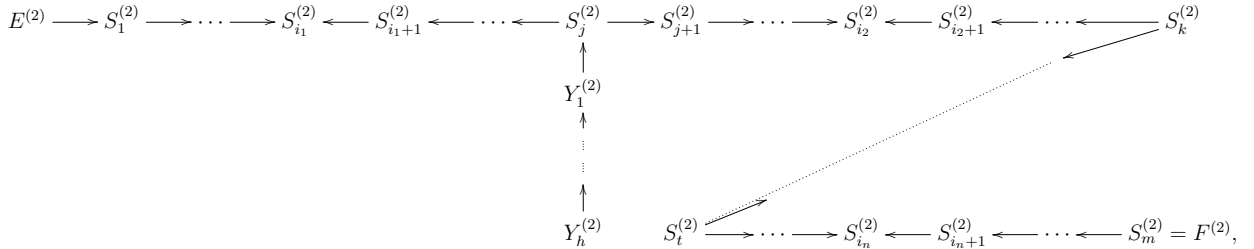
\[
\scalebox{0.7}{
$\xymatrixrowsep{0.4pc}\xymatrix{E^{(2)}\ar[r]&S^{(2)}_1\ar[r]&\cdots\ar[r]&S^{(2)}_{i_1}&S^{(2)}_{i_1+1}\ar[l]&\ar[l]\cdots&\ar[l]S^{(2)}_j\ar[r]&S^{(2)}_{j+1}\ar[r]&\cdots\ar[r]&S^{(2)}_{i_2}&S^{(2)}_{i_2+1}\ar[l]&\ar[l]\cdots&\ar[l]S^{(2)}_k\ar[ld]\\
&&&&&&&&&&&\ar@{..}[lllldddddd]&\\
&&&&&&\ar[uu]Y^{(2)}_1&&&&&&\\
&&&&&&\ar[u]&&&&&&\\
&&&&&&\ar@{..}[u]&&&&&&\\ 
&&&&&&\ar@{..}[u]&&&&&&&\\
&&&&&&&&&&&&&\\
&&&&&&\ar[uu]Y^{(2)}_h&\ar[ru]S^{(2)}_t\ar[r]&\cdots\ar[r]&S^{(2)}_{i_n}&S^{(2)}_{i_n+1}\ar[l]&\ar[l]\cdots&\ar[l]S^{(2)}_m=F^{(2)},}$
}
\] where each $A \to B$ denotes that $B$ is obtained from $A$ via one of the following operations: 
in-splitting, out-splitting, in-amalgamation at an interior vertex, or out-amalgamation at an interior vertex.

Since $Y^{(2)}_h$ satisfies the hypothesis of Lemma \ref{lm:normal-cover}, there exist a meteor graph $G$ of length three in normal form and a meteor graph $M$ of length three such that the following conditions hold:

$(1)$ $M$ is obtained from $G$ by a sequence of in-splittings and out-splittings;

$(2)$ $Y^{(2)}_h$ is a subgraph of $M$;
    
$(3)$ $\mathrm{Trail}(M) \setminus \mathrm{Trail}(Y^{(2)}_h) \subset \mathrm{Trail}_{M}(P_{M}, R_{M});$  

$(4)$ Any two distinct trails in $\mathrm{Trail}(M)$ have no common interior vertices;

$(5)$ All trails in $\mathrm{Trail}(M) \setminus \mathrm{Trail}(Y^{(2)}_h)$ start and end in a vertex of index $1$; 

$(6)$ For all trails $\alpha$ and $\beta \in  \mathrm{Trail}(M) \setminus \mathrm{Trail}(Y^{(2)}_h)$, if $$|\alpha| \equiv |\beta| \pmod{(x,z)},$$ then $|\alpha| = |\beta|$;

$(7)$ For a trail $\alpha \in \mathrm{Trail} (M)\setminus \mathrm{Trail}(Y^{(2)}_h)$ and any $1\le i \le (x,z)$, if $$|\alpha| \equiv i \pmod{(x,z)},$$ then $$|\alpha| \equiv i \pmod{xz}.$$

By condition (1), there exists a sequence of meteor graphs of length three, $G_1, G_2, \ldots, G_t$, such that  \[
 G\to G_t \to  \cdots \to G_2\to G_1\to  M,\] where  each transformation $A \to B$ is either an in-splitting or an out-splitting.
 
 Let $S:=  \mathrm{Trail}(M) \setminus \mathrm{Trail}(Y^{(2)}_h)$. Define 
\begin{center}
$E^{(3)} := (E^{(2)})_{(S, 3)}$, $F^{(3)} := (F^{(2)})_{(S, 3)}$, $S^{(3)}_i := (S^{(2)}_i)_{(S, 3)}$ ($1\le i\le m$),   
\end{center}
and $$Y^{(3)}_l := (Y^{(2)}_l)_{(S, 3)}\ (1\le l\le h)$$ to be the third-type edge extensions of $E^{(2)}$, $F^{(2)}$, $S^{(2)}_i$, and $Y^{(2)}_j$, respectively. It is obvious that $M = Y^{(3)}_h$. From these observations and by Lemma \ref{lm:third-edge-ext}, we obtain the following transformation of meteor graphs of length three: 
\[
\scalebox{0.67}{
$\xymatrixrowsep{0.4pc}\xymatrix{E^{(3)}\ar[r]&S^{(3)}_1\ar[r]&\cdots\ar[r]&S^{(3)}_{i_1}&S^{(3)}_{i_1+1}\ar[l]&\ar[l]\cdots&\ar[l]S^{(3)}_j\ar[r]&S^{(3)}_{j+1}\ar[r]&\cdots\ar[r]&S^{(3)}_{i_2}&S^{(3)}_{i_2+1}\ar[l]&\ar[l]\cdots&\ar[l]S^{(3)}_k\ar[ld]\\
&&&&&&&&&&&\ar@{..}[lllldddddd]&\\
&&&&&&\ar[uu]Y^{(3)}_1&&&&&&\\
&&&&&&\ar[u]&&&&&&\\
&&&&&&\ar@{..}[u]&&&&&&\\ 
&&&&&&\ar@{..}[u]&&&&&&&\\
&&&&&&&&&&&&&\\
&G\ar[r]&G_s\ar[r]&\cdots\ar[r]&G_2\ar[r]&G_1\ar[r]&M=\ar[uu]Y^{(3)}_h&\ar[ru]S^{(3)}_t\ar[r]&\cdots\ar[r]&S^{(3)}_{i_n}&S^{(3)}_{i_n+1}\ar[l]&\ar[l]\cdots&\ar[l]S^{(3)}_m=F^{(3)},}$
}
\] where  $A \to B$ denotes that $B$ is obtained from $A$ via one of the following operations: in-
splitting, out-splitting, in-amalgamation at an interior vertex, or out-amalgamation at an interior vertex. 

Since  both $E^{(3)}$ and $F^{(3)}$ are in quasi-normal form and satisfy the hypotheses of Corollary \ref{cor:normal-cover}, there exist two meteor graphs $E'$ and $F'$ of length three in normal form that can be transformed into $E^{(3)}$ and $F^{(3)}$,
respectively, via a finite sequence of in-splittings and out-splittings.
We then have 
\begin{center}
 $\text{tr.deg}(E', G) = 1$ and $\text{tr.deg}(G, F') \le \text{tr.deg}(E, F) -1 = n -1$.   
\end{center}
By the induction hypothesis, we obtain that $E'\approx G$ and $G\approx F'$, and hence $$E'\approx F',$$ since $\approx$
is an equivalence relation. By the definition of third-type edge extension graphs $E^{(3)}$ and $F^{(3)}$, it follows that 
\begin{center}
$N_{F^{(3)}}- N_{E^{(3)}} = N_{F^{(2)}}- N_{E^{(2)}}$, 
\end{center}
and $$N^{E^{(3)}}_i = N^{E^{(2)}}_i \text{ and }
N^{F^{(3)}}_i = N^{F^{(2)}}_i$$ for $i = 1, 2$. By the proof of Corollary \ref{cor:normal-cover}, taking into account that both $E^{(3)}$ and $F^{(3)}$ are in quasi-normal form, we have 
\begin{center}
$N_{E'} = N_{E^{(3)}}$ and $N^{E'}_i = N^{E^{(3)}}_i$ for $i= 1, 2$,  
\end{center}
and 
\begin{center}
$N_{F'} = N_{F^{(3)}}$ and $N^{F'}_i = N^{F^{(3)}}_i$ for $i= 1, 2$.
\end{center}
From these equalities, and since $E'\approx F'$, we get
$$N_i^{E^{(2)}} =N_i^{E^{(3)}} = N_i^{E'}=N_i^{F'} = N_i^{F^{(3)}}=N_i^{F^{(2)}}$$
for $i=1,2$, and moreover
$$N_{F^{(2)}} -N_{E^{(2)}} =  N_{F^{(3)}} -N_{E^{(3)}} = N_{F'} -N_{E'}.$$
Using again that $E'\approx F'$, these identities immediately imply that 
$$E^{(2)}\approx F^{(2)}.$$
Moreover, by the definitions of the first-edge and second-edge extensions of meteor graphs of length three, we always have
\begin{center}
$E \approx E^{(1)}\approx E^{(2)}$ and $F \approx F^{(1)}\approx F^{(2)}$.    
\end{center}
Since $\approx$ is an equivalence relation, it follows that $E \approx F$, thus finishing the proof.
\end{proof}

\section{Application: Williams' Conjecture}\label{sec5}
In this section, based on Theorems \ref{thm:normal-form} and \ref{numtheo}, we prove that both Williams' Conjecture and Hazrat's Graded Morita Equivalence Conjecture hold for meteor graphs of length three whose cycle lengths are pairwise coprime (Theorem \ref{thm21}). Combining this result with \cite[Theorem 4.3]{dohaznam}, we immediately obtain that both conjectures hold for graphs with disjoint cycles containing exactly three cycles with pairwise coprime lengths (Corollary \ref{cor22}). \medskip

Along this section, we will work with meteor graphs of length three.
We will use the notation introduced in Definition \ref{def31} for a meteor graph of length three in normal form. We also will denote $|P_E|=x$, $|Q_E|=y$ and $|R_E|=z$. Note in particular that all trails from $P_E$ to $Q_E$ start at $p_1^E$ and end at $q_1^E$, all trails from $Q_E$ to $R_E$ start at $q_1^E$ and end at $r_1^E$, and all trails from $P_E$ to $R_E$ start at $p_1^E$ and end at $r_1^E$.

We begin this section with the following useful lemma.

\begin{lemma}\label{lem}
Let $E$ be a meteor graph of length three in normal form with a unique chain of cycles $P_E > Q_E >   R_E$. Then, for any element $t \in T_E$, there exist integers $j_1$ and $j_2 \in \mathbb{Z}$ such that $t$ can be uniquely written as
\[
\sum_{0 \leq i<x} a_i p_1^E(j_1-i)+\sum_{0 \leq i<y} b_i q_1^E(j_2 -i) + \sum \limits_{0 \leq i < z}^{} c_i r_1^E(i).
\]
where $a_i, b_i,c_i \in \mathbb{N}$.
\end{lemma}
\begin{proof}
To simplify the notation, we will skip the superindex $E$ from $p_1^E,q_1^E,r_1^E$ in the course of this proof.

We can see that for any interior vertex $z_1$ in a trail from $P_E$ to $Q_E$, there exist a natural number $j_1$ such that $z_1 = q_1(j_1)$, for any interior vertex $z_2$ in a trail from $Q_E$ to $R_E$, there exists a natural number $j_2$ such that $z_2= r_1(j_2)$, and for any  interior vertex $z_3$ in a trail from $P_E$ to $R_E$, there exists a natural number $j_3$ such that $z_3= r_1(j_3)$.
Hence, we can write $t$ in the form: \begin{equation}
 t= \sum \limits_{v \in (P_E)^0}^{} \left( \sum \limits_{a \in \mathbb{Z}}^{}e_{v,a} v(a) \right)+ \sum \limits_{w \in (Q_E)^0}^{} \left( \sum \limits_{b \in \mathbb{Z}}^{}f_{w,b} w(b) \right) + \sum \limits_{u \in (R_F)^0}^{} \left( \sum \limits_{c \in \mathbb{Z}}^{}g_{u,c} u(c) \right) 
\end{equation}
where $e_{v,a} , f_{w,b}, g_{u,c} \in \mathbb{N}$ and $e_{v,a} = f_{w,b} = g_{z,c}=0$ for cofinitely many indices $a,b$ and $c$, respectively. 
Let $S_E$ be the submonoid of $T_E$ generated by the set $$\{w(b),u(c) \mid w \in (Q_E)^0, u\in (R_E)^0 , b,c \in \mathbb{Z}\}.$$ For any $v \in (P_E)^0$ we have  
\begin{center}
 $v(a) =p_1(a+l_{v}) $ for all $a \in \mathbb{Z}$
\end{center}
where $l_{v}$ is the length of the path in $P_E$ from $v$ to $p_1$, since $E$ is in normal form (note that $l_{p_1}=0$). Moreover  
$$p_1 = p_1(x) + \sum \limits_{i}e_i q_1(t_i) + \sum \limits_{i} f_i r_1(l_i),$$
for $e_i,f_i,t_i,l_i \in \mathbb{N}$. Let $M_1$ be an integer satisfying \[
    M_1\ge \max \{ a+ l_{v} + x \mid  v \in (P_E)^0 \mbox{ and } a \in \mathbb{Z} \mbox{ with } e_{v,a} \ne 0\} 
.\]
Let $v \in (P_E)^0$ with $e_{v,a} \ne 0$, and let $k$ be the unique positive integer such that 
\begin{center}
$ a+l_{v} + k x \le M_1\,\,$ and $\,\, a+l_{v} + (k+1) x > M_1$. 
\end{center}
We then have $$v(a) = p_1(a+ l_{v} + kx) + \phi (v) $$ for some $\phi (v) \in S_E.$ From these observations and equation (2), we obtain that 
\begin{equation}
 t = \sum \limits_{0 \le i < x }^{} a_i p_1(M_1 - i ) +  t_1
\end{equation}
for some $t_1 \in S_E$. 
Let $U_E$ be the submonoid of $T_E$ generated by the set 
$$\{u(c) \mid u \in (R_E)^0 , c \in \mathbb{Z}\}.$$ 
For any $w \in Q_E$, we have  
\begin{center}
 $w(b) =q_1(b+l_{w})$ for all $b  \in \mathbb{Z}$,   
\end{center}
and 
$$q_1 = q_1(y) + \sum \limits_{i} f_i r_1(l_i),$$
where $l_{w}$ is the length of the path in $Q_E$ from $w$ to $q_1$, and $f_i,l_i \in \mathbb{N}$. Let $M_2$ be an integer satisfying \[
    M_2\ge \max \{ b+ l_{w} + y \mid  w \in (Q_E)^0 \mbox{ and } b \in \mathbb{Z} \mbox{ with } f_{w,b} \ne 0\} 
.\]
Let $w$ be a vertex on $Q_E$ with $f_{w,b} \ne 0$, and let $k$ be the unique positive integer such that 
\begin{center}
$a+l_{w} + k y \le M_2\,\, $ and $\,\,  b+l_{w} + (k+1) y > M_2$.    
\end{center}
We then obtain that $$w(b) = q_1(b+ l_{y} + ky) + \phi_2 (u)$$ for some $\phi_2(u) \in U_E.$ Moreover, for every  $u \in (R_E)^0$ and $c \in \mathbb{Z}$, we have $$u( c) = u'(c+1),$$ where  $u'$ is the range of the edge in $R_E$ with source $u$.
From these observations and equation (3), we obtain that 
\[t = \sum \limits_{0 \le i < x }^{} a_i p_1(M_1 - i ) + \sum_{0 \leq i<y} b_i q_1(M_2 -i) + t_2\] 
for some $t_2 \in U_E$. Note that $t_2$ can be written as $$t_2 =\sum \limits_{c \in \mathbb{Z}}^{}g_{r,c} r_1(c).$$ We also have $r_1(i) = r_1(i+z)$ for all $i\in \mathbb{Z}$.   It follows that  \[
 t = \sum \limits_{0 \le i < x}^{} a_i p_1(M_1 - i ) + \sum_{0 \leq i<y} b_i q_1(M_2 -i) + \sum \limits_{0 \leq i < z}^{} c_i r_1(i) 
.\]
We now show that this representation is unique. Suppose that it can also be written as \[
t = \sum \limits_{0 \le i < x}^{} a_i' p_1(M_1 - i ) + \sum_{0 \leq i<y} b_i' q_1(M_2 -i) + \sum \limits_{0 \leq i < z}^{} c_i' r_1(i)
.\]
Consider the $\mathbb{Z}$-order ideal  $\langle q_1,r_1 \rangle$ of  $T_E$ generated by $q_1$ and $r_1$. Clearly, $\langle q_1,r_1 \rangle = S_E$. By \cite[Lemma 3.5]{alfi}, the natural inclusion $T_{Q_E \Rightarrow R_E}\hookrightarrow T_E $ descends to an isomorphism $$T_E /\left\langle q_1,r_1\right\rangle \cong T_{P_E},$$ where $Q_E \Rightarrow R_E$ denotes the subgraph of $E$ obtained by deleting the cycle $P_E$ together with all trails originating from $P_E$, and $T_{P_E}$ is the talented monoid of the cycle $P_E$. Now notice that we have a $\mathbb Z$-isomorphism $T_{P_E} \cong \mathbb{N}^{x}$. Passing to this quotient, the image of two presentations of $t$  satisfy \[
\sum \limits_{0 \le i < x}^{} a_i p_1(M_1 - i ) =  \sum \limits_{0 \le i < x }^{} a_i' p_1(M_1 - i ). 
\]
Hence, we have $a_i = a_i'$ for all $0 \le i <x$. It follows, using the cancellativity of $T_E$, that \[
 t'=\sum_{0 \leq i<y} b_i q_1(M_2 -i) + \sum \limits_{0 \leq i < z}^{} c_i r_1(i)= \sum_{0 \leq i<y} b_i' q_1(M_2 -i) + \sum \limits_{0 \leq i < z}^{} c_i' r_1(i) 
.\]
Notice that we may regard $t'$  as an element in  the talented monoid $T_{Q_E \Rightarrow R_E}$. Consider the $\mathbb{Z}$-order ideal  $\langle r_1 \rangle$ of  $S_E$ generated by $r_1$. It is clear that $\langle r_1 \rangle = U_E$. By \cite[Lemma 3.5]{alfi}, the natural inclusion $T_{R_E}\hookrightarrow T_{Q_E \Rightarrow R_E} $ induces an isomorphism $$T_{Q_E \Rightarrow R_E} /\left\langle r_1\right\rangle \cong T_{Q_E},$$ 
where $T_{Q_E}$ is the talented monoid of $Q_E$. Just as before, we first obtain 
$b_i = b_i'$ for all $0 \le i <y$, and, after cancelling the term $\sum_{0 \leq i<y} b_i q_1(M_2 -i)$, that
\[t_1  = \sum \limits_{0 \leq i < z}^{} c_i r_1(i)=  \sum \limits_{0 \leq i < z}^{} c_i' r_1(i) 
.\]
Viewing $t_1$ as an element in the monoid $T_{R_E} \cong \mathbb{N}^{z}$, we conclude that $c_i = c_i'$ for all $0 \le i < z,$ thus finishing the proof.
\end{proof}

Consequently, we obtain the following useful corollaries. 

\begin{cor}\label{coro1-lm51}
Let $E$ be a meteor graph of length three in normal form  with a unique chain of cycles $P_E > Q_E >  R_E$. 
Let $t$ be an element of $T_E$ of the form
$$t = p_1^E(a)+ \sum_{j=1}^{y} \ell_{j} q^{E}_1(s + j)  + \sum_{j=1}^{z} t_{j} r_1^E(j),$$ where $s \in \mathbb{Z}$, $\ell_{j}, t_j \in \mathbb{N}$, and $a$ is a negative integer divisible by $x$. Then $t$ can be rewritten in the form:
\[t = p_1^E + \sum \limits_{i=1}^{y} a_i q_1^E(I+i) + \sum \limits_{i=1}^{z} b_i r_1^E(i),\] where
$I\in \mathbb{Z}$ and $a_i, b_i\in \mathbb{N}$.
\end{cor}
\begin{proof}
It immediately follows from the algorithm established in the proof of Lemma \ref{lem}  by taking $M_1 = 0$.
\end{proof}

Following \cite[Subsection 2.2, p. 325]{alfi}, a nonzero element $\alpha$ in a commutative monoid $M$ is called an {\it atom} if $\alpha = y +z$ then $y= 0$ or $z =0$.

The following corollary describes all atoms of a meteor graph of length three in normal form.

\begin{corollary}\label{coro4.2}
Let $E$ be a meteor graph of length three in normal form  with a unique chain of cycles $P_E > Q_E >   R_E$. Then, a nonzero element $\beta \in T_E$ is an atom in $T_E$ if and only if  there exists an integer $k$ such that $ 0 \le k < z $ and  $\beta= r_1^E(k)$. \end{corollary}
\begin{proof}  
($\Longrightarrow$). Assume that $\beta$ is an atom in $T_E$. By Lemma \ref{lem}, $\beta$ can be uniquely written in the form: 
\[
\beta=   \sum_{0 \leq i<x} a_i p_1^E(j_1-i)+\sum_{0 \leq i<y} b_i q_1^E(j_2 -i) + \sum \limits_{0 \leq i < z}^{} c_i r_1^E(i), \] 	where $a_{i}, b_{i},c_i \in \mathbb{N}$. Since $\beta$ is an atom, we get that only one coefficient in $a_i,b_i,c_i$ is $1$, all other coefficients will be $0$.  From the proof of Lemma \ref{lem}, we obtain that 
\[ p_1^E(i)= p_1^E(i+x)+ s,\] where $s$ is a non-zero element  in $T_E$, and \[ q_1^E(j)= q_1^E(j+y) +t, \] where $t$ is a non-zero element in $T_E$. This implies that  $p_1^E(i)$ and $q_1^E(j)$ cannot be atoms for all $i,j \in \mathbb{Z}$, and so $\beta=r_1^E(k)$ for some $0 \le k <z$.

$(\Longleftarrow)$. It immediately follows from Lemma \ref{lem}, thus finishing the proof.
\end{proof}

Theorem \ref{numtheo} enables us to reduce a meteor graph of length three with pairwise coprime cycle lengths to a simpler form.

\begin{lemma}\label{lm:redu-SSE}
Let $E$ be a meteor graph of length three in normal form with a unique chain of cycles $P_E > Q_E > R_E$ such that  $|P_E|$, $|Q_E|$ and $|R_E|$ are pairwise coprime. Let $E_{rd}$ be the meteor graph of length three in normal form with a unique chain of cycles $P_{E_{rd}} > Q_{E_{rd}} > R_{E_{rd}}$ satisfying the following conditions:

$(1)$ $|P_E| = |P_{E_{rd}}|$, $|Q_E| = |Q_{E_{rd}}|$, and $|R_E| = |R_{E_{rd}}|$;

$(2)$ $N_1^E(1) = N_1^{E_{rd}}(1)$, $N_2^E(1) = N_2^{E_{rd}}(1)$, and $N_E(1) = N_{E_{rd}}(1)$;

$(3)$ Every trail from $P_{E_{rd}}$ to $Q_{E_{rd}}$ has length $xy$;

$(4)$ Every trail from $Q_{E_{rd}}$ to $R_{E_{rd}}$ has length $yz$;

$(5)$ Every trail from $P_{E_{rd}}$ to $R_{E_{rd}}$ has length $xz$.\\
Then $E$ is strongly shift equivalent to $E_{rd}$; that is, their adjacency matrices are strongly shift equivalent.
\end{lemma}
\begin{proof}
It is easy to see that $E\approx E_{rd}$. Then, by Theorem \ref{numtheo}, it follows that $E$ is strongly shift equivalent to $E_{rd}$, thus finishing the proof.
\end{proof}

We refer to the graph $E_{rd}$, introduced in Lemma \ref{lm:redu-SSE}, as the {\it reduction graph} associated with a meteor graph $E$ of length three whose cycle lengths are pairwise coprime.

We are in a position to state the main result of this section, showing that Williams' conjecture and Hazrat's Graded Morita Equivalence Conjecture hold for meteor graphs of length three.

\begin{theorem}\label{thm21}
\setcounter{equation}{0}
Let $E$ and $F$ be essential graphs, where $E$ is a meteor graph of length three with a unique chain of cycles $P_E > Q_E > R_E$ such that  $|P_E|$, $|Q_E|$ and $|R_E|$ are pairwise coprime, and let $A_E$ and $A_F$ be the adjacency matrices of $E$ and $F$, respectively.  Let $K$ be an arbitrary field. Then the following are equivalent:
	
$(1)$ The Leavitt path algebras $L_K(E)$ and $L_K(F)$ are graded Morita equivalent;
 
$(2)$ There is an order-preserving $\mathbb{Z}\left[x, x^{-1}\right]$-module isomorphism $K_{0}^{\mathrm{gr}}(L_K(E)) \rightarrow K_{0}^{\mathrm{gr}}(L_K(F))$;

$(3)$ The talented monoids $T_{E}$ and $T_{F}$ are $\mathbb{Z}$-isomorphic;

$(4)$ $A_E\sim_{SE} A_F$;

$(5)$ $A_E\sim_{SSE} A_F$. 	
\end{theorem}
\begin{proof}    
$(1) \Longrightarrow (2)$. By \cite[Theorem 2.3.7]{hazbk}, the graded Morita equivalence gives rise to invertible bimodules, which in turn induce an isomorphism on the level of graded $K^{\gr}_{0}$. 

$(2) \Longleftrightarrow (3)$. By \cite{hazli}, the positive cone of the graded Grothendieck group $K_{0}^{\operatorname{gr}}(L_K(E))$ is $\mathcal{V}^{\operatorname{gr}}(L_K(E))$ and $\mathcal{V}^{\operatorname{gr}}(L_K(E)) \cong T_{E}$ as $\mathbb{Z}$-monoids. This implies the equivalence.

$(3) \Longleftrightarrow (4)$. It immediately follows from Theorem \ref{h99} and $(2) \Longleftrightarrow (3)$. 

$(5) \Longrightarrow (1)$. It immediately follows from \cite[Proposition 15]{hazd}.

$(4) \Longrightarrow (5)$.
 Assume that $A_E \sim _{SE} A_F$.  By Proposition \ref{propreserve}, we obtain that $F$ is a meteor graph of length three with a unique chain of cycles $P_F \Rightarrow Q_F  \Rightarrow R_F$ such that   $$x:=|P_E|=|P_F|,\ y:=|Q_E|= |Q_F|,\text{ and } z:=|R_E|= |R_F|.$$
Since $E$ and $F$ are shift equivalent, we have that $T_E$ and $T_F$ are $ \mathbb{Z}$-isomorphic to each other. Assume that $\varphi:T_E \to T_F$ is a $ \mathbb{Z}$-isomorphism. By Theorems \ref{willimove} and \ref{thm:normal-form}, we may assume that $E$ and $F$ are meteor graphs of length three in normal form. By Lemma \ref{lm:redu-SSE}, we may also assume, without loss of generality, that they are reduction graphs.

We claim that
\begin{equation}
\varphi\left(p_1^E\right)={ }^{a} p_1^F+ \sum_{j=1}^{y} \ell_{j}\left({}^{s+j} q^{F}_1\right)  + \sum_{j=1}^{z} t_{j}\left({}^{j} r_1^F\right) 
\end{equation}
and
\begin{equation}
\varphi\left(q_1^E\right)={ }^{b} q^{F}_1 + \sum_{j=1}^{z} v_{j}\left({}^{j} r_1^F\right)  
\end{equation}
for some $a, b, s \in \mathbb{Z}$ and $\ell_{j}, t_j ,v_j\in \mathbb{N}$. Indeed, by Lemma \ref{lem}, there exist integers $j_1$ and $j_2 \in \mathbb{Z}$ such that $$\varphi\left(q_1^E\right) = \sum \limits_{0 \le i <x }a'_{i} p_1^F(j_1 - i)+\sum \limits_{0 \leq i<y} b'_{i} q_1^F(j_2 -i ) + \sum \limits_{i=1}^{z} c'_i r_1^F(i),$$ where $a'_i, b'_i, c'_i \in \mathbb{N}$. By Corollary \ref{coro4.2}, we obtain that  $r_1^E$ and $r_1^F$ are atoms  in $T_E$ and $T_F$, respectively, and so we have \begin{equation}
  \varphi (r_1^E)\ =\ ^tr_1^F.  
\end{equation}
From (3), $\varphi$ induces an isomorphism $$\bar{\varphi}: T_{E} /\left\langle r_1^E\right\rangle\rightarrow T_{F} /\left\langle  r_1^F\right\rangle,$$ passing to the quotient. Hence, we obtain that $$\bar{\varphi}\left(q_1^E\right)=\sum \limits_{0 \le i <x }a'_{i} p_1^F(j_1 - i)+  \sum \limits_{0 \leq i<y} b'_{i} q_1^F(j_2 -i ).$$ Since $q_1^E$ is an atom in $T_{E} /\left\langle r_{1}^{E}\right\rangle $ (by \cite[Lemma 4.2]{dohaznam}), it follows that $\bar{\varphi}\left(q_1^E\right)$ is also an atom in $T_{F} /\left\langle r_{1}^{F}\right\rangle$. Again by \cite[Lemma 4.2]{dohaznam}, we deduce that $a'_i=0$ for all $0 \le i <x $, and $$b'_0+b'_1+\cdots +b'_{y-1}=1,$$ which shows the identity $(2)$.  By Lemma \ref{lem}, we obtain that $$\varphi\left(p_1^E\right) = \sum \limits_{0 \le i <x }a''_{i} p_1^F(j_1 - i)+\sum \limits_{0 \leq i<y} b''_{i} q^F_1(j_2 -i) + \sum \limits_{i=1}^{z} c''_i r^F_1(i),$$
where $a''_i, b''_i, c''_i\in \mathbb{N}$.
From $(3)$ and $(2)$, $\varphi$ induces an isomorphism $$\bar{\varphi}: T_{E} /\left\langle q_1^E, r_{1}^{E}\right\rangle\rightarrow T_{F} /\left\langle q_1^F, r_{1}^{F}\right\rangle,$$ passing to the quotient. Consequently, we have  $$\bar{\varphi}\left(p_1^E\right)=\sum\limits_{i=0}^{x-1} a''_{i} p_1^F\left(j_{1}-i\right).$$ Since $p_1^E$ is an atom in $T_{E} /\left\langle q_1^E, r_{1}^{E}\right\rangle $, it follows that $$a''_0+a''_1+\cdots +a''_{x-1}=1,$$ which establishes the identity $(1).$

We note that, for any fixed $k\in \mathbb{Z}$, the map $T_{F} \longrightarrow T_{F}$ given by $x \longmapsto$$ ^{- k}x$ is a $\mathbb{Z}$-isomorphism. Then, since $x$, $y$ and $z$ are pairwise coprime, we may choose an integer $k$ satisfying 
\begin{center}
 $k \equiv a \pmod{x},$\quad $k \equiv b \pmod{y},$\quad $ k \equiv t \pmod{z},$   
\end{center}
 and such that
$a - k <0$, $b-k<0$ and $t-k <0$. Replacing the right-hand side of formula $(1)$ by its image under this isomorphism
and using $r^F_1(i) = r^F_1(i + z)$ for all $i \in \mathbb{Z}$, we have
$$\varphi(p_1^E) = p_1^F(a-k) + \sum \limits_{i=1}^{y} 
\ell_i q_1^F(s-k+i) + \sum \limits_{i=1}^{z} t'_i r_1^F(i),$$ where $t'_i\in \mathbb{N}$. Since $a - k$ is a negative multiple of $x$, Corollary \ref{coro1-lm51} implies that
\begin{align}
    \varphi(p_1^E) & = p_1^F + \sum \limits_{i=1}^{y} a_i q_1^F(I+i) + \sum \limits_{i=1}^{z} b_i r_1^F(i), 
\end{align} where $a_i,b_i\in \mathbb{N}$ and $I \in \mathbb{Z}$. Replacing the right-hand side of formula $(2)$ by its image under the above isomorphism, we have $$\varphi\left(q_1^E\right)= q^{F}_1(b-k) + \sum_{j=1}^{z} v_{j} r_1^F(j-k).$$ Since $b - k$ is a negative multiple of $y$, we obtain as before \begin{align}
\varphi(q_1^E) & = q_1^F + \sum \limits_{i=1}^{z} c_i r_1^F(i),
\end{align} where $c_i\in \mathbb{N}$. 
Applying the above isomorphism to the right-hand side of formula $(3)$ and using the relation $r^F_1(i) = r^F_1(i + z)$ for all $i \in \mathbb{Z}$, we obtain that
\begin{align}
\varphi(r_1^E) & = r_1^F.
\end{align}
Using that $E$ and $F$ are reduced, it is straightforward to verify that
\begin{align}
p_1^E& = p_1^E(x) +  N_1^E q_1^E(xy) + N^E r_1^E(xz), \\
q_1^E& = q_1^E(y) + N_2^E r_1^E(yz), \\ 
p_1^F& = p_1^F(x) +   N_1^F q_1^F(xy) +  N^F r_1^F(xz), \\
q_1^F & = q_1^F(y) +  N_2^F r_1^F(yz), 
\end{align}
where $N_1^E$ denotes the number of trails from $P_E$ to $Q_E$, $N_2^E$ denotes the number of trails from $Q_E$ to $R_E$, $N_1^F$ denotes the number of trails from $P_F$ to $Q_F$, $N_2^F$ denotes the number of trails from $Q_F$ to $R_F$, $N^E$ denotes the number of trails from $P_E$ to $R_E$, and $N^F$ denotes the number of trails from $P_F$ to $R_F$. 

By substituting $(8)$ into the left-hand side of $(5)$, we obtain that
\begin{align*}
            \varphi(q_1^E)  & = \varphi(q_1^E(y) + N_2^E r_1^E(yz) )  \tag{using (8)}\\
                       & = \varphi(q_1^E(y) ) + N_{2}^E \varphi(r_1^E(yz)) \\
                       & = q_1^F(y) + \sum \limits_{i=1}^{z}  c_ir_1^F(y+i) +  N_2^E r_1^F(yz). \tag{using (5) and (6) }
        \end{align*}
On the other hand, by applying $(10)$ to the right-hand side of (5), we have \begin{align*}
              q_1^F + \sum \limits_{i=1}^{z} c_i   r_1^F(i)  = q_1^F(y)  +  N_2^F r_1^F(yz) + \sum \limits_{i=1}^{z} c_i   r_1^F(i).
         \end{align*}
It follows that \[
   q_1^F(y) + \sum \limits_{i=1}^{z}  c_ir_1^F(y+i) +  N_2^E r_1^F(yz)=  q_1^F(y)  +  N_2^F r_1^F(yz) + \sum \limits_{i=1}^{z} c_i   r_1^F(i),\]   
and so \begin{equation}
      \sum \limits_{i=1}^{z}  c_ir_1^F(y+i) +  N_2^E r_1^F(yz)=   N_2^F r_1^F(yz) + \sum \limits_{i=1}^{z} c_i   r_1^F(i)
    .\end{equation}
Since $r_1^F(i) = r_1^F(i + z)$ for all $i$, equation $(11)$ may be rewritten in the form:
\[\sum^z_{i=1}t_i r^F_1(i) = \sum^z_{i=1}t'_i r^F_1(i).\] By Lemma \ref{lem}, it follows that $t_i = t'_i$ for all $1\le i\le z$. We have
\[\sum \limits_{1 \le j \le z}  t_j= \sum \limits_{1 \le j \le z,  }  c_j + N_2^E \]
and 
\[\sum \limits_{1 \le j \le z }  t'_j= N_2^F + \sum \limits_{1 \le j \le z } c_j.\]
Therefore, we obtain that
\[
      \sum \limits_{1 \le j \le z,  }  c_j + N_2^E = N_2^F + \sum \limits_{1 \le j \le z } c_j,
      \]
which immediately implies that
\begin{center}
 $N_2^E = N_2^F:=N_2$.    
\end{center}
By the same argument, we also have 
\begin{center}
$N_{1}^E = N_{1}^F:=N_1 $.    
\end{center}
Then, using $(11)$, we obtain that 
\begin{equation}   
\sum \limits_{i=1}^{z} c_{i}r_1^F(i+y) =    \sum \limits_{i=1}^{z} c_i   r_1^F(i).  
\end{equation}
For each  $i \in \mathbb{Z}$, by comparing the coefficients of $r_1^F(i)$ on both sides, we have \[ c_{i-y}  = c_i,\] where  the indices of the sequence $(c_j)$ are taken modulo $z$.  We also have $c_{i-z} = c_i$ for all $i \in \mathbb{Z}$, and so  $c_{i-1} = c_{i-(y,z)} = c_i$ for all $i$. It follows that $$c_1=c_2= \cdots = c_z :=C.$$ 
We next analyze the left-hand side of $(4)$ and receive that \begin{align*}
          \varphi(p_1^E)  =&  \varphi(p_1^E(x)  + N_1 q_1^E(xy) + N^E  r_1^E(xz) ) \tag{using (7) } \\
          =& \varphi(p_1^E(x)) +  N_1 \varphi(q_1^E(xy)) + N^E  \varphi(r_1^E(xz))  \\ 
          =&  p_1^F(x) + \sum \limits_{i=1}^{y} a_i q_1^F(I+x+i)  + \sum \limits_{i=1}^{z} b_i r_1^F(x+i) +  N_1 \left( q_1^F(xy) + \sum \limits_{i=1}^{z} Cr_1^F(xy+i) \right) \\
           & +  N^E r_1^F(xz).  \tag{using (4),(5) and (6)}
      \end{align*}
Now we analyze the right-hand side of $(4)$. By applying equation (9), we obtain that 
\begin{align*} 
        p_1^F + \sum \limits_{i=1}^{y} a_i q_1^F(I+i) + \sum \limits_{i=1}^{z} b_i r_1^F(i) = & p_1^F(x) +  N_1 q_1^F(xy) +  N^F r_1^F(xz)  \\                                                                          & + \sum \limits_{i=1}^{y} a_i q_1^F(I+i) + \sum \limits_{i=1}^{z} b_i r_1^F(i),
    \end{align*}
which yields 
\begin{align*}
p_1^F(x) + &\sum \limits_{i=1}^{y} a_i q_1^F(I+x+i)  + \sum \limits_{i=1}^{z} b_i r_1^F(x+i) +  N_1 \left( q_1^F(xy) + \sum \limits_{i=1}^{z} C r_1^F(xy+i) \right) 
          \\ &   +  N^E r_1^F(xz)  =   p_1^F(x) +  N_1 q_1^F(xy) +  N^F r_1^F(xz)                                                                         + \sum \limits_{i=1}^{y} a_i q_1^F(I+i) + \sum \limits_{i=1}^{z} b_i r_1^F(i).
  \end{align*}
Hence,\begin{align*}
 \sum \limits_{i=1}^{y} a_i q_1^F(I+x+i) & + \sum \limits_{i=1}^{z} b_i r_1^F(x+i) +  N_1 \left( \sum \limits_{i=1}^{z} C r_1^F(xy+i) \right) 
           +  N^E r_1^F(xz)  \\ &=    N^F r_1^F(xz)                                                                         + \sum \limits_{i=1}^{y} a_i q_1^F(I+i) + \sum \limits_{i=1}^{z} b_i r_1^F(i)
  \end{align*}
 Let
 \[
     X_1:=   N_1   \left( \sum \limits_{j=1}^{z} C r_1^F(xy+i-I) \right)   +  N^E r_1^F(xz-I) + \sum \limits_{i=1}^{z} b_i r_1^F(x+i-I)
 \]
 and
 \[X_2:= N^F  r_1^F(xz-I) + \sum \limits_{i=1}^{z} b_i r_1^F(i-I).\] 
Applying the automorphism $x\mapsto {}^{-I}x$, we then have 
 \begin{equation}
  \sum \limits_{i=1}^{y} a_i q_1^F(x+i)     + X_1  =    \sum \limits_{i=1}^{y} a_iq_1^F(i) + X_2 .
  \end{equation} 
For each  $i \in \mathbb{Z}$, by comparing the coefficients of $q_1^F(j)$ where $j\equiv i \pmod{y}$ on both sides, we obtain that \begin{equation}
    a_{i-x}  =  a_i,
\end{equation}
where the indices of the sequence $(a_j)$ are taken modulo $y$. We note that $a_{i-y} = a_i$ for all $i$, and hence,  $a_{i-1}= a_{i-(x,y)}= a_i$. This implies that $$a_1=a_2= \cdots = a_y := A.$$
Equation \((8)\) can be generalized as
\begin{equation}
q_1^F(l) = q_1^F(l+y) +  N_2 r_1^F(l)
\end{equation}
for all $l \in \mathbb{Z}$. Let $m_y(n)$ and $u_y(n)$ be the remainder and quotient, respectively, when $n$ is divided by $y$. Define
\[T(i)=
\begin{cases}
N_2 \left( \sum \limits_{j=0}^{u_y(i+x)-1} r_1^F(m_y(i+x) +jy ) \right),
& \text{if } u_y(i+x) \ge 1,\\[2mm]
0, & \text{ if } u_y(i+x) =0,
\end{cases}
\] for all $i\in \mathbb{Z}$. We claim that \[q_1^F(m_y(i+x) )  = q^F_1(i+x) + T(i)  \tag{16}\]
for all $1 \le i \le y.$ Indeed, using formula $(15)$, we get that
\begin{align*}
q_1^F(l)  =& q_1^F(l +y)  + N_2 r_1^F(l) = q_1^F(l +2y) + N_2 r_1^F(l + y) + N_2 r_1^F(l)\\
          =& q_1^F(l +2y) + N_2(r_1^F(l) + r_1^F(l + y))\\ 
          = & \cdots\\
          = &  q_1^F(l + ky) + N_2(r_1^F(l) + \cdots + r_1^F(l + (k-1)y))      
\end{align*}
for all positive integers $k$ and for all integers $l$. We also observe that
\[ i + x = u_y(i + x)y + m_y(i + x).\]
It follows immediately from these observations that $(16)$ holds, as desired.

By adding $A (T(1) + T(2) + \cdots +  T(y) ) $ to both sides of $(13)$, we obtain that \begin{align*}
   \sum \limits_{i=1}^{y} A \left(  q_1^F(x+i) + T(i) \right)      +X_1  =    \sum \limits_{i=1}^{y} A q_1^F(i) + A \left( \sum \limits_{i=1}^y  T(i) \right)  +X_2 .
\end{align*}
By formula $(16)$, it follows that \[
\sum \limits_{i=1}^y A q_1^F (m_y (i+x))  + X_1 = A \left( \sum \limits_{i=1}^y q_1^F(i) \right) +  A \left( \sum \limits_{i=1}^y  T(i) \right) + X_2,
\]
equivalently \[\sum \limits_{i=0}^{y-1} A q_1^F (i)  + X_1 = A \left( \sum \limits_{i=1}^y q_1^F(i) \right) +  A \left( \sum \limits_{i=1}^y  T(i) \right) + X_2,\]
and so \[ A q_1^F + X_1 = A q_1^F(y) + A\left( \sum \limits_{i=1}^y  T(i) \right) + X_2.\]
Using formula $(8)$, we have \[
Aq_1^F(y) +A N_2 r_1^F(yz) + X_1 = A q_1^F +X_1=  Aq_1^F(y) + A \left( \sum \limits_{i=1}^y T(i) \right) +X_2
,\]
and hence \[AN_2 r_1^F(yz) + X_1 =    A \left( \sum \limits_{i=1}^y T(i) \right) + X_2.\]
Since $r_1^F(i) = r_1^F(i + z)$ for all $i$, the above equation may be rewritten in the form:
\[\sum^z_{i=1}g_i r^F_1(i) = \sum^z_{i=1}g'_i r^F_1(i).\] By Lemma \ref{lem}, it follows that $g_i = g'_i$ for all $1\le i\le z$. We have
\[\sum \limits_{1 \le j \le z}  g_j= AN_2 + N_1 z C + N^E + \sum \limits_{i=1}^z b_i\]
and 
\[\sum \limits_{1 \le j \le z }  g'_j= N^F + A\left( \sum \limits_{i=1}^y u_y(i+x)  \right) N_2 + \sum \limits_{i=1}^z b_i.\] Hence, we obtain that 
\[ AN_2 + N_1 z C + N^E + \sum \limits_{i=1}^z b_i = N^F + A\left( \sum \limits_{i=1}^y u_y(i+x)  \right) N_2 + \sum \limits_{i=1}^z b_i. \tag{17}\]
We also note that $u_y(i+x) = \floor*{\dfrac{x+i}{y} }$. We write $x=ky+r$, where $0 \le r \le y-1$. We then have  $$u_y(x+1)= u_y(x+2)= \cdots = u_y(x+y-r-1) = k$$ and $$u_y(x+y-r)=\cdots = u_y(x+y) = k+1,$$ and so $$\sum \limits_{i=1}^y u_y(i+x) =k(y-r-1) + (k+1)(r+1)= ky +r +1 =x +1.$$  Applying this to $(17),$ we obtain that \[N^E - N^F = (Ax)N_2 - (Cz)N_1= \frac{1}{y} (-C(yz)N_1+A(xy)N_2).\]
This shows that $E\approx F$. 
Now it follows from Theorem \ref{numtheo} that $A_E \approx _{SSE} A_F$, thus finishing the proof.
\end{proof}

We end this section with the following corollary, showing that Williams' conjecture and Hazrat's Graded Morita Equivalence Conjecture hold for graphs of disjoint cycles which contain three cycles of coprime lengths.

\begin{cor}\label{cor22}
Let $E$ and $F$ be essential graphs, where $E$ is a graph with disjoint cycles which contains exactly three cycles whose lengths are pairwise coprime, and let $A_E$ and $A_F$ be the adjacency matrices of $E$ and $F$, respectively.  Let $K$ be an arbitrary field. Then the following conditions are equivalent:
	
$(1)$ The Leavitt path algebras $L_K(E)$ and $L_K(F)$ are graded Morita equivalent;
	
$(2)$ There is an order-preserving $\mathbb{Z}\left[x, x^{-1}\right]$-module isomorphism $K_{0}^{\mathrm{gr}}(L_K(E)) \rightarrow K_{0}^{\mathrm{gr}}(L_K(F))$;
	
$(3)$ The talented monoids $T_{E}$ and $T_{F}$ are $\mathbb{Z}$-isomorphic;

$(4)$ $A_E\sim_{SE} A_F$;
	
$(5)$ $A_E\sim_{SSE} A_F$. 	
\end{cor}
\begin{proof}
If $E$ contains a chain of cycles of length three, then the statement immediately follows from Theorem \ref{thm21}. Otherwise, the statement  follows from \cite[Theorem 4.3]{dohaznam}, thus finishing the proof.
\end{proof}

\section{Acknowledgements}
The first author was partially supported by the Spanish State Research Agency (grant No.\  PID2023-147110NB-I00). The second and third authors were supported by the Vietnam Academy of Science and Technology under grant CBCLCA.01/26-28. The authors would like to thank Roozbeh Hazrat and Adam Dor-On for their valuable comments and suggestions, which helped shape the final version of this article.

\end{document}